\documentclass[12pt]{article}
\usepackage{amsthm,amsfonts,amssymb,amscd}

\begin{document}

%\begin{flushleft}
%ÓÄÊ 512.76
%\end{flushleft}

\begin{center}
\large \bf Birational geometry of varieties,\\
fibred into complete intersections of codimension two
\end{center}\vspace{0.5cm}

\centerline{A.V.Pukhlikov}\vspace{0.5cm}

\parshape=1
3cm 10cm \noindent {\small \quad\quad\quad \quad\quad\quad\quad
\quad\quad\quad {\bf }\newline In this paper we prove the
birational superrigidity of Fano-Mori fibre spaces $\pi\colon V\to
S$, every fibre of which is a complete intersection of type
$d_1\cdot d_2$ in the projective space ${\mathbb P}^{d_1+d_2}$,
satisfying certain conditions of general position, under the
assumption that the fibration $V\slash S$ is sufficiently twisted
over the base (in particular, under the assumption that the
$K$-condition holds). The condition of general position for every
fibre guarantees that the global log canonical threshold is equal
to one. This condition bounds the dimension of the base $S$ by a
constant that depends on the dimension $M$ of the fibre only (as
the dimension $M$ of the fibre grows, this constant grows as
$\frac12 M^2$). The fibres and the variety $V$ itself may have
quadratic and bi-quadratic singularities, the rank of which is
bounded from below.

Bibliography: 37 items.}

AMS classification: 14E05, 14E07

Key words: Fano variety, Mori fibre space, birational map,
birational rigidity, linear system, maximal singularity, quadratic
singularity, bi-quadratic singularity.\vspace{0.3cm}

\section*{Introduction}

{\bf 0.1. The main results.} Let us fix a pair of integers
$(d_1,d_2)$, such that $d_2\geqslant d_1\geqslant 2$ and
$d_2\geqslant 27$. Set $M=d_1+d_2-2$ and denote by the symbol
${\mathbb P}$ the complex projective space ${\mathbb P}^{M+2}$
with homogeneous coordinates $(x_0:\dots:x_{M+2}$). The aim of
this paper is to prove birational rigidity of Fano-Mori fibre
spaces $\pi\colon V\to S$, every fibre of which is a complete
intersection of type $d_1\cdot d_2$ in ${\mathbb P}$, satisfying
certain natural conditions of general position, under the
assumption that the fibration $V\slash S$ is sufficiently twisted
over the base. When the base $S$ is fixed, a majority of families
of fibre spaces satisfies the condition of twistedness. The
condition of general position for every fibre bounds the dimension
of the base $S$ by a constant that depends on $M$ only (as the
dimension $M$ of the fibre grows, that constant grows as $\frac12
M^2$). The key property, which must be satisfied for every fibre
of the fibration $V\slash S$ in order for our proof to work, is
the boundedness, in a certain sense, of the singularities of every
divisor. Now let us give precise statements.

Let $S$ be a non-singular projective rationally connected variety
of positive dimen\-sion. By a {\it Fano-Mori fibre space} over $S$
we mean a surjective morphism $\pi\colon V\to S$, every fibre of
which is irreducible, reduced and of dimension $\mathop{\rm
dim}V-\mathop{\rm dim}S\geqslant 3$, where the variety $V$ is
projective, factorial, with at most terminal singularities, and
moreover, for the Picard numbers the equality $\rho(V)=\rho(S)+1$
holds and the anti-canonical class $(-K_V)$ is ample on a general
(and thus on every) fibre. (We assume the base of the fibre space
to be rationally connected and its total space to be factorial
from the beginning because the total spaces of all fibre spaces
considered in this paper are rationally connected and factorial.)

Let $\pi'\colon V'\to S'$ be a {\it rationally connected fibre
space}, that is, a surjective morphism of non-singular projective
varieties, the general fibre of which is an irreducible rationally
connected variety and the base $S'$ is rationally connected. We
say that a birational map (if there are such maps, in particular,
$\mathop{\rm dim}V=\mathop{\rm dim} V'$)
$$
\chi\colon V\dashrightarrow V'
$$
is {\it fibre-wise}, if there is a rational dominant map
$\beta\colon S\dashrightarrow S'$, such that
$\pi'\circ\chi=\beta\circ\pi$, that is, the following diagram of
maps commutes:
$$
\begin{array}{rcccl}
   & V & \stackrel{\chi}{\dashrightarrow} & V' & \\
\pi\!\!\!\!\! & \downarrow &   &   \downarrow & \!\!\!\!\!\pi' \\
   & S & \stackrel{\beta}{\dashrightarrow} & S'.
\end{array}
$$
We emphasize that the property of being fibre-wise does not mean
that $\beta$ is birational: $\chi$ maps the fibres a fibre of the
fibration $\pi$ onto, generally speaking, a closed subset of a
fibre of the fibration $\pi'$, so that the inverse map $\chi^{-1}$
may well not be fibre-wise (more precisely, both maps $\chi$,
$\chi^{-1}$ are fibre-wise if and only if the map $\beta$ is
birational).

Now let $\pi'\colon V'\to S'$ be a {\it Mori fibre space}, that
is, the singularities of $V'$ and $S'$ are terminal and ${\mathbb
Q}$-factorial, the anticanonical class of the variety $V'$ is
relatively ample and the equality $\rho(V')=\rho(S')+1$ holds.

{\bf Definition 0.1.} A Fano-Mori fibre space $V\slash S$ is {\it
birationally rigid}, is for every birational map $\chi\colon
V\dashrightarrow V'$ onto the total space of (any) Mori fibre
space $V'\slash S'$ there is a birational map $\beta\colon
S\dashrightarrow S'$, such that $\pi'\circ\chi=\beta\circ\pi$,
that is, $\chi$ maps the fibre of general position of the
projection $\pi$ birationally onto the fibre of general position
of the projection $\pi'$. If, moreover, the restriction of $\chi$
onto the general fibre of the projection $\pi$ is always a
biregular isomorphism, then the Fano-Mori fibre space $V\slash S$
is {\it birationally superrigid}.

For $d\in{\mathbb Z}_+$ let ${\cal P}_{d,M+3}$ be the linear space
of homogeneous polynomials of degree $d$ in the variables
$x_0,\dots,x_{M+2}$ and
$$
{\cal F}={\cal P}_{d_1,M+3}\times{\cal P}_{d_2,M+3}
$$
--- the corresponding space of pairs $(f_1,f_2)$. We say that
the set $\{f_1=f_2=0\}\subset{\mathbb P}$ of their common
zeros is a complete intersection of codimension 2 with good
singularities, if the scheme-theoretic intersection of
hypersurfaces $\{f_1=0\}$ and $\{f_2=0\}$ is an irreducible
reduced subvariety $F=F(f_1,f_2)\subset{\mathbb P}$ of codimension
2, and moreover, for every point $o\in F$ one of the following
three cases takes place:

(1) the point $o\in F$ is non-singular,

(2) the point $o$ is non-singular on one of the hypersurfaces
$\{f_1=0\}$, $\{f_2=0\}$, and the variety $F$ has at the point $o$
a quadratic singularity of rank $\geqslant 5$,

(3) the point $o$ is a quadratic singularity on both hypersurfaces
$\{f_1=0\}$ and $\{f_2=0\}$, and if
$$
\begin{array}{ccccccccc}
f_1 & = & f_{1,2} & + & \dots & + & f_{1,d_1}, & & \\
f_2 & = & f_{2,2} & + & \dots &   &            & + & f_{2,d_2}, \\
\end{array}
$$
are decompositions of the non-homogeneous presentations of the
polynomials $f_1$, $f_2$ with respect to some system of affine
coordinates $z_1,\dots, z_{M+2}$ with the origin at the point $o$,
into homogeneous components, $\mathop{\rm deg\,} f_{i,j}=j$, then
the quadratic forms $f_{1,2}$ and $f_{2,2}$ satisfy the following
condition of general position. In order to make its formulation
more convenient, let us give the following definition, which will
be useful later.

{\bf Definition 0.2.} The {\it rank} of a set of quadratic forms
$q_1,\dots, q_k$ in a certain set of variables is the number
$$
\mathop{\rm rk} (q_1,\dots, q_k)=\mathop{\rm min}\{ \mathop{\rm
rk} (\lambda_1 q_1+\dots+\lambda_k q_k)\,|\, (\lambda_1,\dots,
\lambda_k)\neq (0,\dots, 0)\}
$$
(the minimum is taken over all tuples $(\lambda_1,\dots, \lambda_k)\in
{\mathbb C}^k\setminus \{0\}$).

Now the condition of general position for the forms $f_{1,2}$ and
$f_{2,2}$ takes the following form:
$$
\mathop{\rm rk} (f_{1,2}, f_{2,2})\geqslant 7.
$$

The points of type (3) we call bi-quadratic singularities of $F$.
Quadratic and bi-quadratic singularities will be considered in
detail in \S 1.

If $F$ is a complete intersection of codimension 2 with good
singularities, then the inequality
$$
\mathop{\rm codim}(\mathop{\rm Sing}F\subset F)\geqslant 4
$$
holds, so that by Grothendieck's theorem \cite{CL} the variety $F$
is a factorial complete intersection, $\mathop{\rm Pic}F={\mathbb
Z}H$, where $H$ is the class of a hyperplane section, and
$K_F=-H$, that is to say, $F$ is a primitive Fano variety. It is
easy to check (see Subsection 1.7), that the singularities of the
variety $F$ are terminal. The open set of pairs $(f_1,f_2)\in{\cal
F}$, such that the set of common zeros $\{f_1=f_2=0\}$ is a
complete intersection of codimension 2 with good singularities
will be denoted by the symbol ${\cal F}_{\rm bq}$.

Now let $S$ be a non-singular projective rationally connected
variety of positive dimension and $\pi_X\colon X\to S$ a locally
trivial fibre space with the fibre ${\mathbb P}$. A subvariety
$V\subset X$ of codimension 2 is a {\it fibration into complete
intersections of type} $d_1\cdot d_2$, if the base $S$ can be
covered by Zariski open subsets $U$, over which the fibration
$\pi_X$ is trivial, $\pi^{-1}_X(U)\cong U\times{\mathbb P}$, and
for each of them there is a map
$$
\Phi_U\colon U\to{\cal F}_{\rm bq},
$$
such that for every point $s\in U$ the subvariety of common zeros
$\{\Phi_U(s)=0\}\subset{\mathbb P}$ is $V\cap\pi^{-1}_X(s)$ in the
sense of the trivialization above. In other words, over $U$ the
variety $V$ is given by the pair of equations
$$
f_1(x_*,s)=0,\,\,f_2(x_*,s)=0,
$$
the coefficients of which are regular functions on $U$, and for
every $s\in U$
$$
f_1(x_*,s),f_2(x_*,s)\in{\cal F}_{\rm bq}.
$$
The restriction of the projection $\pi_X$ onto $V$ we denote by
the symbol $\pi$. As we will show in Subsection 1.7, the
singularities of the variety $V$ are terminal, the variety $V$
itself is factorial, and
$$
\mathop{\rm Pic}V={\mathbb Z}K_V\oplus\pi^*\mathop{\rm Pic}S,
$$
that is, $\pi\colon V\to S$ is a Fano-Mori fibre space over $S$.

In order to state the first main result of the present paper,
recall one more well known definition.

{\bf Definition 0.3.} A primitive Fano variety $F$ (that is to
say, a factorial projective variety with terminal singularities,
the Picard group of which is generated by the canonical class) is
{\it divisorially canonical}, if for every $n\geqslant 1$ and
every effective divisor $D_F\in|-nK_F|$, the pair
$(F,\frac{1}{n}D_F)$ is canonical, that is, for every exceptional
divisor $E$ over $F$ the inequality
$$
\mathop{\rm ord}\nolimits_ED_F\leqslant n\cdot a(E,F)
$$
holds. Now we are ready to state the first main result, a
sufficient condition for birational superrigidity of the Fano-Mori
fibre space $V\slash S$ constructed above. We somewhat abuse the
notations and identify the pair $(f_1,f_2)$ and the corresponding
subvariety $\{f_1=f_2=0\}$, so that we write, for example,
$F\in{\cal F}_{\rm bq}$ for a fibre of the fibre space $\pi$.

{\bf Theorem 0.1.} {\it Assume that the Fano-Mori fibre space
$\pi\colon V\to S$, constructed above, satisfies the following
properties:

{\rm (i)} for every point $s\in S$ the corresponding fibre
$F_s=\pi^{-1}(s)\in{\cal F}_{\rm bq}$ is a divisorially canonical variety,

{\rm (ii)} for every effective divisor
$$
D\in|-nK_V+\pi^*Y|,
$$
$n\geqslant 1$, the class $Y$ is pseudoeffective on $S$,

{\rm (iii)} for every mobile family of irreducible rational curves
$\overline{\cal C}$ on the base $S$, sweeping out a dense open
subset $S$, and a curve $\overline{C}\in\overline{\cal C}$ no
positive multiple of the class of the algebraic cycle
$$
-(K_V\cdot\pi^{-1}(\overline{C}))-F
$$
(where $F$ is the class of a fibre of the projection $\pi$) is
effective, that is to say, it is not rationally equivalent to an
effective cycle of dimension $\mathop{\rm dim}F$ on $V$.

Then for every rationally connected fibre space $V'\slash S'$
every birational map $\chi\colon V\dashrightarrow V'$ is
fibre-wise, and the fibre space $V\slash S$ itself is birationally
superrigid.}

We emphasize that the two claims of Theorem 0.1 are independent of
each other: the fact that every birational map $\chi$ onto the
total space of every rationally connected fibre space is
birational {\it does not imply} birational superrigidity (and
birational rigidity), and, similarly, birational suprrigidity does
not imply the first claim about the map being fibre-wise.
Therefore, Theorem 0.1 is an essential improve\-ment of Theorem
1.1 in \cite{Pukh15a} (as we will see below, in \S 1, the proof of
Theorem 0.1 works without any changes for Fano-Mori fibre spaces,
considered in \cite{Pukh15a}, and improves the results of that
paper accordingly). The assumptions of Theorem 0.1 are somewhat
different from those of Theorem 1.1 in \cite{Pukh15a}; how they
are related, we discuss in \S 1. Note for the moment, that the
assumption (ii) in Theorem 0.1 is the well known $K$-{\it
condition} for the fibration $V\slash S$: the anticanonical class
$(-K_V)$ is not contained in the interior of the pseudoeffective
cone:
$$
(-K_V)\not\in\mathop{\rm Int} A^1_+V.
$$
Theorem 0.1 motivates the following general conjecture.

{\bf Conjecture 0.1.} {\it Assume that an arbitrary Fano-Mori
fibre space satisfies the $K$-condition and the fibre of general
position is a birationally rigid (respectively, birationally
superrigid) Fano variety. Then this fibre space is birationally
rigid (respe\-ctively, birationally superrigid) and every
birational map from its total space onto the total space of an
arbitrary rationally connected fibre space (if there are such
birational maps) is fibre-wise. In particular, on the total space
of the original fibre space there are no structures of a
rationally connected fibre space over a base of dimension higher
than that of the base of the original fibre space.}

Apart from Theorem 0.1, the present paper contains two more main
results, the meaning of which is that they make it possible to
construct very large classes of Fano-Mori fibre spaces, satisfying
the assumptions of Theorem 0.1. Precise statements of those
results is given below.\vspace{0.3cm}

%%%%%%%%%%%%%%%%%%%%%%%%%%%%%%%%%%%%%%%%%%%%%%%%%%%%%%%%%%%%%%%%%%
%%%%%%%%%%%%%%%%%%%%   subsection 0.2

{\bf 0.2. Conditions of general position.} Let $o\in{\mathbb P}$
be an arbitrary point and $(z_1,\dots,z_{M+2}$) a system of affine
coordinates with the origin at that point. Abusing the notations,
we use the same symbols $f_1,f_2$ for the corresponding
non-homogeneous polynomials in $z_*$. Assume that $f_1,f_2$ vanish
at the point $o$ and write down
$$
f_i(z_*)=f_{i,1}+f_{i,2}+\dots+f_{i,d_i},
$$
where $i=1,2$, the polynomials $f_{i,j}$ are homogeneous of degree
$j$. For instance, the system of equations $f_1=f_2=0$ defines in
a neighborhood of the point $o$ a non-singular complete
intersection of codimension 2 if and only if the linear forms
$f_{1,1},f_{2,1}$ are linearly independent. Let us put the pairs
of indices $(i,j)$ in the lexicographic order:
$(i_1,j_1)<(i_2,j_2)$, if $j_1<j_2$ or $j_1=j_2$, but $i_1<i_2$.
The polynomials $f_{i,j}|_{\{f_{1,1}=f_{2,1}=0\}}$ with
$j\geqslant 2$, placed in that order, form the sequence
$$
f_{1,2}|_{\{f_{1,1}=f_{2,1}=0\}},\quad\dots,\quad
f_{2,d_2}|_{\{f_{1,1}=f_{2,1}=0\}}
$$
which we denote by the symbol ${\cal S}$. Removing the last
$k\geqslant 1$ polynomials, we obtain the sequence ${\cal S}[-k]$.
Note that in the case of a quadratic singularity the linear
subspace $\{f_{1,1}=f_{2,1}=0\}$ has codimension 1, and in the
case of a bi-quadratic singularity it is the whole space ${\mathbb
C}^{M+2}_{z_1,\dots,z_{M+2}}$.

Now let us state the conditions of general position (the
regularity conditions), which will be needed in order to prove the
condition (i) of Theorem 0.1 for a sufficiently large class of
complete intersections of codimension 2. For a linear subspace
${\cal L}\subset\{f_{1,1}=f_{2,1}=0\}$ the restriction of the
sequence ${\cal S}[-k]$ onto ${\cal L}$ (that is, the sequence,
consisting of the restrictions of all polynomials in ${\cal
S}[-k]$ onto ${\cal L}$) is denoted by the symbol ${\cal
S}[-k]|_{\cal L}$.

{\bf The regularity condition at a non-singular point.}

Assume that the linear forms $f_{1,1}$ and $f_{2,1}$ are linearly
independent. The point $o$ is {\it regular}, if the following
condition (R1) is satisfied.

(R1) For any linear subspace
$$
{\cal L}\subset\{f_{1,1}=f_{2,1}=0\}
$$
of codimension 2 the sequence ${\cal S}[-5]|_{{\cal L}}$ is
regular, that is, the set of its common zeros in ${\mathbb
P}({\cal L})\cong{\mathbb P}^{M-3}$ has dimension 2.

{\bf The regularity condition at a quadratic point.}

Assume now that the linear forms $f_{1,1}$ and $f_{2,1}$ are
linearly dependent, but are not both identically zero:
$f_{1,1}=\alpha_1\tau$ and $f_{2,1}=\alpha_2\tau$, where
$\tau(z_*)$ is a non-zero linear form and
$(\alpha_1,\alpha_2)\neq(0,0)$. We say that the point $o$ is {\it
regular} in that case, if the following condition (R2) is
satisfied.

(R2) The quadratic form
$$
(\alpha_2f_{1,2}-\alpha_1f_{2,2})|_{\{\tau=0\}}
$$
is of rank $\geqslant 9$ and for any linear subspace ${\cal
L}\subset\{\tau=0\}$ of codimension 1 (with respect to the
hyperplane $\{\tau=0\}$) the sequence ${\cal S}[-4]|_{\cal L}$ is
regular, that is, the set of its common zeros in ${\mathbb
P}({\cal L})\cong{\mathbb P}^{M-2}$ is of dimension 3.

{\bf The regularity conditions at a bi-quadratic point.}

Assume now that both linear forms $f_{1,1}$ and $f_{2,1}$ vanish
identically. Let $\widetilde{\mathbb P}\to{\mathbb P}$ be the blow
up of the point $o$ and $E_{\mathbb P}$ the exceptional divisor.
The affine coordinates $(z_1,\dots,z_{M+2})$ generate homogeneous
coordinates
$$
(z_1:z_2:\dots:z_{M+2})
$$
on $E_{\mathbb P}\cong{\mathbb P}^{M+1}$. The point $o$ is {\it
regular}, if the following conditions (R$2^2$.1), (R$2^2$.2) and
(R$2^2$.3) are satisfied. The second and third conditions depend
on the value of the degree $d_1$, the first one is common for all
values of $d_1$.

(R$2^2$.1) The inequality
$$
\mathop{\rm rk} (f_{1,2},f_{2,2})\geqslant 13
$$
holds, and moreover, the rank of at least one of the two quadratic
forms $f_{1,1}$ and $f_{2,2}$ is at least 18.

We will see below (see Subsection 1.7), that the condition
(R$2^2$.1) implies that the system of equations
$$
f_{1,2}=f_{2,2}=0
$$
defines an irreducible reduced complete intersection $E\subset
E_{\mathbb P}$ of codimension 2, satisfying the inequality
$$
\mathop{\rm codim}(\mathop{\rm Sing}E\subset E)\geqslant 10.
$$

Now assume that $d_1\geqslant 4$.

(R$2^2$.2) For any linear subspace ${\cal L}\subset{\mathbb
C}^{M+2}$ of codimension 2 the system of equations
$$
f_{1,2}|_{\cal L}=f_{2,2}|_{\cal L}=f_1|_{\cal L}=f_2|_{\cal L}=0
$$
defines an irreducible reduced complete intersection of
codimension 4 in ${\cal L}\cong{\mathbb C}^M$.

(R$2^2$.3) For every linear subspace ${\cal L}\subset {\mathbb
C}^{M+2}$ of codimension $\mathop{\rm codim}{\cal L}\in \{2,3\}$
the sequence
$$
{\cal S}[-\mathop{\rm codim}{\cal L}-1]|_{\cal L}
$$
is regular.

Assume that $d_1=3$. In that case the second and third regularity
conditions take the following form.

(R$2^2$.2) For every linear subspace ${\cal L}\subset{\mathbb
C}^{M+2}$ of codimension 2 the system of equations
$$
f_{1,2}|_{\cal L}=f_{1,3}|_{\cal L}=f_{2,2}|_{\cal L}=
f_{2,3}|_{\cal L}=f_2|_{\cal L}=0
$$
defines an irreducible reduced complete intersection of
codimension 5 in ${\cal L}\cong{\mathbb C}^M$.

(R$2^2$.3) For every linear subspace ${\cal L}\subset {\mathbb
C}^{M+2}$ of codimension $\mathop{\rm codim}{\cal L}\in \{2,3\}$
the sequence
$$
{\cal S}[-\mathop{\rm codim}{\cal L}]|_{\cal L}
$$
is regular.

Assume, finally, that $d_1=2$. In that case the second and third
regularity conditions are as follows.

(R$2^2$.2) For every linear subspace ${\cal L}\subset{\mathbb
C}^{M+2}$ of codimension 2 the system of equations
$$
f_{1,2}|_{\cal L}=f_{2,2}|_{\cal L}=f_{2,3}|_{\cal L} =f_2|_{\cal
L}=0
$$
defines an irreducible reduced complete intersection of
codimension 4 in ${\cal L}\cong{\mathbb C}^M$, and the condition
(R$2^2$.3) is the same as in the case $d_1=3$.

We say that the pair of polynomials $(f_1,f_2)\in{\cal F}$ is {\it
regular}, if every point of the set of their common zeros is
regular in the sense of the corresponding condition (R1),(R2) or
(R$2^2$). The set of regular pairs will be denoted by the symbol
${\cal F}_{\rm reg}$. For $(f_1,f_2)\in{\cal F}_{\rm reg}$ the
corresponding complete intersection $F(f_1,f_2)$ is also said to
be regular. It is clear that a regular complete intersection is a
complete intersection of of codimension 2 with good singularities,
so that it is a primitive Fano variety.

The next claim is the second main result of the present paper.

{\bf Theorem 0.2.} {\it For $d_1\neq 3$ the complement ${\cal
F}\backslash{\cal F}_{\rm reg}$ is of codimension at least
$$
\frac12(M^2-17M+64)
$$
in ${\cal F}$. For $d_1=3$ the codimension of the complement
${\cal F}\backslash{\cal F}_{\rm reg}$ in ${\cal F}$ is at least}
$$
\frac12(M^2-19M+82).
$$

Theorem 0.2 implies that if for $d_1\neq 3$
$$
\mathop{\rm dim}S<\frac12(M^2-17M+64)
$$
(for $d_1=3$ the right hand side of the last inequality should be
replaced by $\frac12(M^2-19M+82)$), then for any ambient
projective bundle $\pi_X\colon X\to S$ a general subvariety
$V\subset X$ of codimension 2, which locally over $S$ is given by
the pair of equations
$$
f_1(x_*,s)=f_2(x_*,s)=0
$$
for $(f_1(s),f_2(s))\in{\cal F}$, is a Fano-Mori fibre space into
complete intersections of type $d_1\cdot d_2$, because we may
assume that
$$
(f_1(s),f_2(s))\in{\cal F}_{\rm reg}\subset{\cal F}_{\rm bq}.
$$
The following fact, which is the third main result of the present
paper, allows to apply Theorem 0.1 to the Fano-Mori fibre space
$\pi\colon V\to S$.

{\bf Theorem 0.3.} {\it A regular complete intersection $F\in{\cal
F}_{\rm reg}$ is divisorially canonical.}

The following example demonstrates that it is easy to check the
conditions (ii) and (iii) of Theorem 0.1.

{\bf Example 0.1.} Assume that $d_1\neq 3$ and $m$ is an integer,
satisfying the inequality
$$
1\leqslant m\leqslant\frac12(M^2-17M+62).
$$
Consider the direct product $X={\mathbb P}^m\times{\mathbb
P}^{M+2}$ and a general subvariety $V$ of codimension 2, which is
a complete intersection of two hypersurfaces of bidegree
$$
(l_1,d_1)\quad\mbox{and}\quad (l_2,d_2),
$$
respectively. Denote the projection of the variety $V$ onto
${\mathbb P}^m$ by the symbol $\pi$. Obviously, $\pi\colon
V\to{\mathbb P}^m$ is a Fano-Mori fibre space, the fibres of which
are complete intersections of type $d_1\cdot d_2$. Let $H_S$ and
$H_{\mathbb P}$ be the classes of hyperplanes in ${\mathbb P}^m$
and ${\mathbb P}$, respectively. The same symbols are used for
their pull backs on $X$ and subsequently, for their restrictions
onto $V$. Taking into account these notations, we have
$$
-K_V=(m+1-l_1-l_2)H_S+H_{\mathbb P}.
$$
It is easy to see that the conditions (ii) and (iii) of Theorem
0.1 are satisfied if the inequality
$$
((-K_V)\cdot\pi^{-1}(L)\cdot H^M_{\mathbb P})\leqslant 0
$$
holds, where $L\subset{\mathbb P}^m$ is a line. Simple
computations show that this inequality is equivalent to the
estimate
$$
l_1\left(1-\frac{1}{d_1}\right)+
l_2\left(1-\frac{1}{d_2}\right)\geqslant
m+1.
$$
If the integral parameters $l_1,l_2\in{\mathbb Z}_+$ satisfy the
last inequality, then the Fano-Mori fibre space $\pi\colon
V\to{\mathbb P}^m$ is birationally superrigid and $V$ has no other
structures of a Mori fibre space, apart from the projection $\pi$
(up to a fibre-wise birational equivalence). Note that in order to
get the condition (iii), the inequality
$$
l_1\left(1-\frac{1}{d_1}\right)+l_2\left(1-\frac{1}{d_2}\right)>
m
$$
is sufficient, and in order to get the condition (ii), the inequality
$l_1+l_2\geqslant m+1$ is sufficient, if
$$
H_{\mathbb P}\not\in\mathop{\rm Int}A^1_+V.
$$
On the other hand, for $l_1+l_2\leqslant m$ the anti-canonical
class $-K_V$ is ample and the projection onto ${\mathbb P}$ gives
a structure of a Fano-Mori fibre space on $V$, which is
``transversal'' to the original structure $\pi$, so that for
$l_1+l_2\leqslant m$ the fibration $\pi\colon V\to{\mathbb P}^m$
is not birationally rigid and our numerical conditions for
$l_1,l_2$, that provide birational superrigidity, are close to
optimal ones.

{\bf Remark 0.1.} The assumption that $d_2\geqslant 27$, which
bounds from below the dimension of complete intersections
considered in this paper (since $M=d_1+d_2-2$), is needed for the
only reason that for smaller dimensions the estimate in Theorem
0.2 is given by certain particular values, and not by a uniform
simple formula, which has the form of a quadratic polynomial in
$M$. What is essential, is the assumption on the rank of one of
the quadratic forms $f_{1,2},f_{2,2}$ in the condition (R$2^2$.1):
if it is removed, we have to exclude from consideration the
bi-quadratic points, which bounds from above the admissible
dimension of the base $S$ much stronger. For these reasons, in
order to have simple formulations of the main results, we assume
that $d_2\geqslant 27$, so that for the dimension we also have
$M\geqslant 27$.\vspace{0.3cm}

%%%%%%%%%%%%%%%%%%%%%%%%%%%%%%%%%%%%%%%%%%%%%%%%%%%%%%%%%%%%%%%%%%%
%%%%%%%%%%%%%%%%%%%%   subsection 0.3

{\bf 0.3. The structure of the paper.} In \S 1 we prove Theorem
0.1. In fact, we will show a much more general fact: the
birational superrigidity of Fano-Mori fibre spaces, satisfying
certain additional assumptions. Those assumptions are
automatically satis\-fied for fibrations into Fano hypersurfaces
of index 1, considered in \cite{Pukh15a}, so that the proof of
Theorem 0.1 immediately implies an essential improvement of the
main theorem of that paper: not only that birational maps onto
rationally connected fibre spaces are fibre-wise, but also the
birational superrigidity in full, that is, the uniqueness of the
structure of a Mori fibre space up to fibre-wise birational
modifications, which are biregular on the generic fibre. We will
prove Theorem 0.1 in the form in which it is stated, and then
identify the conditions that are used in the arguments. The
general fact that follows from the proof of Theorem 0.1, applies
also to Fano-Mori fibre spaces, the fibres of which are multiple
projective spaces of index 1, recently studied in
\cite{Pukh19a,Pukh20a}: these varieties are realized as
hypersurfaces in the weighted projective space, for that reason
they have hypersurface singularities, which may be assumed to be
quadratic with the rank bounded from below.

In \S 2 we prove Theorem 0.2. Except from the condition
(R$2^2$.2), one can estimate the codimension of the set of pairs
$(f_1,f_2)\in{\cal F}$, that do not satisfy the regularity
conditions at at least one point, using the routine techniques
that were used many times and described, for instance, in
\cite{Pukh13a}. However, the condition (R$2^2$.2) forms an
essentially harder problem. We will consider the general question
of estimating the codimension of the set of tuples of polynomials,
the set of common zeros of which is reducible or non-reduced.
Comparing the estimates for the codimension for violation of each
of the regularity conditions, we obtain the claim of Theorem 0.2.

The remaining part of the paper, \S\S 3-5, is given to the proof
of Theorem 0.3. This is the hardest part of the paper. Let us
describe its main stages. In \S 3 Theorem 0.3 is shown in the
assumption that certain claims, global and local, are true.
Assuming that these claims hold, we assume that the pair
$(F,\frac{1}{n}D_F)$ is not canonical, where $D_F\sim n(-K_F)$,
and show that this assumption leads to a contradiction. In order
to obtain the contradiction, we use the technique of hypertangent
divisors \cite[Chapter 3]{Pukh13a}, based on the regularity
conditions for the fibre $F$. Thus we reduce Theorem 0.3 to a
number of facts of global and local geometry.

The global claims that belong to the geometry of complete
intersections of two quadrics in the projective space, are shown
in \S 4. They are needed for the local analysis of the fibres $F$
at bi-quadratic points. For the proof, we use the elementary, but
geometrically non-trivial, technique, developed in \cite[\S
4]{Pukh2018b} for the case of one quadric. In the case of two
quadrics it is essentially more difficult to study maximal linear
subspaces on their intersection.

The local facts, characterizing the behaviour of a non-canonical
singularity of the pair $(F,\frac{1}{n}D_F)$ at a quadratic and a
bi-quadratic points under their blow up, are shown in \S 5. Note
that the local fact for a quadratic singularity simplifies
radically the proof of divisorial canonicity for Fano
hypersurfaces of index 1 with quadratic singularities, given in
\cite{Pukh15a} (Theorem 1.4), and also for mupltiple projective
spaces in \cite{Pukh19a}. A similar radical simplification is
obtained for the main result of \cite{Pukh09b}, which comes as
Theorem 3.2 in \cite[Chapter 7]{Pukh13a}. A detailed discussion of
how the local facts shown in \S 5 improve the known results is
given in Subsection 3.10.\vspace{0.3cm}

%%%%%%%%%%%%%%%%%%%%%%%%%%%%%%%%%%%%%%%%%%%%%%%%%%%%%%%%%%%%%%%%%
%%%%%%%%%%%%%%%%%%   subsection 0.4

{\bf 0.4. Historical remarks.} A majority of results on birational
rigidity are about the absolute case (Fano varieties, considered
as Fano-Mori fibre spaces over a point) and about Fano-Mori fibre
spaces over ${\mathbb P}^1$, see \cite{Pukh13a} and the
bibliography in that book. Except for Sarkisov's theorem on conic
bundles \cite{S80,S82}, up to a recent time the only result of
birational rigidity type for fibrations over a base of dimension
higher than 1 was the theorem on Fano direct products
\cite{Pukh05} (reproduced in \cite[Chapter 7]{Pukh13a}). Note that
Theorem 0.3 allows to use Fano complete intersections of
codimension 2 and index 1 as the class of varieties that can be
taken as direct factors of Fano direct products, preserving
birational rigidity. Of the highest interest, however, are
theorems on birational rigidity of {\it general} Fano-Mori fibre
spaces, realizing the principle that sufficient ``twistedness''
over the base implies birational (super)rigidity and the
uniqueness of the structure of a Mori fibre space on the given
variety. The first result of that type was obtained in
\cite{Pukh15a} for fibrations into primitive Fano hypersurfaces
and double spaces over a fixed base. One of the main challenges in
the work with bibrations over a base of dimension higher than 1 is
that, on one hand, singularities of fibres over some points and
subvarieties become worse, and on the other hand, one may need to
use birational modifications of the base. For one-dimensional
fibres (conics) these problems were within reach 40 years ago.
That is why the uniqueness of the structure of a conic bundle,
which is sufficiently twisted over the base (in other words, with
sufficiently large degenerations), is an old theorem. In
\cite{Pukh15a} it was discovered how to get around this trouble,
considering fibrations into Fano varieties with bounded
singularities, stable with respect to blow ups, and so working
with fibre spaces, the base of which can be blown up, preserving
the properties of a Fano-Mori fibre space. This approach is used
in the present paper, too, however, singularities of complete
intersections can degenerate much worse. For that reason the
technique, developed in \cite{Pukh15a}, needs to be essentially
improved, which is what we do in the present paper (\S\S 2,4 and
5). As a reward, we simplify the proof of the main result in
\cite{Pukh15a} and make the results obtained there stronger.

Theorem 0.3 can be stated in terms of {\it global canonical
thresholds}: it states that for $F\in{\cal F}_{\rm reg}$ the
inequality $\mathop{\rm ct}(F)\geqslant 1$ holds. In fact, in
order to use it in the proof of Theorem 0.1, it is sufficient for
every fibre $F$ of the fibre space $\pi\colon V\to S$ to satisfy,
instead of the condition (i), two weaker conditions: the equality
$\mathop{\rm lct}(F)=1$ for the global {\it log canonical}
threshold and the inequality $\mathop{\rm mct}(F)\geqslant 1$ for
the {\it mobile canonical} threshold; no changes in the proof
(given in \S 1) are needed. These weaker conditions follow from
the divisorial canonicity. One should take into account that the
technique of the proof of Theorem 0.3 (in \S\S 3-5) gives at once
the divisorial canonicity, and replacing that condition by the two
weaker ones in no way simplifies the proof. For that reason the
condition (i) is given in the strongest version.

Computing and estimating the global log canonical threshold became
recently a very popular topic: se,, for instance,
\cite{ChShr2008,ChParkWon2014} and many more papers in the last
2--3 years. The reason is the applications to complex differential
geometry, to such questions as existence of the
K\"{a}hler-Einstein metric and the $K$-stability.

Another topic, popular in the recent years and related to the
theory of birational rigidity, is investigation of groups of
birational self-maps of rationally connected varieties. After the
remarkable result about the Jordan property of the groups of
birational self-maps \cite{ProkhShr2014,ProkhShr2016} the papers
in that direction that was started by Serre's paper
\cite{Serre2009} appear one after another. Note that Theorem 0.1
implies the fact that is so standard for the theory of birational
rigidity that we did not state it explicitly: for a general fibre
space $V\to S$ into complete intersections of codimension 2 the
group of birational self-maps $\mathop{\rm Bir}V$ is trivial.

Theorem 0.1 on birational superrigidity of fibrations into
complete intersections of codimension 2 can be seen as a fragment
of birational classification of rationally connected varieties. In
this context we would like to mention the new approach to the
proof of stable non-rationality via the ``decomposition of the
diagonal'', developed by C.Voisin and used in many recent papers
by many authors, see, for instance,
\cite{C-ThPir2016,HasKreTsch2016,Tot2016,AuBoPir2018,HasPirTsch2019,Schreied2019}.
This method makes it possible to show stable non-rationality of a
very general variety in a given family. Another approach, based on
using the Grothendieck ring, was recently suggested in
\cite{NicShind2019}, see also \cite{KontsevichTschinkel2019}. Note
that the results, obtained by the method of maximal singularities,
in the first place, the theorem on Fano direct products
\cite{Pukh05}, are ``on the opposite pole'' from the papers listed
above: one can take divisorially canonical primitive Fano
varieties (and not ${\mathbb P}^N$) as a direct factor, and gets a
birational classification, which is stable with respect to such
direct products.

Finally, let us mention one more important point. There are two
types of the techniques, used in the theory of birational
rigidity: the linear and quadratic ones. The linear technique is
based on the study of singularities of a general divisor in the
linear system that defines a birational map; it is that technique
that we use in the present paper. For the first time this
technique was developed and used in \cite{Pukh05}. The quadratic
technique is based on the study of singularities of the {\it
self-intersection} of a mobile linear system, defining a
birational map, that is to say, of the scheme-theoretic
intersection of two general divisors in that linear system. The
quadratic technique goes back to the paper of V.A.Iskovskikh and
Yu.I.Manin on the three-dimensional quartic \cite{IM}; almost all
results on the birational rigidity in the absolute case and for
Fano-Mori fibre spaces over ${\mathbb P}^1$ were obtained by means
of that technique. Among the recent papers, where it is applied in
the proof of birational rigidity, we mention
\cite{Krylov2018,AhmKry2017,EvansPukh2}. In
\cite{Pukh16a,Pukh2018b} both techniques were used. It seems that
the quadratic technique can relax the upper bound for the
dimension of the base $S$ of a Fano-Mori fibre space. However, the
attempts to apply it to fibrations over a base of dimension
$\geqslant 2$ meet considerable difficulties. The most that was
possible to do in that direction via the quadratic technique, is
the theorem on birational geometry of fibrations into double
spaces of index 1 in \cite{Pukh2017}.\vspace{0.3cm}

{\bf 0.5. Acknowledgements.} The work was supported by The
Leverhulme Trust (Research Project Grant RPG-2016-279). The author
is grateful to the members of Divisions of Algebraic Geometry and
Algebra at Steklov Institute of Mathematics for the interest to
his work, and also to the colleagues in Algebraic Geometry
research group at the University of Liverpool for general support.
The author is also grateful to the referees for their work on the
paper and a number of useful suggestions.

%%%%%%%%%%%%%%%%%%%%%%%%%%%%%%%%%%%%%%%%%%%%%%%%%%%%%%%%%%%%%%%%%
%%%%%%%%%%%%%%%%%%%%%%%%%%%%%%%%%%%%%%%%%%%%%%%%%%%%%%%%%%%%%%%%%
%%%%%%%%%%%%%%%%%%%   SECTION 1

\section{Birationally rigid Fano-Mori fibre spaces}

In this section we prove Theorem 0.1. In Subsection 1.1 we explain
the main idea of the proof and introduce its most important
constructions: the mobile linear system $\Sigma$ and the mobile
family ${\cal C}$ of irreducible rational curves on $V$. In
Subsection 1.2 we start our study of maximal singularities of the
mobile linear system $\Sigma$. In Subsection 1.3 we apply
fiber-wise modifications to the fibration $V/S$, in order to make
the centres of all maximal singularities cover divisors on the
base. In Subsection 1.4 we show the first claim of Theorem 0.1,
that a birational map in the rationally connected case is
fibre-wise, and in Subsection 1.5 we prove the second claim of
this theorem, the claim on birational rigidity. Finally, in
Subsection 1.6, looking at the properties  of the fibration $V/S$,
which were used in the proof of Theorem 0.1, we state in the most
general form a theorem on birational superrigidity of Fano-Mori
fibre spaces, of which Theorem 0.1 is a particular case (for
fibrations into complete intersections of codimension 2), and the
proof of which is a word for word the same as the arguments of
Subsections 1.1-1.5.\vspace{0.3cm}

{\bf 1.1. Set up of the problem and the plan of the proof.} Fix a
Fano-Mori fibre space $\pi\colon V\to S$, satisfying the
conditions of Theorem 0.1. For the fibration $\pi'\colon V'\to S'$
we will simultaneously consider two options:

--- it is a rationally connected fibre space, that is, the
base $S'$ and the fibre of general position are rationally
connected; in that case we say that the {\it rationally connected
case} is under consideration,

--- it is a Mori fibre space; in that case we say the {\it case
of a Mori fibre space} is under consideration.

Let us fix a fibration $\pi'\colon V'\to S'$, which belong to
either of the two cases, where $\mathop{\rm dim}V'=\mathop{\rm
dim}V$. Assume that there is a birational map $\chi\colon V\to
V'$. Fix it, too. In the rationally connected case we have to show
that the map $\chi$ is fibre-wise, that is, there exists a
rational dominant map of the base $\beta\colon S\dashrightarrow
S'$, such that the diagram
$$
\begin{array}{ccccc}
& V &\stackrel{\chi}{\dashrightarrow}& V' &\\
\pi \!\!\!\! & \downarrow & & \downarrow &\!\!\!\! \pi'\\
& S & \stackrel{\beta}{\dashrightarrow} & S' &\\
\end{array}
$$
commutes. In the case of a Mori fibre space we have to show that,
in addition, the map $\beta$ is birational, so that for a point
$s\in S$ of general position the map $\chi$ induces a birational
isomorphism of the fibres $F_s=\pi^{-1}(s)$ and
$F'_{\beta(s)}=(\pi')^{-1}(\beta(s))$. This isomorphism is
biregular due to birational superrigidity of the fibres $F_s$.

{\bf Remark 1.1.} First of all, note that the inverse birational
map $\chi^{-1}\colon V\dashrightarrow V'$ can be fibre-wise only
in the case when $\chi$ is fibre-wise, and the corresponding map
of the base $\beta\colon S\dashrightarrow S'$ is birational.
Indeed, assume that there exists a rational dominant map
$\beta'\colon S'\dashrightarrow S$, such that
$\beta'\circ\pi'=\pi\circ\chi^{-1}$, where $\beta'$ is not
birational. Then $\mathop{\rm dim}S'>\mathop{\rm dim}S$ and for a
point $s\in S$ of general position the map
$$
\pi'\circ\chi|_{F_s}\colon F_s\dashrightarrow(\beta')^{-1}(s)
$$
fibres the Fano variety $F_s$ over a positive-dimensional variety
$(\beta')^{-1}(s)$ into rationally connected vaieties, birational
to the fibres $F'_t$ of the projection $\pi'$ for
$t\in(\beta')^{-1}(s)$. However, by the condition (i) of Theorem
0.1 every fibre $F_s$ is a birationally superrigid variety and for
that reason has no structures of a rationally connected fibre
space over a positive-dimensional base. This contradiction proves
our claim.

By the remark above, we assume that $\chi^{-1}\colon
V\dashrightarrow V'$ is not fibre-wise. However, $\chi$ still can
be fibre-wise, but in the case of a Mori fibre space this
assumption must lead to a contradiction. Now let us describe the
two main objects, which are of key importance in the proof of
Theorem 0.1.

The first object is the mobile linear system $\Sigma$, connected
with the map $\chi$. In the rationally connected case it is the
same system as in \cite[\S 2]{Pukh15a}: let us consider a very
ample system $\overline{\Sigma'}$ on $S'$ and denote by the symbol
$\Sigma'$ its $\pi'$-pull back on $V'$. Now
\begin{equation}\label{08.05.2020.1}
\Sigma=(\chi^{-1}_*)\Sigma'\subset|-nK_V+\pi^*Y|
\end{equation}
is its strict transform on $V$, where $n\in{\mathbb Z}_+$. In the
case of a Mori fibre space (recall that we consider both cases
simultaneously) for $\Sigma'$ we take a very ample complete linear
system
$$
|-mK'+(\pi')^*Y'|
$$
on $V'$, where $m'\geqslant 1$, the symbol $K'$ stands for the
canonical class $K_{V'}$ and $Y'$ is some very ample divisor on
$S'$. We define the linear system $\Sigma$ by the formula
(\ref{08.05.2020.1}). This is a mobile linear system on $V$, but
in the case of a Mori fibre space we have $n\geqslant 1$. For
uniformity of notations in the rationally connected case we set
$m=0$. The following claim is the key fact in the proof of Theorem
0.1.

{\bf Theorem 1.1.} {\it The inequality $n\leqslant m$ holds.}

It is obvious that Theorem 1.1 immediately implies the claim of
Theorem 0.1 in the rationally connected case: we get the equality
$n=0$, so that the linear system $\Sigma$ is pulled back from the
base $S$ and the map $\chi$ is fibre-wise. It is not very hard to
complete the proof of Theorem 0.1 in the case of a Mori fibre
space, if Theorem 1.1 is shown. So our main purpose is to prove
that theorem. Starting from this moment (and up to the end of the
proof of Theorem 1.1), we assume that the inequality $n>m$ holds;
in particular, in the rationally connected case the map $\chi$ is
not fibre-wise.

The second main object of our proof is the family of curves ${\cal
C}$ on $V$. It is constructed in the same was as in \cite[\S
2]{Pukh15a}. Recall the construction. Let
$\varphi\colon\widetilde{V}\to V$ be the resolution of
singularities of the map $\chi$ and ${\cal E}_{\rm exc}$ the set
of all prime $\varphi$-exceptional divisors $E$ on
$\widetilde{V}$, such that their image on $V'$ is divisorial and
the prime divisor $[\chi\circ\varphi](E)\subset V'$ on $V'$ covers
the base of the fibration $\pi'$, that is,
$$
\pi'([\chi\circ\varphi](E))=S'.
$$
Now let us consider a family of (irreducible rational) curves
${\cal C'}$ on $V'$, contracted by the projection $\pi'$, sweeping
out a dense open subset in $V'$, not intersecting the set where
the rational map
$$
[\chi\circ\varphi]^{-1}\colon V'\dashrightarrow\widetilde{V}
$$
is not determined and intersecting every divisor
$[\chi\circ\varphi](E)$, $E\in{\cal E}_{\rm exc}$, transversally
at points of general position, see \cite[Subsection 2.1]{Pukh15a}.
The curves $C'\in{\cal C'}$ lie in the fibres of the projection
$\pi'$. Let
$$
{\cal C}=\varphi_*\circ[\chi\circ\varphi]^{-1}(\cal C')
$$
be the strict transform of the family ${\cal C'}$ on $V$. This is
a mobile family of (irreducible rational) curves on $V$, sweeping
out a dense open subset of the variety $V$. By Remark 1.1, the
curves $C\in{\cal C}$ are not contained in the fibres of the
projection $\pi$, so that the image
$$
\pi_*{\cal C}=\overline{\cal C}
$$
of that family on the base $S$ is a mobile family of curves on
$S$. For $C\in{\cal C}$ we write $\overline{C}$ for the image
$\pi(C)$. The linear system $\Sigma\subset|-nK_V+\pi^*Y|$ and the
family of curves ${\cal C} $ are the main elements of our
construction. Since the system $\Sigma$ is mobile, by the
condition (ii) of Theorem 0.1 the divisor $Y$ is pseudoeffective,
so that for the curve $C\in{\cal C}$ and its image $\overline{C}$
on $S$ we have
$$
(C\cdot\pi^*Y)=(\overline{C}\cdot Y)\geqslant 0.
$$
The main idea of the proof of Theorem 1.1 is, assuming the
inequality $n>m$, to obtain a contradiction with the condition
(iii) of Theorem 0.1. Obviously, for a general divisor
$D\in\Sigma$ the class of the cycle
$$
(D\circ\pi^{-1}(\overline{C}))\sim
-n(K_V\cdot\pi^{-1}(\overline{C}))+ (\overline{C}\cdot Y)F
$$
is effective. We will show that the scheme-theoretic intersection
in the left hand side contains sufficiently many fibres of the
projection $\pi$, subtracting which we obtain a cycle that is
still effective, which contradicts the condition (iii) of Theorem
0.1.\vspace{0.3cm}

%%%%%%%%%%%%%%%%%%%%%%%%%%%%%%%%%%%%%%%%%%%%%%%%%%%%%%%%%%%%%%%
%%%%%%%%%%%%%%%%%%   subsection 1.2

{\bf 1.2. Maximal singularities.} Let us present the set ${\cal
E}_{\rm exc}$ of all prime $\varphi$-exceptional divisors, the
image of which on $V'$ is divisorial and covers the base $S'$, as
a disjoint union
$$
{\cal E}_{\rm exc}={\cal E}_S\sqcup{\cal E}_{\rm div}\sqcup{\cal
E},
$$
where $E\in{\cal E}_S$ (respectively, ${\cal E}_{\rm div}$ and
${\cal E}$) if and only if the centre $\varphi(E)$ of the divisor
$E$ on $V$ covers the base $S$, that is,
$$
\pi[\varphi(E)]=S
$$
(respectively, covers a prime divisor and an irreducible closed
subset of codimension $\geqslant 2$ on $S$).

Setting $\widetilde{K}=K_{\widetilde{V}}$ and omitting for
simplicity of notations the symbol $\varphi^*$, we get
$$
\widetilde{K}=K_V+\sum\limits_{E\in{\cal E}_{\rm exc}}a_EE+(\dots)
$$
and for the strict transform $\widetilde{\Sigma}$ of the linear
system $\Sigma$ on $\widetilde{V}$ we have
$$
\widetilde{\Sigma}\subset\left|-n\widetilde{K}+\left(\pi^*Y-
\sum\limits_{E\in{\cal E}_{\rm
exc}}\varepsilon(E)E\right)+(\dots)\right|,
$$
where $\varepsilon(E)=b_E-na_E$, $b_E=\mathop{\rm
ord}\nolimits_E\varphi^*\Sigma$ and the symbol $(\dots)$ in both
cases means a linear combination of prime $\varphi$-exceptional
divisors, which are not in the set ${\cal E}_{\rm exc}$ (that is,
either their image on $V'$ is not divisorial, or it is divisorial,
but does not cover the base $S'$ --- those prime divisors are
inessential for our construction). If $\varepsilon(E)>0$, that is,
the {\it Noether-Fano inequality}
$$
b_E>na_E
$$
holds, then the prime divisor $E$ is a {\it maximal singularity in
the strong sense} of the linear system $\Sigma$ (there can be
maximal singularities, satisfying the Noether-Fano inequality,
among the inessential prime divisors in $(\dots)$). The following
claim is true.

{\bf Proposition 1.1.} {\it For at least one $E\in{\cal E}_{\rm
exc}$ the inequality $\varepsilon(E)>0$ holds, that is to say,
maximal singularities in the strong sense do exist.}

{\bf Proof:} this is Proposition 2.1 in \cite{Pukh15a}.

The set of all maximal singularities in the strong sense we denote
by the symbol ${\cal M}$. The following claim somewhat improves
Proposition 2.2 in \cite{Pukh15a}.

{\bf Proposition 1.2.} {\it The inclusion ${\cal M}\subset{\cal
E}$ holds, that is to say, the image of the centre of a maximal
singularity $E\in{\cal M}$ on $V$ under the projection $\pi$ is of
codimension $\geqslant 2$ on $S$.}

{\bf Proof.} Assume the converse: the centre of a maximal
singularity $E$ covers either the base $S$, or a prime divisor on
the base. Let $s\in\pi[\varphi(E)]$ be a point of general
position. Since the linear system $\Sigma$ has no fixed
components, the fibre $F_s=\pi^{-1}(s)$ is not contained in its
base set. For a general divisor $D\in\Sigma$ the pair
$(V,\frac{1}{n}D)$ is not canonical and, moreover, $E$ is a non
canonical singularity of that pair, and $F_s$ is not contained in
the support $|D|$ of the divisor $D$. Denote by the symbol $D_F$
the restriction
$$
D|_{F_s}=(D\circ F_s),
$$
then $D_F\sim-nK_{F_s}$, and if $\mathop{\rm
codim}(\pi[\varphi(E)]\subset S)=1$, then the pair
$(F_s,\frac{1}{n}D_F)$ is not log canonical by inversion of
adjunction, and if $\pi[\varphi(E)]=S$, then the pair
$(F_s,\frac{1}{n}D_F)$ is not canonical --- in any case this
contradicts the condition (i) of Theorem 0.1, which completes the
proof of the Proposition. Q.E.D.

{\bf Remark 1.2.} (i) In the case when the centre $\varphi(E)$
covers the base $S$ (this case is considered in \cite[Proposition
2.2]{Pukh15a}), the divisor $D_F$ is a general divisor of a {\it
mobile} linear system $\Sigma|_{F_s}$, therefore in order to
exclude this case it is sufficient to assume that the {\it mobile}
canonical threshold of the fibre $\geqslant 1$; in the case when
$\varphi(E)$ covers a divisor on $S$, it is sufficient to assume
that $\mathop{\rm lct}(F_s)=1$, which is weaker than the condition
(i) of Theorem 0.1. The divisorial canonicity is stronger than
each of these two conditions.

(ii) The proof of Proposition 1.2 works as it is for any maximal
singularity, that is, a prime $\varphi$-exceptional divisor
$E\subset\widetilde{V}$, satisfying the Noether-Fano inequality
$$
\mathop{\rm ord}\nolimits_E\varphi^*\Sigma>n\cdot a(E,V).
$$
The property that the image on $V'$ is divisorial and covers $S'$
is not used in the proof. This remark will be needed below in the
proof of birational superrigidity of the fibration $V\slash S$.

Now, following \cite[Subsection 2.2]{Pukh15a}, let us consider a
fibre-wise modification of the fibre space $V\slash S$. We need to
ensure that the centre of every maximal singularity in the strong
sense covers at least a divisor on the base of the fibre
space.\vspace{0.3cm}

%%%%%%%%%%%%%%%%%%%%%%%%%%%%%%%%%%%%%%%%%%%%%%%%%%%%%%%%%%%
%%%%%%%%%%%%%%%%   subsection 1.3

{\bf 1.3. Fibre-wise birational modification.} Let $\sigma_S\colon
S^+\to S$ be a composition of blow ups with non-singular centres,
such that the centres of all singularities $E\in{\cal E}$ on
$V^+=V\times_SS^+$ cover a prime divisor on $S^+$. This is the
minimal requirement for the birational morphism $\sigma_S$. Later
we will need an additional property: the base set of the strict
transform of the system $\Sigma$ on $V^+$ does not contain fibres
of the modified Fano-Mori fibre space.

Set $\sigma\colon V^+\to V$ and $\pi_+\colon V^+\to S^+$ to be the
projections of the fibre product, so that the following diagram
commutes:
$$
\begin{array}{ccccc}
& V^+ &\stackrel{\sigma}{\to}& V &\\
\pi_+ \!\!\!\!& \downarrow\phantom{i} & & \downarrow &\!\!\!\!\pi\\
& S^+ & \stackrel{\sigma_S}{\to} & S. &\\
\end{array}
$$
The existence of such modification is shown in \cite[Subsection
2.2]{Pukh15a}. The morphism $\pi_+\colon V^+\to S^+$ is a
Fano-Mori fibre space: the variety $V^+$ is factorial and its
singularities are terminal. This follows easily from the
construction of the fibre space $V\slash S$ in Subsection 0.1, see
Subsection 1.7.

Now let us introduce some new notations. Let ${\cal T}$
(respectively, $\overline{\cal T}$) be the set of prime
$\sigma$-(respectively, $\sigma_S$-)exceptional divisors on $V^+$
(respectively, on $S^+$). The projection $\pi_+$ gives a bijection
${\cal T}\to\overline{\cal T}$. For $T\in{\cal T}$ and
$\overline{T}=\pi_+(T)\in\overline{\cal T}$ we have the equality
of discrepancies
$$
a_T=a(T,V)=a(\overline{T},S).
$$
Besides, let us define the set $\overline{\cal T}_{\rm div}$ of
prime divisors on $S$, which are $\pi$-images of subvarieties
$\varphi(E)\subset V$ for $E\in{\cal E}_{\rm div}$. Respectively,
let ${\cal T}_{\rm div}$ be the set of all prime divisors on $V$
of the form $\pi^{-1}(\overline{T})$, where
$\overline{T}\in\overline{\cal T}_{\rm div}$.

Let $\Sigma^+$ be the strict transform of the linear system
$\Sigma$ on $V^+$. We have:
$$
\Sigma^+\subset \left|-nK_V+\pi^*Y-\sum_{T\in{\cal T}}b_T T\right|
$$
for some $b_T\in{\mathbb Z}_+$ (again we omit the pull back symbol
$\sigma^*$). Similar to \cite[Subsection 2.3]{Pukh15a}, we assume
that the resolution $\varphi\colon\widetilde{V}\to V$ factors
through $\sigma$, that is to say, is of the form
$\psi=\sigma\circ\varphi_+$, where
$\varphi_+\colon\widetilde{V}\to V^+$ is a sequence of blow ups
with non-singular centres. By construction, there is a map
$$
\lambda\colon{\cal E}\to{\cal T},
$$
where $\lambda(E)=T\in{\cal T}$ is the uniquely determined
$\sigma$-exceptional divisor, for which
$$
\varphi_+(E)\subset T,
$$
and we have $\pi_+[\varphi_+(E)]=\overline{T}=\pi_+(T)$. In a
similar way, we define the map
$$
\lambda_{\rm div}\colon{\cal E}_{\rm div}\to{\cal T}_{\rm div},
$$
For the discrepancies we have the obvious equality
$$
a_E=a^+_E+a_T\cdot\mathop{\rm ord}\nolimits_E\varphi^*_+T,
$$
where $T=\lambda(E)$ and $a^+_E=a(E,V^+)$. It may happen that
$a^+_E=0$: this happens precisely when $E=T$ (and $a_E=a_T$). For
exceptional divisors $E\in{\cal E}_S\sqcup{\cal E}_{\rm div}$ we
have the equality $a_E=a^+_E$. Furthermore, the mobile linear
system $\widetilde{\Sigma}$ is a subsystem of the complete linear
system
$$
\left|\varphi^*_+\left(-nK_V+\pi^*Y-\sum_{T\in{\cal T
}}b_TT\right)-\sum_{E\in{\cal E}_{\rm exc}}b^+_E E-(\dots)\right|,
$$
where $b^+_E\in{\mathbb Z}_+$ and the symbol $(\dots)$ has the
same meaning as before. For $E\in{\cal E}$ we have the equality
$$
b_E=b^+_E+b_T\cdot\mathop{\rm ord}\nolimits_E\varphi^*_+T,
$$
where $T=\lambda(E)$. Again, $b^+_E=0$ if and only if $E=T$.

{\bf Proposition 1.3.} {\it For every $E\in{\cal E}_{\rm exc}$
the following inequality holds:}
$$
b^+_E\leqslant n\cdot a^+_E.
$$

{\bf Proof} repeats the proof of Proposition 1.2 almost word for
word: assume the converse, that is, that the inequality
$b^+_E>n\cdot a^+_E$ holds. Then $b^+_E>0$, and for that reason
$E\neq T$ for $T=\lambda(E)$, that is, the centre of $E$ has the
codimension $\geqslant 2$ on $V^+$, so that $a^+_E>0$.
Furthermore, the pair $(V^+,\frac{1}{n}D^+)$ is not canonical for
$D^+\in\Sigma^+$ by our assumption. Since the linear system
$\Sigma^+$ is mobile, we may assume that for a point
$s\in\overline{T}$ of general position the fibre
$F_s=\pi^{-1}_+(s)$ is not contained in the support of the divisor
$D^+$, so that the pair $(F_s,\frac{1}{n}D^+_F)$ is not log
canonical, where $D^+_F=(D^+\circ F)\sim-nK_{F_s}$. This
contradicts the assumption about the properties of the fibres of
the fibre space $V\slash S$ (which are also the fibres of the
fibre space $V^+\slash S^+$). Q.E.D. for the proposition.

{\bf Corollary 1.1.} {\it For a maximal singularity $E\in{\cal
M}$ the inequality
$$
b_T>n\cdot a_T
$$
holds, where $T=\lambda(E)$}.

{\bf Proof:} this follows immediately from the inequality $b_E>n\cdot
a_E$, the previous proposition and the explicit formulas for
$b_E$ and $a_E$, given above. Q.E.D. for the corollary.

Now let us go back to the mobile family of irreducible rational
curves ${\cal C}$ on $V$, constructed in Subsection 1.1. The
family ${\cal C}$ is the strict transform of the mobile family
${\cal C'}$ on $V'$. In the rationally connected case the linear
system $\Sigma'$ is pulled back from the base $S'$, so that for a
divisor $D'\in\Sigma'$ and a curve $C'\in{\cal C'}$ we have the
equality
$$
(C'\cdot D')=0.
$$
In the case of a Mori fibre space we have the equality
$$
(C'\cdot D')=-m(C'\cdot K')\geqslant m,
$$
so that for $l>m$ we get
$$
(C'\cdot [D'+lK])=(l-m)(C'\cdot K')\leqslant-(l-m).
$$

Let us denote by the symbols $\widetilde{\cal C}$, ${\cal C}^+$
and $\overline{\cal C^+}$ the strict transforms of the family
${\cal C}$ on $\widetilde{V}$, $V^+$ and the image of the family
${\cal C}^+$ on $S^+$, respectively. The corresponding symbols
will be used for the general curve:
$\widetilde{C}\in\widetilde{\cal C}$, $C^+\in{\cal C}^+$ and
$\overline{C^+}(=\pi_+(C^+))\in{\overline{\cal C^+}}$. The family
$\overline{\cal C^+}$ of curves on $S^+$ is mobile and sweeps out
a dense open subset of the base, so that for a general divisor
$D^+\in\Sigma^+$ we get a well defined algebraic cycle of the
scheme-theoretic intersection
$$
(D^+\circ\pi^{-1}_+(\overline{C^+})),
$$
which is an effective cycle of dimension $\mathop{\rm dim}F$ on
$V^+$. Estimating the class of that cycle, we will complete the
proof of Theorem 1.1.\vspace{0.3cm}

%%%%%%%%%%%%%%%%%%%%%%%%%%%%%%%%%%%%%%%%%%%%%%%%%%%%%%%%%%%
%%%%%%%%%%%%%%%%   subsection 1.4

{\bf 1.4. Restriction onto the inverse image of a curve.}
Obviously,
$$
(D^+\circ\pi^{-1}_+(\overline{C^+}))
\sim -n(\sigma^*K_V\cdot\pi^{-1}_+(\overline{C^+}))+bF,
$$
where
$$
b=\left(\left[\sigma^*_S Y-\sum_{\overline{T}\in\overline{\cal T
}}b_T\overline{T}\right]\cdot\overline{C^+}\right).
$$
The following claim is of key importance.

{\bf Proposition 1.4.} (i) {\it In the rationally connected case
the following inequality holds:}
$$
b\leqslant -2n.
$$

(ii) {\it In the case of a Mori fibre space the following
inequality holds:}
$$
b\leqslant -n.
$$

{\bf Proof.} Since $(\pi_+)_*C^+=\overline{C^+}$, we get
$$
b=\left(\left[\pi^*Y-\sum_{\overline{T}\in\overline{\cal T
}}b_TT\right]\cdot C^+\right)
$$
(the symbol $\sigma^*$ is omitted, as usual). Here $(C^+\cdot
T)\neq 0$ only for the divisors $T\in\lambda({\cal E})$, since the
curve $C'\in{\cal C'}$ does not meet the set where the map
$(\chi\circ\varphi)^{-1}$ is not defined by construction. Since
the families ${\cal C}$ and $\widetilde{\cal C}$ sweep out dense
open subsets of the varieties $V$ and $\widetilde{V}$,
respectively, we get the inequalities
$$
(K_V\cdot C)<0\quad\mbox{and}\quad
(\widetilde{K}\cdot\widetilde{C})<0
$$
for $C=\sigma(C^+)\in{\cal C}$ and
$\widetilde{C}\in\widetilde{\cal C}$. On the other hand, for a
general divisor $\widetilde{D}\in\widetilde{\Sigma}$ in the
rationally connected case we have the equality
$$
(\widetilde{D}\cdot\widetilde{C})=(D'\cdot C')=0,
$$
and in the case of a Mori fibre space the equality
$$
(\widetilde{D}\cdot\widetilde{C})=(D'\cdot C')=-m(K'\cdot C').
$$
Since
$$
\widetilde{D}\sim-n\widetilde{K}+\left[\pi^*Y-\sum_{E\in{\cal E
_{\rm exc}}}(b_E-na_E)E\right|+(\dots),
$$
where $(\dots)$ is a linear combination of divisors that do not
meet $\widetilde{C}$, we get the inequality
$$
\left(\left[\pi^*Y-\sum_{E\in{\cal E _{\rm
exc}}}(b_E-na_E)E\right]\cdot\widetilde{C}\right)\leqslant-(n-m)
$$
(recall that in the rationally connected case we set $m=0$).
However, for $E\in{\cal E}_S\sqcup{\cal E}_{\rm div}$ we have
$$
b_E-na_E=b^+_E-na^+_E\leqslant 0
$$
(see Proposition 1.2), and for $E\in{\cal E}$ the equality
$$
b_E-na_E=(b^+_E-na^+_E)+\mathop{\rm ord}\nolimits_E\varphi^*_+T
\cdot(b_T-na_T)
$$
holds, where $T=\lambda(E)\in{\cal T}$, so that taking into
account Proposition 1.3, we obtain the inequality
$$
\left(\left[\pi^*Y-\sum_{T\in\lambda(\cal
E)}(b_T-na_T)\sum_{E\in\lambda^{-1}(T)}(\mathop{\rm
ord}\nolimits_E\varphi^*_+T)E\right]\cdot\widetilde{C}\right)\leqslant-(n-m).
$$
However, for every $T\in\lambda({\cal E})$ the divisor
$$
\varphi^*_+T-\sum\limits_{E\in\lambda^{-1}(T)}(\mathop{\rm
ord}\nolimits_E\varphi^*_+T)E
$$
(where the strict transform $\widetilde{T}$ of the divisor $T$ on
$\widetilde{V}$ can also come into the right hand sum, if
$\widetilde{T}\in\lambda^{-1}(T)$, and in that case, obviously,
$\mathop{\rm ord}_{\widetilde{T}}\varphi^*_+T=1$) is effective and
has the zero intersection with $\widetilde{C}$, since it does not
contain any irreducible components which are divisorial on $V'$.
Therefore, the intersection number
$$
\left(\left[\pi^*Y-\sum_{T\in\lambda(\cal
E)}(b_T-na_T)\varphi^*_+T\right]\cdot\widetilde{C}\right)=
\left(\left[\pi^*Y-\sum_{T\in{\cal T}}(b_T-na_T)T\right]\cdot
C^+\right)
$$
does not exceed the number $-(n-m)$ (recall that for $T\in{\cal
T}\backslash\lambda({\cal E})$ we have the equality $(T\cdot C^+)=0$).
Therefore,
$$
b=\left(\left[\pi^*Y-\sum_{T\in{\cal T}}b_T T\right]\cdot
C^+\right)\leqslant-(n-m)-n\left(\left[\sum_{T\in{\cal
T}}a_TT\right]\cdot C^+ \right)\leqslant-2n+m.
$$
This completes the proof of Proposition 1.4. Q.E.D.

Therefore,
$$
\sigma_*(D^+\circ\pi^{-1}_+(\overline{C^+}))
\sim-n(K_V\cdot\pi^{-1}(\overline{C}))+bF,
$$
where $b\leqslant -n-(n-m)$, and since the algebraic cycle in the
left hand side is effective, the class
$$
-(K_V\cdot\pi^{-1}(\overline{C}))-F
$$
turns out to be effective, which contradicts the condition (iii)
of Theorem 0.1. We assumed that $n>m$ and obtained a
contradiction. This completes the proof of Theorem 0.1. Q.E.D.

In particular, in the rationally connected case (where $m=0$) we
get the equality $n=0$, that is, the map $\chi$ is fibre-wise: the
first (``rationally connected'') claim of Theorem 0.1 ias shown.
Let $\beta\colon S\dashrightarrow S'$ be the rational dominant
map, satisfying the equality $\beta\circ\pi=\pi'\circ\chi$.

Now we start the proof of the second claim of Theorem 0.1, of the
birational superrigidity of the Fano-Mori fibre space $V\slash S$.

First of all, we make the following observation.

{\bf Remark 1.3.} It is easy to see from the proof of Proposition
1.2 or Proposition 1.3, that if a prime $\varphi_+$-exceptional
divisor $E\subset\widetilde{V}$ is a singularity of the linear
system $\Sigma^+$, that is, the inequality
$$
b^+_E=\mathop{\rm ord}\nolimits_E\varphi^*_+\Sigma^+>0
$$
holds, and moreover, the closed subset
$$
\pi^{-1}_+(\pi_+[\varphi_+(E)])\subset V^+
$$
(consisting of all fibres of the projection $\pi_+$ that contain
at least one point of the centre of $E$ on $V^+$) is not contained
entirely in the base set $\mathop{\rm Bs}\Sigma^+$, then the
inequality
$$
b_E\leqslant na_E
$$
holds. Indeed, our assumption implies that for a point of general
position
$$
s\in\pi_+[\varphi_+(E)]\subset S^+
$$
the fibre $F_s$ is not contained in the support of a general
divisor $D^+\in\Sigma^+$, so that the proof of Proposition 1.3 (or
1.2), based on the inversion of adjunction, works without any
changes.\vspace{0.3cm}

%%%%%%%%%%%%%%%%%%%%%%%%%%%%%%%%%%%%%%%%%%%%%%%%%%%%%%%%%%%%%
%%%%%%%%%%%%%%%%%   subsection 1.5

{\bf 1.5. Proof of birational rigidity.} We have to show that the
map $\beta$ is birational. As above, let us assume that this is
not the case. Since a general fibre of a Mori fibre space is
irreducible, this implies that $\mathop{\rm dim} S'<\mathop{\rm
dim} S$ and a general fibre of the map $\beta$ is an irreducible
variety of positive dimension. Let us show that this is
impossible.

Using the notations of Subsection 1.3, let us assume in addition,
that the sequence of blow ups $\sigma_S\colon S^+\to S$ resolves
the singularities of the map $\beta$, so that the composite map
$$
\beta_+=\beta\circ\sigma_S\colon S^+\to S'
$$
is a morphism (the general fibre of which is irreducible). Denote
the composite map
$$
\chi\circ\sigma\colon V^+\dashrightarrow V'
$$
by the symbol $\chi_+$. Since we completed the study of the
rationally connected case, we are in the case of a Mori fibre
space, so that
$$
\Sigma'=|-mK'+(\pi')^*Y'|
$$
is a very ample complete linear system on $V'$ and
$m\geqslant 1$. Its strict transform on $V$,
$$
\Sigma\subset|-nK_V+\pi^*Y|,
$$
by Theorem 1.1 satisfies the inequality $n\leqslant m$. Consider
the strict transform $\Sigma^+$ of the system $\Sigma$ on $V^+$.
In addition to the assumptions about the birational morphism
$\sigma_S$, we will assume that $\sigma_S$ {\it flattens}
singularities of the linear system $\Sigma$ over $S$ in the
following sense: no fibre $F_s$, $s\in S^+$, of the projection
$\pi_+$ is contained in the base set $\mathop{\rm Bs}\Sigma^+$.
Now we get:
$$
\Sigma^+\subset|-nK^++\pi^*_+Y^+|,
$$
where $K^+=K_{V^+}$, and $Y^+$ is some divisor on the base $S^+$.
Again, let us consider the resolution
$\varphi_+\colon\widetilde{V}\to V^+$ of singularities of the map
$\chi_+=\chi\circ\sigma\colon V^+\dashrightarrow V'$ and let
$\varphi'=\chi_+\circ\varphi_+\colon\widetilde{V}\to V'$ be the
corresponding birational morphism. By the symbol ${\cal E'}$ we
denote the set of all prime $\varphi'$-exceptional divisors on
$\widetilde{V}$, covering the base $S'$, and by the symbol ${\cal
E}^+$ the set of all prime $\varphi_+$-exceptional divisors on
$\widetilde{V}$ (we no longer require the image of $E\in{\cal
E}^+$ on $V'$ to be divisorial). For the canonical class
$\widetilde{K}=K_{\widetilde{V}}$ we have the equality:
\begin{equation}\label{09.04.2020.1}
\widetilde{K}=K^++\sum_{E\in{\cal E}^+}a^+_E E=K'+\sum_{E'\in{\cal
E'}}a(E')E'+(\dots),
\end{equation}
where $a^+_E$ and $a(E')$ are the discrepancies with respect to
$V^+$ and $V'$, respectively, and in the brackets $(\dots)$ there
is an effective linear combination of $\varphi'$-exceptional prime
divisors on $\widetilde{V}$, the image of which on $V'$ does not
cover the base $S'$.

Let $F'=F'_t=(\pi')^{-1}(t)$ be a fibre of general position of the
projection $\pi'$. Obviously, the divisors in $(\dots)$ do not
meet the strict transform $\widetilde{F'}$ on $\widetilde{V}$.
Denote by the symbol $G$ the subvariety
$$
\pi^{-1}_+(\beta^{-1}_+(t))\subset V^+.
$$
Obviously, its strict transform $\widetilde{G}$ on
$\widetilde{V}$ is $\widetilde{F'}$. For the strict transform
$\widetilde{\Sigma}$ of the system $\Sigma^+$ on $\widetilde{V}$ we have:
$$
\widetilde{\Sigma}\subset\left|-nK^++\pi^*_+Y^+-\sum_{E\in{\cal
E}^+}b^+_EE\right|
$$
(omitting the pull back symbols $\varphi^*_+$, as usual), so that
the following equality of divisorial classes holds:
$$
-mK'+(\pi')^*Y'=-nK^++\pi^*_+Y^+-\sum_{E\in{\cal E}^+}b^+_EE.
$$
From here easy computations, using the formula
(\ref{09.04.2020.1}), lead us to the equality
\begin{equation}\label{09.04.2020.3}
(m-n)K^++(\pi^*_+Y^+-(\pi')^*Y')=\sum_{E\in{\cal
E}^+}(b^+_E-ma^+_E)E+m\sum_{E'\in{\cal E'}}a(E')E'+(\dots),
\end{equation}
where the symbol $(\dots)$ has the same meaning as in the formula
(\ref{09.04.2020.1}). Since by construction no fibre of the
projection $\pi_+$ is contained in $\mathop{\rm Bs}\Sigma^+$, by
Remark 1.3 we get:
$$
b^+_E\leqslant na^+_E\leqslant ma^+_E.
$$
Let $F=F_s\subset G$ be a general fibre of the projection $\pi_+$
and $\widetilde{F}$ its strict transform on $\widetilde{V}$. Since
$$
\pi^*_+Y^+|_F=0\quad \mbox{and} \quad
(\pi')^*Y'|_{\widetilde{G}}=(\dots)|_{\widetilde{G}}=0,
$$
restricting the equality (\ref{09.04.2020.3}) onto
$\widetilde{F}$, we get:
$$
(m-n)K^+|_{\widetilde{F}}=\sum_{E\in{\cal
E^+}}(b^+_E-ma^+_E)E|_{\widetilde{F}}+m\sum_{E'\in{\cal
E'}}a(E')E'|_{\widetilde{F}}.
$$
Applying $(\varphi_+)_*$ to both sides, we obtain an equality of
divisorial classes on $F$, where on the left, if $n<m$, we have a
negative divisor, and on the right an effective divisor
$$
m\sum_{E'\in{\cal E'}}a(E')(\varphi_+)_*E'|_{\widetilde{F}}.
$$
Therefore, we conclude that $n=m$. Now, restricting the equality
(\ref{09.04.2020.3}) onto $\widetilde{G}$, we obtain:
\begin{equation}\label{09.04.2020.4}
\pi^*_+Y^+|_{\widetilde{G}}+\sum_{E\in{\cal
E^+}}(na^+_E-b^+_E)E|_{\widetilde{G}}=\sum_{E'\in{\cal
E'}}a(E')E'|_{\widetilde{G}}.
\end{equation}
Applying $(\varphi_+)_*$ again and taking into account that all
discrepancies $a(E')$ are strictly positive, we get that the
divisorial class
$$
(\pi^*_+Y^+)|_G=\pi^*_+(Y^+|_{\beta^{-1}_+(t)})
$$
is effective. Therefore, in the equality (\ref{09.04.2020.4}) we
have effective divisors on the left and on the right. However, all
divisors $E'|_{\widetilde{G}}$ are exceptional for the birational
morphism $\varphi'|_{\widetilde{G}}\colon\widetilde{G}\to F'$, so
that any linear combination of them with positive coefficients is
a fixed divisor. We conclude that on the variety
$\overline{G}=\beta^{-1}_+(t)=\pi_+(G)\subset S^+$ there is a
finite set of prime divisors $Y_i$, $i\in I$, such that
$$
Y^+|_{\overline{G}}\sim\sum_{i\in I}m_iY_i
$$
for some $m_i\geqslant 1$, and moreover, for any $m^*_i\in{\mathbb
Z}_+$ the complete linear system $\left|\sum\limits_{i\in
I}m^*_iY_i\right|$ consists of a single divisor, the divisor
$\sum_{i\in I}m^*_iY_i$ itself. Consider the ${\mathbb Q}$-vector
space
$$
{\cal A}=\mathbb {\rm Pic\,}\widetilde{G}\otimes{\mathbb Q}
$$
and its subspace ${\cal D}\subset{\cal A}$, generated by the
classes of irreducible components of the divisors
$E'|_{\widetilde{G}}$ for all $E'\in{\cal E'}$. Obviously, there
is the equality
$$
{\cal A}={\mathbb Q}K'\oplus{\cal D},
$$
since the Picard number of the fibre $F'$ is 1. On the other hand,
replacing in (\ref{09.04.2020.4})
$$
\pi^*_+Y^+|_{\widetilde{G}}\quad\mbox{by} \quad
(\varphi_+|_{\widetilde{G}})^*(\pi_+|_G)^*\sum_{i\in I}m_iY_i,
$$
we obtain, since the divisor in the right hand side is fixed, a
physical equality (and not just a linear equivalence) of two
effective divisors, which implies that the hyperplane ${\cal
D}\subset{\cal A}$ is generated by the classes of irreducible
components of the divisors $E|_{\widetilde{G}}$, $E\in{\cal E^+}$,
and the divisors
$$
(\varphi_+|_{\widetilde{G}})^*(\pi_+|_G)^*Y_i,\quad i\in I.
$$
Now let us show that the equality
$$
\mathop{\rm dim}\nolimits_{\mathbb Q}{\cal A}\slash{\cal D}=1
$$
is impossible. Indeed, the class of any ample divisor
$\overline{\Delta}$ on the variety $\overline{G}$ in the space
$\mathop{\rm Pic}\overline{G}\otimes{\mathbb Q}$ can not belong to
a subspace, generated by $Y_i$, $i\in I$, because it is mobile.
Therefore,
$$
\Delta=(\pi_+|_G)^*\overline{\Delta}\not\in{\cal D}
$$
(we omit the symbol $(\varphi_+|_{\widetilde{G}})^*$ of pulling
back onto $\widetilde{G}$ again). This implies that
$$
{\cal A}={\mathbb Q\Delta\oplus{\cal D}}.
$$
However a curve $\Gamma\subset F$ of general position in a fibre
$F$ of general position of the projection $\pi_+|_G$ does not
intersect the set where the map
$$
(\varphi_+|_{\widetilde{G}})^{-1}\colon
G\dashrightarrow\widetilde{G}
$$
is not defined and for that reason its inverse image
$\widetilde{\Gamma}=\Gamma$ on $\widetilde{G}$ has the zero
intersection with any divisor, the class of which is contained in
${\cal D}$, and with each divisor that is pulled back from the
base $\overline{G}$. Therefore, the curve $\widetilde{\Gamma}$ has
the zero intersection with any divisor on $\widetilde{G}$, which
is absurd. This contradiction shows that the rational dominant map
$\beta\colon S\dashrightarrow S'$ is birational and, therefore,
$\chi$ maps birationally a fibre of general position $F_s$ onto a
fibre of general position $F'_{\beta(s)}$. By the birational
superrigidity of the variety $F_s$, terminality and ${\mathbb
Q}$-factoriality of singularities of the variety $F'_{\beta(s)}$
and the equality $\rho(F'_{\beta(s)})=1$, we conclude that the
birational map
$$
\chi|_{F_s}\colon F_s\dashrightarrow F'_{\beta(s)}
$$
is an isomorphism. Q.E.D. for Theorem 0.1.\vspace{0.3cm}

%%%%%%%%%%%%%%%%%%%%%%%%%%%%%%%%%%%%%%%%%%%%%%%%%%%%%%%%%%%%
%%%%%%%%%%%%%%%%%%%%%%%%   subsection 1.6

{\bf 1.6. A generalization of Theorem 0.1.} In the proof of
Theorem 0.1 (including the proof of Theorem 1.1) we never used the
explicit construction of the fibre space $V\slash S$, but we used
its properties that followed from its construction. The claim of
Theorem 0.1 is true for any Fano-Mori fibre space, satisfying
those properties. Let $\pi\colon V\to S$ be an arbitrary Fano-Mori
fibre space over a non-singular rationally connected base $S$,
that is, the variety $V$ is factorial, has terminal singularities
the equality
$$
\mathop{\rm Pic}V={\mathbb Z}K_V\oplus\pi^*\mathop{\rm Pic}S
$$
holds and the anticanonical class $(-K_V)$ is relatively ample.
(The varieties $V$ and $S$ are assumed to be projective.)

{\bf Definition 1.1.} A Fano-Mori fibre space $\pi\colon V\to S$,
every fibre of which is irreducible and reduced, is {\it stable
with respect to fibre-wise birational modifications}, if for any
birational morphism $\sigma_S\colon S^+\to S$, where $S^+$ is
non-singular, the corresponding morphism
$$
\pi_+\colon V^+=V\times_SS^+\to S^+
$$
is a Fano-Mori fibre space (that is, $V^+$ is factorial and its
singularities are terminal).

Now the proof Theorem 0.1, given in Subsections 1.1-1.5, yields
the following general fact.

{\bf Theorem 1.2.} {\it Assume that the Fano-Mori fibre space
$\pi\colon V\to S$, every fibre of which is irreducible and
reduced, is stable with respect to fiber-wise birational
modifications. Assume, in addition, that the following conditions
are satisfied:

{\rm (i)} for any point $s\in S$ the global log canonical
threshold of the corresponding fibre $F_s=\pi^{-1}(s)$ is
$\mathop{\rm lct}(F_s)=1$ and its mobile canonical threshold
satisfies the inequality $\mathop{\rm mct}(F_s)\geqslant 1$,

{\rm (ii)} the $K$-condition holds, that is, for every effective
divisor
$$
D\in|-nK_V+\pi^*Y|
$$
where $n\geqslant 1$, the divisorial class $Y$ is pseudoeffective
on $S$,

{\rm (iii)} for every mobile family of irreducible rational curves
$\overline{\cal C}$ on the base $S$, sweeping out a dense open
subset of $S$, and a curve $\overline{C}\in\overline{\cal C}$ no
positive multiple of the algebraic cycle
$$
-(K_V\cdot\pi^{-1}(\overline{C}))-F,
$$
where $F$ is the class of a fibre of the projection $\pi$, is
effective.

Then for any rationally connected fibre space $V'\slash S'$ every
birational map $\chi\colon V\dashrightarrow V'$ (if there are such
maps) is fibre-wise, and the fibre space $V\slash S$ is
birationally superrigid.}\vspace{0.3cm}

%%%%%%%%%%%%%%%%%%%%%%%%%%%%%%%%%%%%%%%%%%%%%%%%%%%%%%%%%%%%%%%%%
%%%%%%%%%%%%%%%%%%%%%   subsection 1.7

{\bf 1.7. Quadratic and bi-quadratic singularities.} Let us prove
that Fano-Mori fibre spaces $V/S$, constructed in Subsections
0.1-0.2, are stable with respect to fibre-wise birational
modifications. The proof is not hard, but it is convenient to do
it in a more general context.

First of all, recall that an algebraic variety ${\cal X}$ is a
variety with at most quadratic singularities of rank $\geqslant
r$, if in a neighborhood of every point $o\in {\cal X}$ this
variety can be realized as a hypersurface in a non-singular
variety ${\cal Y}$ with a local equation at the point $o$ of the
form
$$
0=\beta_1(u_*)+\beta_2(u_*)+\dots,
$$
where $(u_*)$ is a system of local parameters on ${\cal Y}$ at the
point $o$, and either the linear form $\beta_1\not\equiv 0$, that
is, the point $o\in{\cal X}$ is non-singular, or $\beta_1\equiv 0$
and then the quadratic form $\beta_2$ is of rank $\geqslant r$. In
\cite[Subsection 3.1]{Pukh15a} it was shown that the property to
have quadratic singularities of rank $\geqslant r$ is stable with
respect to blow ups in the following sense: let $B\subset{\cal X}$
be an irreducible subvariety, then there is an open set ${\cal
U}\subset{\cal Y}$, such that

1) ${\cal U}\cap B\neq\emptyset$, where the variety ${\cal U}\cap
B$ is non-singular,

2) for its blow up $\sigma_B\colon {\cal U}_B\to {\cal U}$ the
strict transform of the intersection ${\cal X}\cap {\cal U}$,
which is a hypersurface ${\cal X}_B\subset{\cal U}_B$, is again a
variety with at most quadratic singularities of rank $\geqslant
r$.

It is easy to see that if ${\cal X}$ is a variety with quadratic
singularities of rank $\geqslant 5$, then it is factorial and its
singularities are terminal. In the proof of stability with respect
to blow ups, given in \cite{Pukh15a}, the following obvious fact
is used. Let ${\cal X}$ be a hypersurface in a non-singular
variety ${\cal Y}$ and ${\cal Z}\subset{\cal Y}$ a non-singular
hypersurface, and $o\in{\cal X}\cap{\cal Z}$ is some point. Let
$\beta(u_*)=0$ be a local equation of ${\cal X}$ at the point $o$
with respect to the system of parameters $(u_*)$ at that point. If
the equation
$$
\beta|_{\cal Z}=0
$$
defines a hypersurface that has at the point $o$ a quadratic
singularity of rank $\geqslant r$, then ${\cal X}$ also has at the
point $o$ a singularity of that type. It follows from this
observation that if ${\cal X}\cap({\cal Y}\backslash{\cal Z})$ has
at most quadratic singularities of rank $\geqslant r$ and the same
is true for the restriction ${\cal X}|_{\cal Z}$, then the
hypersurface ${\cal X}$ has at most quadratic singularities of
rank $\geqslant r$. If we take for the hypersurface ${\cal Z}$ the
exceptional divisor of the blow up of the subvariety ${\cal U}
\cap B$, then we get precisely the stability of singularities of
that type with respect to blow ups.

Now let us consider bi-quadratic singularities.

Fix a pair $(r_1,r_2)\in{\mathbb Z}^2_+$ with $r_2\geqslant
r_1+2$.

Let us define the class of varieties with at most quadratic
singularities of rank $\geqslant r_1$ and bi-quadratic
singularities of rank $\geqslant r_2$: a variety ${\cal X}$
belongs to that class, if in a neighborhood of any point
$o\in{\cal X}$ this variety can be realized as a complete
intersection of codimension 2 in a non-singular variety ${\cal Y}$
with local equations at the point $o$ of the form
$$
\begin{array}{c}
0=\beta_{1,1}(u_*)+\beta_{1,2}(u_*)+\dots,\\
0=\beta_{2,1}(u_*)+\beta_{2,2}(u_*)+\dots,\\
\end{array}
$$
where $(u_*)$ is a system of local parameters on ${\cal Y}$ at the
point $o$, and precisely one of the following three cases takes
place:

\begin{itemize}

\item the linear forms $\beta_{1,1}(u_*)$ and $\beta_{2,1}(u_*)$
are linearly independent, and then ${\cal X}$ is non-singular at
the point $o$,

\item $\beta_{1,1}=\alpha_1\tau$ and $\beta_{2,1}=\alpha_2\tau$,
where $\tau(u_*)$ is non-zero linear form and
$(\alpha_1,\alpha_2)\neq(0,0)$, where, moreover, the rank of the
quadratic form
$$
(\alpha_2\beta_{1,2}-\alpha_1\beta_{2,2})|_{\{\tau=0\}}
$$
is at least $r_1$, and then the variety ${\cal X}$ has at the
point $o$ a quadratic singularity of rank $\geqslant r_1$,

\item $\beta_{1,1}\equiv\beta_{2,1}\equiv 0$ and the rank of
the pair of quadratic forms $\mathop{\rm
rk}(\beta_{1,2},\beta_{2,2})\geqslant r_2$.

\end{itemize}

It is easy to see that
$$
\mathop{\rm codim}(\mathop{\rm Sing}{\cal X}\subset{\cal X})
\geqslant\mathop{\rm min}(r_1-1,r_2-3).
$$

For brevity we will say that ${\cal X}$ has {\it singularities of
type} $(r_1,r_2)$.

{\bf Theorem 1.3.} {\it Assume that ${\cal X}$ has singularities
of type $(r_1,r_2)$ and $B\subset{\cal X}$ is an irreducible
subvariety. There exists an open set ${\cal U}\subset{\cal X}$,
such that ${\cal U}\cap B\neq\emptyset$, where ${\cal U}\cap B$ is
a non-singular subvariety and the blow up $\sigma_B\colon {\cal
U}_B\to {\cal U}$ along $B$ gives a variety ${\cal U}_B$ with
singularities of type $(r_1,r_2)$.}

{\bf Proof.} By the stability of quadratic singularities, shown in
\cite{Pukh15a}, it is sufficient to consider the case when $B$ is
entirely contained in the set of bi-quadratic singularities. Let
${\cal U}$ be an open set, such that the subvariety ${\cal U}\cap
B$ is non-singular, the ranks $\mathop{\rm rk}\beta_{1,2}$,
$\mathop{\rm rk}\beta_{2,2}$ and $\mathop{\rm
rk}(\beta_{1,2},\beta_{2,2})$ are constant along ${\cal U}\cap B$.
Considering ${\cal U}\cap B$ as a subvariety of a non-singular
variety ${\cal Y}$, given by a pair of equations, let us choose a
system of local parameters, with respect to which $B$ is given by
a system of equations
$$
u_1=\dots=u_k=0.
$$
Blowing up $B\subset{\cal Y}$ with the exceptional divisor
$E_B\subset\widetilde{\cal Y}$, we realize the blow up ${\cal
U}_B$ as a subvariety of codimension 2 in $\widetilde{\cal Y}$,
defined by the pair of equations
$\widetilde{\beta_1}=\widetilde{\beta_2}=0$. It is sufficient to
check that the singularities of ${\cal U}_B$ on the intersection
${\cal U}_B\cap E_B$ are either quadratic of rank $\geqslant r_1$,
or bi-quadratic of rank $\geqslant r_2$. It is clear that the
system of equations
$$
\widetilde{\beta_1}|_{E_B}=\widetilde{\beta_2}|_{E_B}=0
$$
defines an irreducible reduced subvariety of codimension 2 in
$E_B$, fibred over $B$, where the fibre over a point $b\in B$ is a
complete intersection of two quadrics of rank $\geqslant r_2$ in
${\mathbb P}^{k-1}$. Denote this fibre by the symbol $E_b$. If
$p\in E_b$ is a non-singular point, then ${\cal U}_B$ is also
non-singular at the point $p$. If $p\in E_b$ is a quadratic
singularity, then its rank $\geqslant r_2-2\geqslant r_1$, because
the rank of a quadratic form drops by at most 2 when the form is
restricted onto a hyperplane. In that case ${\cal U}_B$ is either
non-singular at the point $p$, or has a quadratic singularity of
rank $\geqslant r_1$. Finally, if $p\in E_b$ is a bi-quadratic
singularity, then the rank of every quadratic form in the pencil
$$
\lambda_1\widetilde{\beta_1}|_{{\mathbb
P}^k}+\lambda_2\widetilde{\beta_2}|_{{\mathbb P}^k}
$$
is at least $r_2$. Therefore, if $p\in {\cal U}_B$ is a singular
point, then it is either a quadratic singularity of rank
$\geqslant r_2-2\geqslant r_1$, or a bi-quadratic singularity of
rank $\geqslant r_2$, which is what we claimed. Q.E.D.

Note that the singularities of type $(r_1,r_2)$ satisfy the same
principle as the quadratic singularities: if for a non-singular
divisor ${\cal Z}\subset{\cal Y}$ the restriction ${\cal X}|_{\cal
Z}$ has singularities of type $(r_1,r_2)$, then the same is also
true for singularities of ${\cal X}$ at the points of intersection
${\cal X}\cap{\cal Z}$. We used this principle in the proof of the
last theorem.

If $r_1\geqslant 5$ and $r_2\geqslant 7$, then by the theorem
which we have just shown, the singularities of type $(r_1,r_2)$
are factorial and terminal, as we claimed in Subsection 0.1. For
the fibres of the fibre space $V/S$ in Theorem 0.1 we have
$(r_1,r_2)=(5,7)$, which guarantees factoriality and terminality
of singularities of every fibre, and by construction of the fibre
space $V/S$ as a sub-fibration of a locally trivial bundle
$\pi\colon X\to S$ with the fibre ${\mathbb P}$ also factoriality
and terminality of singularities of the variety $V$ itself. If
$\sigma_S\colon S^+\to S$ is a birational morphism, then
$V^+=V\mathop{\times}\nolimits_S S^+$ is a sub-fibration of the
locally trivial bundle ${\mathbb P}$-bundle
$$
\pi^+_X\colon X^+=X\mathop{\times}\nolimits_S S^+\to S^+
$$
and for that reason the variety $V^+$ is factorial and its
singularities are terminal. This proves that the Fano-Mori fibre
space $V/S$, constructed in Subsection 0.1, is stable with respect
to fibre-wise birational modifications. (The values
$(r_1,r_2)=(9,13)$ in Subsection 0.2 are needed for the proof of
divisorial canonicity in \S\S 2-3.)

Let us make a few more concluding remarks.

{\bf Remark 1.4.} If the base $S$ is one-dimensional, that is,
$S={\mathbb P}^1$, then the claim of Theorem 0.1 can be improved:
the birational map $\chi\colon V\dashrightarrow V'$ onto the total
space of the Mori fibre space $V'/S'$ is a biregular isomorphism,
see \cite[Theorem 1, (iv)]{Pukh2018b}; the proof is word for word
the same as in \cite[Subsection 1.5]{Pukh2018b}.

{\bf Remark 1.5.} The conditions (ii) and (iii) in Theorem 0.1
(and in Theorem 1.2) can be replaced by one stronger condition (as
it was done in \cite{Pukh15a}): the class
$$
-N(K_V\cdot \pi^{-1}(\overline{C}))-F
$$
is not effective for all $N\geqslant 1$. That condition is easier
to check (see Example 0.1). However, it seems that the
$K$-condition is close to a criterial one (in the
three-dimensional case see \cite{Grin03a,Grin06} for varieties
with a pencil of del Pezzo surfaces of degree 1) and is one of the
fundamental conditions for Fano-Mori fibre spaces.

{\bf Remark 1.6.} In \cite[Subsection 2.3]{Pukh15a} there is a misprint in
the displayed formula (2.2): the class $\pi^* Y$ is missing. This
does not affect the result of the computations as the missing
class is accounted for and the final result (at the end of
Subsection 2.3) is correct.

%%%%%%%%%%%%%%%%%%%%%%%%%%%%%%%%%%%%%%%%%%%%%%%%%%%%%%%%%%%%%%%%%%
%%%%%%%%%%%%%%%%%%%%%%%%%%%%%%%%%%%%%%%%%%%%%%%%%%%%%%%%%%%%%%%%%%
%%%%%%%%%%%%%%%%   SECTION 2

\section{Regular complete intersections}

In this section we prove Theorem 0.2 and obtain upper bounds for
multiplicities of singular points of certain subvarieties of
regular complete intersections $F\in {\cal F}_{\rm reg}$. In
Subsection 2.1 the proof of Theorem 0.2 is reduced to an estimate
for the codimension of the complement to the set of tuples of
homogeneous polynomials, defining an irreducible reduced complete
intersection of the corresponding codimension in the projective
space (Theorem 2.1). In Subsection 2.2 we prove that theorem. In
Subsecti\-on 2.3, using the technique of hypertangent divisors,
based on the regularity conditi\-ons, we obtain upper bounds for
the ratio of multiplicity to the degree for prime divisors on the
sections of the variety $F$ by linear subspaces in the projective
space ${\mathbb P}$.\vspace{0.3cm}

{\bf 2.1. Proof of Theorem 0.2.} Let $o\in{\mathbb P}$ be an
arbitrary point and ${\cal F}(o)\subset{\cal F}$ the closed set of
pairs $(f_1,f_2)$, such that $f_1(o)=0$, $f_2(o)=0$. Let us fix a
system of affine coordinates $z_1,\dots,z_{M+2}$ with the origin
at the point $o$. For a pair of homogeneous polynomials
$(f_1,f_2)\in{\cal F}(o)$ denote by the same symbols $f_1$ and
$f_2$ the corresponding non-homogeneous polynomials in the
variables $z_*$. Write down
$$
f_i(z_*)=f_{i,1}+f_{i,2}+\dots+f_{i,d_i},
$$
$i=1,2$, where $f_{i,j}(z_*)$ is a homogeneous polynomial of
degree $j$. Define the subsets ${\cal B}(?)\subset{\cal F}(o)$ for
$$
?\in\{1,2,2^2.1,2^2.2,2^2.3\}
$$
in the following way. The subset ${\cal B}(1)$ consists of pairs
$(f_1,f_2)$, such that the linear forms $f_{1,1}$ and $f_{2,1}$
are linearly independent, but the condition (R1) is not satisfied.
The subset ${\cal B}(2)$ consists of pairs $(f_1,f_2)\in{\cal
F}(o)$, such that
$$
\mathop{\rm dim}\langle f_{1,1},f_{2,1}\rangle=1,
$$
but the condition (R2) is not satisfied. Finally, the subsets
${\cal B}(2^2.*)$ consist of pairs $(f_1,f_2)\in{\cal F}(o)$, such
that
$$
f_{1,1}\equiv f_{2,1}\equiv 0,
$$
but the corresponding condition (R$2^2$.*) is not satisfied. Set
${\cal B}=\cup {\cal B}(?)$. Since the point $o\in{\mathbb P}$ is
arbitrary and $\mathop{\rm codim}({\cal F}(o)\subset{\cal F})=2$,
we get the obvious inequality
$$
\mathop{\rm codim}(({\cal F\backslash{\cal F}_{\rm
reg}})\subset{\cal F})\geqslant\mathop{\rm codim}({\cal B}\subset{\cal
F}(o))-M.
$$
Therefore, in order to prove Theorem 0.2, it is sufficient to
estimate from below the codimensions
$$
\mathop{\rm codim}({\cal B}(?)\subset{\cal F}(o))
$$
and choose the worst of these estimates. In order to get the
estimates we will use the following well known facts and methods:

--- the standard properties of the binomial coefficient,

--- the fact that the set of quadratic forms of rank $\leqslant r$
in $N$ variables is of codimension $\frac12(N-r)(N-r+1)$ in the
space of all quadratic forms,

--- the ``projection method'' of estimating the codimension
of the set of non-regular sequences, see \cite[Chapter 3, Section
1.3]{Pukh13a}.

Let us consider first the {\bf non-singular case}.

Let $\xi_1(z_*)$ and $\xi_2(z_*)$ be linearly independent linear
forms. Set
$$
{\cal F}(o,\xi_1,\xi_2)=\{(f_1,f_2)\in{\cal
F}(o)\,|\,f_{1,1}=\xi_1,f_{2,1}=\xi_2\}
$$
and ${\cal B}(1,\xi_1,\xi_2)={\cal B}(1)\cap{\cal
F}(o,\xi_1,\xi_2)$. Obviously,
$$
\mathop{\rm codim}({\cal B}(1)\subset{\cal F}(o))=\mathop{\rm
codim}({\cal B}(1,\xi_1,\xi_2)\subset{\cal F}(o,\xi_1,\xi_2)),
$$
so that we may assume that the linear forms $f_{1,1}$ and
$f_{2,1}$ are fixed. Let
$$
{\cal L}\subset\{\xi_1=\xi_2=0\}
$$
be a subspace of codimension 2 and
$$
{\cal B}(1,\xi_1,\xi_2,{\cal L})\subset{\cal B}(1,\xi_1,\xi_2)
$$
the subset, consisting of pairs $(f_1,f_2)$, such that the
sequence ${\cal S}[-5]|_{\cal L}$ is non-regular. Obviously the
codimension
$$
\mathop{\rm codim}({\cal B}(1,\xi_1,\xi_2) \subset{\cal
F}(o,\xi_1,\xi_2))
$$
is not less than
$$
\mathop{\rm codim}({\cal B}(1,\xi_1,\xi_2,{\cal L})\subset{\cal
F}(o,\xi_1,\xi_2))-2(M-2).
$$
Therefore, it is sufficient to estimate the codimension of the set
${\cal B}(1,\xi_1,\xi_2,{\cal L})$. Now let us apply the
``projection method'' (taking into account that $d_1\leqslant
d_2$). This method, introduced in \cite{Pukh98b}, estimates the
codimension of the set of non-regular sequences, fixing the first
moment when the regularity is violated: in that way, one obtains
an estimate for the set of sequences, in which the regularity is
violated {\it for the first time} in the $k$-th member of the
sequence; after that one has to choose the worst of all estimates.
For recent papers where this technique is described in detail, see
\cite[\S 3]{EvansPukh2}. In out case the codimension of the set
${\cal B}(1,\xi_1,\xi_2,{\cal L})$ is not less than the minimum of
the following set of integers, which it is convenient to break
into two parts. Assume first that $d_2\geqslant d_1+5$. Then the
first part of the set consists of integers
$$
{M-k\choose k},\quad k=2,\dots,d_1,
$$
and the second part --- of integers
$$
{M-d_1\choose d_1+k},\quad k=1,\dots,M-2d_1-3.
$$
From here, using the usual properties of binomial coefficients, it
is easy to conclude that the minimum of that set is attained at
one of the endpoints; comparing the numbers
$$
{M-2\choose 2}\quad \mbox{and}\quad {d_2-2\choose 3},
$$
we see by an elementary check that the first one is smaller than
the second. The same result is obtained in the omitted five cases,
when
$$
d_1\leqslant d_2\leqslant d_1+4.
$$
For instance, for $d_2=d_1$ in the first part of the set
$k=2,\dots, d_1-3$, and the second part consists of just one
number ${M-(d_1-3)\choose 3}$; for $d_2=d_1+1$ in the first part
we have $k=2,\dots,d_1-2$, and the second part is empty etc. The
minimum of the sequence is always
$$
{M-2\choose 2}=\frac12(M-2)(M-3).
$$
From here, taking into account the remarks above, we get the
estimate
$$
\mathop{\rm codim}({\cal B}(1)\subset{\cal
F}(o))\geqslant\frac12(M^2-9M+14).
$$

Now let us consider the {\bf quadratic case}.

Arguing as above, let us fix a pair of proportional linear forms
$\xi_1,\xi_2$, not both identically zero, and define the subsets
${\cal F}(o,\xi_1,\xi_2)$ and ${\cal B}(2,\xi_1,\xi_2)$ by the
same formulas as above in the non-singular case. Write down
$$
\xi_1=\alpha_1\tau\quad\mbox{and}\quad \xi_2=\alpha_2\tau,
$$
where $\tau(z_*)$ is a non-zero linear form and $(\alpha_1,
\alpha_2)\neq(0,0)$. Let ${\cal L}\subset\{\tau=0\}$ be a linear
subspace of codimension 1 and
$$
{\cal B}(2,\xi_1,\xi_2,{\cal L})\subset{\cal B}(2,\xi_1,\xi_2)
$$
the subset, consisting of pairs $(f_1,f_2)$, such that the
sequence ${\cal S}[-4]|_{\cal L}$ is non-regular. Similar to the
non-singular case, the set ${\cal B}(2,\xi_1,\xi_2)$ with respect
to ${\cal F}(o,\xi_1,\xi_2)$ is of codimension at least
$$
\mathop{\rm codim}({\cal B}(2,\xi_1,\xi_2,{\cal L}) \subset{\cal
F}(o,\xi_1,\xi_2))-M,
$$
so that it is sufficient to estimate the last codimension. We do
not give the elementary computations, based on the ``projection
method'': the codimension of the subset ${\cal B}(2)$ with respect
to ${\cal F}(o)$ turns out to be higher than the codimension of
the subset ${\cal B}(2^2.3)$, which we will estimate below.

Finally, let us consider the {\bf bi-quadratic case}.

The identical vanishing of the linear forms $f_{1,1}$ and
$f_{2.1}$ gives $2(M+2)$ independent conditions. Let us first
estimate the codimension of the set ${\cal B}(2^2.1)$, because the
condition (R$2^2$.1) does not depend on $d_1$. If
$$
\mathop{\rm rk}(f_{1,2},f_{2,2})\leqslant 12,
$$
then either $\mathop{\rm rk}f_{1,2}\leqslant 12$, or $\mathop{\rm
rk}f_{1,2}\geqslant 13$ and $f_{2,2}$ lies on the cone in the
space ${\cal P}_{2,M+2}$, the vertex of which is the form
$f_{1,2}$, and the base of which is the closed subset of forms of
rank $\leqslant 12$, whence we get:
$$
\mathop{\rm codim}({\cal B}(2^2.1)\subset{\cal
F}(o))\geqslant\frac12(M-9)(M-10)-1+2(M+2),
$$
because the conditions $\mathop{\rm rk}f_{i,2}\leqslant 17$,
$i=1,2$, together give $(M-14)(M-15)$ independent conditions for
the pair of quadratic forms $(f_{1,2},f_{2,2})$, which is much
higher.

Now let us consider the conditions (R$2^2$.2) and (R$2^2$.3).
Assume first that $d_1\geqslant 4$. Let us estimate the
codimension of the set ${\cal B}(2^2.3)$. Let ${\cal
L}\subset{\mathbb C}^{M+2}$ be a linear subspace of codimension 2
or 3. Let ${\cal B}(2^2.3,{\cal L})\subset{\cal F}(o)$ be the
subset, consisting of pairs $(f_1,f_2)$, such that $f_{1,1}\equiv
f_{2,1}\equiv 0$ and the sequence
$$
{\cal S}[-\mathop{\rm codim}{\cal L}-1]|_{\cal L}
$$
is not regular. The case $\mathop{\rm codim}{\cal L}=3$ gives the
worst estimate for the codimension, so we will consider that case.
Applying the ``projection method'', in the same way as in the
non-singular case, we get the inequality
$$
\mathop{\rm codim}({\cal B}(2^2.3,{\cal L})\subset{\cal F}(o))
\geqslant{M-1\choose 2}+2(M+2),
$$
which, taking into account that ${\cal L}$ varies in the
$3(M-1)$-dimensional Grassmanian, implies the following estimate:
$$
\mathop{\rm codim}({\cal B}(2^2.3)\subset{\cal F}(o))
\geqslant\frac12(M^2-5M+16).
$$

Now let us consider the condition (R$2^2$.2). Here we will need
the following general fact. Let $\underline{m}=(m_1,\dots,m_k)$ be
a tuple of integers, where
$$
2\leqslant m_1\leqslant m_2\leqslant\dots\leqslant m_k,
$$
and
$$
{\cal P}(\underline m)=\prod^k_{i=1}{\cal P}_{m_i,N+1}
$$
the space of tuples $(\underline{g})=(g_1,\dots,g_k)$ of
homogeneous polynomials of degree $m_1,\dots,m_k$ in $N+1$
variables, which we consider as homogeneous polynomials on the
projective space ${\mathbb P}^N$. Let
$$
{\cal B}^*(\underline{m})\subset{\cal P}(\underline{m})
$$
be the set of tuples $(g_1,\dots,g_k)$, such that the scheme of
their common zeros
$$
V(\underline{g})=V(g_1,\dots,g_k)
$$
is not an irreducible reduced subvariety of codimension $k$ in
${\mathbb P}^N$.

{\bf Theorem 2.1.} {\it The following inequality holds:}
$$
\mathop{\rm codim} ({\cal B}^*(\underline{m})\subset{\cal P}(\underline{m}))
\geqslant\frac12(N-k-1)(N-k-4)+2.
$$

{\bf Proof} is given in Subsection 2.2.

Applying Theorem 2.1, we get the following inequality for the
codimension of the subset ${\cal B}(2^2.2,{\cal L})\subset{\cal
F}(o)$, where ${\cal L}$ is a linear subspace of codimension 2 in
${\mathbb C}^{M+2}$, consisting of pairs $(f_1,f_2)$, such that
$f_{1,1}\equiv f_{2.1}\equiv 0$ and the condition (R$2^2$.2) is
not satisfied:
$$
\mathop{\rm codim}({\cal B}(2^2.2,{\cal L})\subset{\cal F}(o))
\geqslant\frac12(M^2-15M+58)+2(M+2),
$$
so that, taking into account that the subspace ${\cal L}$ is an
arbitrary one, we get
$$
\mathop{\rm codim}({\cal B}(2^2.2)\subset{\cal F}(o))
\geqslant\frac12(M^2-15M+66).
$$

Choosing among the estimates for the codimension, obtained above,
the worst one (it corresponds to violation of the condition
(R$2^2$.2)), and recalling that the closed subset ${\cal F}(o)$
has codimension 2 and the point $o$ varies in the
$(M+2)$-dimensional projective space, we complete the proof of
Theorem 0.2 for $d_1\geqslant 4$.

Let us consider now the two remaining cases $d_1=2,3$. The
codimension of the set ${\cal B}(2^2.3)$ is estimated by the
projection method, as above. If $\mathop{\rm codim}{\cal L}=2$,
then the worst estimate of the codimension
$$
\mathop{\rm codim}({\cal B}(2^2.3,{\cal L})\subset{\cal F}(o))
$$
corresponds to the violation of the regularity by the last member
of the sequence ${\cal S}[-\mathop{\rm codim}{\cal L}]|_{\cal L}$
and is given by the number
$$
{d_2-2+2\choose 2}+2(M+2),
$$
which implies that
$$
\mathop{\rm codim}({\cal B}(2^2.3)\subset{\cal F}(o)) \geqslant
{d_2\choose 2}-M+7
$$
(where $d_2=M+2-d_1\geqslant M-1$). This is better than what
Theorem 0.2 claims.

In order to estimate the codimension of the set ${\cal B}(2^2.2)$,
we use Theorem 2.1.

Assume first that $d_1=3$. This case is special because the
condition (R$2^2$.2) requires irreducibility and reducedness of a
complete intersection of codimension 5 (in all other cases this
codimension is equal to 4). Using Theorem 2.1 for a fixed linear
subspace ${\cal L}\subset{\mathbb C}^{M+2}$, we obtain the
inequality
$$
\mathop{\rm codim}({\cal B}(2^2.2,{\cal L})\subset{\cal F}(o))
\geqslant\frac12(M^2-17M+74)+2(M+2).
$$
From this inequality, taking into account that the subspace ${\cal
L}$ is arbitrary, we obtain the estimate
$$
\mathop{\rm codim}({\cal B}(2^2.2)\subset{\cal F}(o))
\geqslant\frac12(M^2-17M+82).
$$
Taking into account that the closed subset ${\cal F}(o)$ is of
codimension 2, and the point $o$ varies in ${\mathbb P}^{M+2}$, we
complete the proof of Theorem 0.2 for $d_1=3$.

In the case $d_1=2$ the codimension of the set ${\cal B}(2^2.2)$
is estimated in the word for word the same way as for
$d_1\geqslant 4$, since in the condition (R$2^2$.2) we have a
complete intersection of codimension 4. Proof of Theorem 0.2 is
complete.\vspace{0.3cm}

%%%%%%%%%%%%%%%%%%%%%%%%%%%%%%%%%%%%%%%%%%%%%%%%%%%%%%%%%%%%%%%
%%%%%%%%%%%%%%%%%%%%%   subsection 2.2

{\bf 2.2. Irreducible reduced complete intersections.} Let us
prove Theorem 2.1. Our arguments are similar to the arguments of
\cite[\S 2]{EvansPukh2}, but are more general and more formal. The
idea of the proof is to use the induction on $k$. The case $k=1$
is almost obvious, the codimension of the set of reducible
polynomials is very easy to estimate. However, in order to make
the step of induction, one has to consider such tuples of
polynomials that define not only irreducible reduced subvarieties
of the required codimension, but also factorial subvarieties,
which allows to add one more polynomial. Our notations are close
to the notations in \cite[\S 2]{EvansPukh2}, whenever it is
possible. Now, let
$$
{\cal P}^{\geqslant j}=\prod^k_{i=j}{\cal P}_{m_i,N+1}
$$
be the space of truncated tuples $(g_j,\dots,g_k)$, which we will
denote by the symbol $g_{[j,k]}$. Instead of the increasing
induction on $k$ we use the equivalent decreasing induction on
$j=k,\dots,1$. It is more convenient if the degree of the added
polynomial is not higher than the degrees of the polynomials which
are already in the tuple. By the symbol ${\cal P}^{\geqslant
j}_{\rm mq}$ we denote the Zariski open subset of the space ${\cal
P}^{\geqslant j}$, consisting of tuples $g_{[j,k]}$, such that the
scheme of their common zeros
$$
V(g_{[j,k]})\subset{\mathbb P}^N
$$
is an irreducible reduced complete intersection of codimension
$k-j+1$ and at most factorial mutli-quadratic singularities (see
\cite[\S 2]{EvansPukh2}). Recall the meaning of the last
condition. Let $o\in V(g_{[j,k]})$ be an arbitrary point and
$w_1,\dots,w_N$ a system of affine coordinates on ${\cal P}^N$
with the origin at that point. Let
$$
g_a=g_{a,1}+g_{a,2}+\dots+g_{a,m_a}
$$
be the decomposition of the non-homogeneous polynomial $g_a$
(which is denoted by the same symbol) into components, homogeneous
in $w_*$. If
$$
\mathop{\rm dim}\langle g_{j,1},\dots,g_{k,1}\rangle=k-j+1,
$$
then in a neighborhood of the point $o$ the scheme $V(g_{[j,k]})$
is a non-singular subvariety of codimension $k-j+1$. Assume that
$$
l=k-j+1-\mathop{\rm dim}\langle
g_{j,1},\dots,g_{k,1}\rangle\geqslant 1.
$$
Let $\varphi_{\mathbb P,o}\colon X\to{\mathbb P}^N$ be the blow up
of the point $o$ with the exceptional divisor $E_X\cong{\mathbb
P}^{N-1}$. Let us break the set of indices $\{j,\dots,k\}$ into
two disjoint subsets $I_1\sqcup I_2$, such that the linear forms
$g_{\alpha,1}$, $\alpha\in I_1$, are linearly independent, and
moreover,
$$
\langle g_{j,1},\dots,g_{k,1}\rangle=\langle
g_{\alpha,1}\,|\,\alpha\in I_1\rangle,
$$
so that for $\gamma\in I_2$ there exist uniquely determined
constants $c_{\gamma\alpha}$, such that
$$
g_{\gamma,1}=\sum_{\alpha\in I_1}c_{\gamma\alpha}g_{\alpha,1}.
$$
Set $\Pi=\{g_{\alpha,1}=0\,|\,\alpha\in I_1\}\subset{\mathbb C}^N$
and construct for $\gamma\in I_2$ the following quadratic forms:
$$
g^*_{\gamma,2}=g_{\gamma,2}-\sum_{\alpha\in
I_1}c_{\gamma\alpha}g_{\alpha,2}.
$$
If the inequality
\begin{equation}\label{15.06.2020.1}
\mathop{\rm rk}(g^*_{\gamma,2}|_{\Pi},\gamma\in I_2) \geqslant
2l+3
\end{equation}
holds, then in a neighborhood of the point $o$ the scheme
$V(g_{[j,k]})$ is an irreducible complete intersection with a
mutli-quadratic singularity of type $2^l$, see \cite{EvansPukh2}:
let $V^+(g_{[j,k]})$ be the strict transform of $V(g_{[j,k]})$ on
$X$ and $Q=V^+(g_{[j,k]})\cap E_X$ the exceptional divisor of the
blow up $V^+(g_{[j,k]})\to V(g_{[j,k]})$ of the point $o$, then
$Q\subset E_X$ is an irreducible reduced complete intersection of
type $2^l$ in the linear subspace
$$
{\mathbb P}(\Pi)=\{g_{\alpha,1}=0 \,|\, \alpha\in I_1\}\subset E_X
$$
of codimension $(k-j+1-l)$, where by \cite[\S 2, Lemma
2.1]{EvansPukh2} the inequality
$$
\mathop{\rm codim}(\mathop{\rm Sing}Q\subset Q)\geqslant 4
$$
holds. If the inequality (\ref{15.06.2020.1}) is true for any
singular point $o\in V(g_{[j,k]})$, then $V(g_{[j,k]})$ is an
irreducible reduced factorial variety of codimension $k-j+1$.

We will prove Theorem 2.1, estimating the codimension of the
complement
$$
\mathop{\rm codim}(({\cal P}^{\geqslant j}\backslash{\cal
P}^{\geqslant j}_{\rm mq})\subset{\cal P}^{\geqslant j})
$$
for $j=k,\dots,1$. So we argue by decreasing induction on $j$.

{\bf The basis of induction.} The closed set of reducible
polynomials of degree $m_k$ on ${\mathbb P}^N$ has the codimension
$$
{N+m_k-1\choose m_k}-N
$$
in the space ${\cal P}_{m_k,N+1}$. If $g_k\in{\cal P}_{m_k,N+1}$
is an irreducible polynomial, then the condition that the
hypersurface $\{g_k=0\}$ has at least one singular point which is
not a quadratic singularity of rank $\geqslant 5$ (that is, a
factorial quadratic singularity), imposes on the coefficients of
$g_k$ at least
$$
\frac12(N-3)(N-4)+1
$$
independent conditions. Since $m_k\geqslant 2$, we get
$$
\mathop{\rm codim}(({\cal P}^{\geqslant k}\backslash{\cal
P}^{\geqslant k}_{\rm mq})\subset{\cal P}^{\geqslant k})
\geqslant\frac12(N-3)(N-4)+1.
$$

{\bf Step of induction: irreducibility.} Let
$$
g_{[j+1,k]}\in{\cal P}^{\geqslant j+1}_{\rm mq}.
$$
The arguments, similar to \cite[Subsection 2.2]{EvansPukh2}, show
that the set of polynomials $g_j\in{\cal P}_{m_j,N+1}$, such that
$V(g_{[j,k]})$ is either reducible, or non-reduced, or is of
``incorrect'' codimension $(k-j)$ in ${\mathbb P}^N$, has in
${\cal P}_{m_j,N+1}$ the codimension at least
\begin{equation}\label{16.06.2020.2}
{N+m_j-1\choose m_j}-N-(k-j).
\end{equation}
(This follows from the factoriality of the variety $V(g_{[j,k]})$
and the Lefschetz theorem: the Picard group of this variety is
generated by the hyperplane section and every divisor is cut out
by a hypersurface in ${\mathbb P}^N$. It is here that we use the
condition
$$
2\leqslant m_1\leqslant m_2\leqslant\dots\leqslant m_k.)
$$
Now we assume that the scheme $V(g_{[j,k]})$ is an irreducible
reduced complete intersection of codimension $(k-j+1)$ in
${\mathbb P}^N$.

{\bf Step of induction: factoriality.} Let us estimate the
codimension of the set of tuples $g_{[j,k]}$, such that the
irreducible reduced variety $V(g_{[j,k]})$ does not satisfy the
condition to have factorial multi-quadratic singularities. It is
easier to do it from zero, that is to say, not assuming that the
variety $V(g_{[j+1,k]})$ has factorial multi-quadratic
singularities. Let $o\in V(g_{[j,k]})$ be a point. By the symbol
${\cal Q}(l)$ we denote the space of quadratic forms on the linear
space $\Pi$. The dimension of the latter is $\mathop{\rm
dim}\Pi=N+j-k+l-1$. Since $|I_2|=l$, instead of tuples of
quadratic forms $(g^*_{\gamma,2}|_{\Pi},\gamma\in I_2)$, varying
independently of each other, we will consider, in order to
simplify the formulas, the tuples
$$
h_{[1,l]}=(h_1,\dots,h_l)\in{\cal Q}(l)^{\times l}.
$$
By the symbol ${\cal R}_{\leqslant a}$ we denote the closed subset
of quadratic forms of rank $\leqslant a$ in ${\cal Q}(l)$,
$$
\mathop{\rm codim}({\cal R}_{\leqslant a}\subset{\cal Q}(l))=
\frac12(N+j-k+l-a)(N+j-k+l-a-1),
$$
and for $e\in\{1,\dots,l\}$ let ${\cal X}_{e,a}\subset{\cal
Q}(l)^{\times e}$ be the closed set of tuples
$h_{[1,e]}=(h_1,\dots,h_e)$, such that
$$
\mathop{\rm rk}h_{[1,e]}\leqslant a.
$$
The following claim is true.

{\bf Lemma 2.1.} {\it The following inequality holds:}
$$
\mathop{\rm codim}({\cal X}_{e,a}\subset{\cal Q}(l)^{\times e})
\geqslant\mathop{\rm codim}({\cal R}_{\leqslant a}\subset{\cal
Q}(l))-e+1.
$$

{\bf Proof} is identical to the proof of Lemma 2.2 in \cite[\S
2]{EvansPukh2}.

In particular, substituting $e=l$ and $a=2l+2$, we get that the
codimension of the closed subset ${\cal X}_{l,2l+2}$ in the space
${\cal Q}(l)^{\times l}$ (that is, precisely of the subset of
tuples of quadratic forms $g^*_{\gamma,2}|_{\Pi},\gamma\in I_2$,
which do not satisfy the inequality (\ref{15.06.2020.1})) is at
least
\begin{equation}\label{16.06.2020.1}
\frac12(N+j-k-l-2)(N+j-k-l-3)-l+1.
\end{equation}
Coming back to the complete intersection $V(g_{[j,k]})$, which we
assume to be irreducible and reduced, let us estimate the
codimension of the set of tuples
$$
g_{[j,k]}\not\in{\cal P}^{\geqslant j}_{\rm mq}.
$$
Let $o\in{\mathbb P}^N$ be a fixed point and
$$
l\in\{1,\dots,k-j+1\}
$$
an integer. The of tuples $g_{[j,k]}$, such that $o\in
V(g_{[j,k]})$ and
$$
\mathop{\rm dim}\langle g_{j,1},\dots, g_{k,1}\rangle=k-j+1-l,
$$
is of codimension
$$
(k-j+1)+l(N+j-k-1+l)
$$
in ${\cal P}^{\geqslant j}$. Violation of the condition
(\ref{15.06.2020.1}) gives the additional codimension
(\ref{16.06.2020.1}). Taking into account that the point $o$
varies in ${\mathbb P}^N$, but keeping the value $l$ fixed, we
obtain a quadratic function in $l$, the minimum of which is
attained at $l=1$ and is equal to
\begin{equation}\label{16.06.2020.3}
\frac12(N+j-k-2)(N+j-k-5)+2.
\end{equation}
This expression estimate from below the codimension of the
complement ${\cal P}^{\geqslant j}\backslash{\cal P}^{\geqslant
j}_{\rm mq}$ in ${\cal P}^{\geqslant j}$, since it is certainly
smaller than (\ref{16.06.2020.2}). Finally, considering
(\ref{16.06.2020.3}) as a function of the variable
$j\in\{1,\dots,k\}$, we see that its minimum is attained at $j=1$
and is equal to
$$
\frac12(N-k-1)(N-k-4)+2.
$$
This completes the proof of Theorem 2.1. Q.E.D.

{\bf Remark 2.1.} The regularity condition (R$2^2$.2) is about a
complete intersection in the affine space ${\cal L}$, but not in
the projective space. Let us explain, how Theorem 2.1 applies to
the affine situation. Consider ${\cal L}$ as an affine chart
${\mathbb C}^M\subset{\mathbb P}^M$. Obviously, the equations
$f_{1,2},f_{2,2}$ (which are homogeneous by construction) and
$f_1,f_2$ (we come back to the original homogeneous polynomials on
${\mathbb P}$), restricted onto ${\mathbb P}^M$, define an
irreducible reduced complete intersection of codimension 4 in
${\mathbb P}^M$ if and only if the condition (R$2^2$.2) is
satisfied. If the same is true for the restrictions of those four
polynomials onto the ``hyperplane at infinity'' ${\mathbb
P}^M\backslash{\mathbb C}^M$, then the condition (R$2^2$.2) holds.
That hyperplane can be identified with the projectivization
${\mathbb P}({\cal L})$, and then we get a set of four homogeneous
polynomials
$$
f_{1,2}|_{\cal L},\, f_{2,2}|_{\cal L},\, f_{1,d_1}|_{\cal L},\,
f_{2,d_2}|_{\cal L}
$$
of degrees $2,2,d_1,d_2$, respectively. If these four polynomials
define an irreducible reduced complete intersection of
codimension4 in ${\mathbb P}({\cal L})\cong{\mathbb P}^{M-1}$,
then the condition (R$2^2$.2) is satisfied. Since the coefficients
of these four polynomials are four disjoint groups of coefficients
of the original polynomials $f_1,f_2$, Theorem 2.1 can be applied
to obtaining an estimate of the codimension of the subset ${\cal
B}(2^2.2)$.\vspace{0.3cm}

%%%%%%%%%%%%%%%%%%%%%%%%%%%%%%%%%%%%%%%%%%%%%%%%%%%%%%%%%%%%%%%%%%%%
%%%%%%%%%%%%%%%%%%%%%   subsection 2.3

{\bf 2.3. Degrees and multiplicities.} From the regularity
conditions we obtain in the usual way (see \cite[Chapter
3]{Pukh13a}) estimates for the multiplicity of singular points of
subvarieties $Y\subset F\in{\cal F}_{\rm reg}$. From this moment,
we fix a regular complete intersection $F\in{\cal F}_{\rm reg}$.
As usual, we denote the ratio of the multiplicity $\mathop{\rm
mult}_oY$ to the degree $\mathop{\rm deg}Y$ (with respect to the
embedding $Y\subset F\subset{\mathbb P}$) by the symbol
$$
\frac{\mathop{\rm mult}_o}{\mathop{\rm deg}}\, Y.
$$
Let us consider first the points that are non-singular on the
complete intersection $F$.

{\bf Proposition 2.1.} {\it Assume that $\Delta\subset F$ is a
section of the variety $F$ by an arbitrary linear subspace of
codimension 2 in ${\mathbb P}$, the point $o\in\Delta$ is
non-singular and $Y\subset\Delta$ is a prime divisor. Then the
following inequality holds:}
\begin{equation}\label{23.06.2020.1}
\frac{\mathop{\rm mult}_o}{\mathop{\rm deg}}\, Y
\leqslant\frac{2}{d_1d_2}.
\end{equation}

{\bf Proof.} Note that the inequality (\ref{23.06.2020.1}) is
optimal: the section of $\Delta$ by any hyperplane, which is
tangent to $\Delta$ at the point $o$, gives an equality. Let us
assume the converse: the inequality
$$
\frac{\mathop{\rm mult}_o}{\mathop{\rm deg}}\, Y>\frac{2}{d_1d_2}
$$
holds, and show that this assumption leads to a contradiction. Our
arguments are based on the technique of hypertangent divisors and
completely similar to the arguments in \cite[Chapter 3, Section
2]{Pukh13a}. We discuss in detail only those fragments of our
arguments, which need to be modified. By the conditions (R2) and
(R$2^2$.1), we have the inequality $\mathop{\rm codim}(\mathop{\rm
Sing}F\subset F)\geqslant 10$, so that for the subvariety
$\Delta\subset F$ we get:
$$
\mathop{\rm codim}(\mathop{\rm Sing}\Delta\subset\Delta) \geqslant
6.
$$
Therefore, $\Delta\subset{\mathbb P}^M$ is a factorial complete
intersection of codimension 2 and $\mathop{\rm Pic}\Delta={\mathbb
Z}H_{\Delta}$, where $H_{\Delta}$ is the class of a hyperplane
section. By the symbol $|H_{\Delta}-2o|$ we denote the pencil of
tangent hyperplanes at the point $o$. Let
$D_{1,1}\in|H_{\Delta}-2o|$ be an arbitrary divisor. By the
condition (R1) we have $\mathop{\rm mult}_oD_{1,1}=2$, so that
$D_{1,1}\neq Y$ and for that reason the algebraic cycle
$(D_{1,1}\circ Y)$ of the scheme-theoretic intersection of
$D_{1,1}$ and $Y$ is well defined. This is an effective cycle of
codimension 2 on $\Delta$, and moreover,
$$
\mathop{\rm mult}\nolimits_o(D_{1,1}\circ Y) \geqslant
2\mathop{\rm mult}\nolimits_oY
$$
and $\mathop{\rm deg}(D_{1,1}\circ Y)=\mathop{\rm deg}Y$, so that
there is an irreducible component $Y_2$ of that cycle (setting
$Y_1=Y$), which satisfies the inequality
$$
\frac{\mathop{\rm mult}_o}{\mathop{\rm deg}}\,
Y_2>\frac{4}{d_1d_2}.
$$
Let $P\ni o$ be a general 7-dimensional subspace in ${\mathbb
P}^M$, containing the point $o$. Then $\Delta_P=\Delta\cap P$ is a
non-singular 5-dimensional complete intersection in
$P\cong{\mathbb P}^7$, so that for the numerical Chow group of
classes of cycles of codimension 2 on $\Delta_P$ we have:
$$
A^2\Delta_P={\mathbb Z}H^2_P,
$$
where $H_P\in\mathop{\rm Pic}\Delta_P$ is the class of a
hyperplane section. Let $D_{2,1}\neq D_{1,1}$ be another element
of the tangent pencil, then we have
$$
(D_{1,1}\circ D_{2,1}\circ\Delta_P)\sim H^2_P,
$$
so that $(D_{1,1}\circ D_{2,1}\circ\Delta_P)= D_{1,1}\cap
D_{2,1}\cap\Delta_P$ is an irreducible subvariety, and for that
reason $(D_{1,1}\circ D_{2,1})=D_{1,1}\cap D_{2,1}$ is also an
irreducible subvariety of codimension 2 on $\Delta$, which by the
condition (R1) satisfies the equality
$$
\mathop{\rm mult}\nolimits_o(D_{1,1}\circ D_{2,1})=4
$$
and for that reason $Y_2\neq (D_{1,1}\circ D_{2,1})$. Since by
construction $Y_2\subset D_{1,1}$, we have $Y_2\not\subset
D_{2,1}$, so that the effective cycle $(D_{2,1}\circ Y_2)$ of the
scheme-theoretic intersection of $D_{2,1}$ and $Y_2$ is well
defined. There is an irreducible component $Y_3$ of that cycle,
satisfying the inequality
$$
\frac{\mathop{\rm mult}\nolimits_o}{\mathop{\rm deg}}\,
Y_3>\frac{8}{d_1d_2}.
$$
Now we apply the technique of hypertangent divisors in the word
for word the same way as in \cite[Chapter 3]{Pukh13a}: define the
hypertangent linear systems in terms of a system of affine
coordinates $z_1,\dots, z_{M+2}$
$$
\Lambda_j=\left|\sum^{\mathop{\rm
min}\{j,d_1-1\}}_{a=1}f_{1,[1,a]}s_{1,j-a}+\sum^j_{a=
1}f_{2,[1,a]}s_{2,j-a}\right|,
$$
where $j=2,\dots,d_2-1$ and the symbol $f_{i,[1,a]}$ stands for
the left segment of length $a$ of the (non-homogeneous) polynomial
$f_i$,
$$
f_{1,[1,a]}=f_{i,1}+\dots+f_{i,a},
$$
and the homogeneous polynomials $s_{i,j-a}(z_*)\in{\cal
P}_{j-a,M+2}$ are arbitrary and independent of each other. If
$j\leqslant d_1-1$, then we choose in the linear system
$\Lambda_j|_{\Delta}$ two general divisors $D_{1,j}$ and
$D_{2,j}$; if $j\geqslant d_1$, then we choose in the linear
system $\Lambda_j|_{\Delta}$ one general divisor $D_{2,j}$. We put
the constructed
$$
d_1+d_2-4=M-2
$$
divisors on $\Delta$ in the lexicographic order (as it was done in
Subsection 0.2) and remove from this sequence the very first and 5
last divisors. We constructed a sequence consisting of $M-8$
hypertangent divisors, which it is convenient to denote,
respectively, by
$$
Z_3,Z_4,\dots,Z_{M-6}.
$$
Now, starting with already constructed subvariety $Y_3$, we
construct by induction a sequence of irreducible subvarieties
$$
Y_3,Y_4,\dots,Y_{M-6},Y_{M-5}
$$
of codimension $\mathop{\rm codim}(Y_k\subset\Delta)=k$ in the
following way: if the subvariety $Y_k$ is already constructed,
then the condition (R1) guarantees that $Y_k\not\subset|Z_k|$, so
that the algebraic cycle of the scheme-theoretic intersection
$(Z_k\circ Y_k)$is well defined and has codimension $k+1$, and we
choose for $Y_{k+1}$ such component of that cycle which has the
maximal value of the ratio $\mathop{\rm mult}_o\slash\mathop{\rm
deg}$. We do not give more details since this technique was used
many times and is well known. The subvariety $Y_{M-2}$ is
three-dimensional and for it the ratio $\mathop{\rm
mult}_o\slash\mathop{\rm deg}$ is strictly higher than:
$$
\frac43\frac{(d_1-2)(d_2-3)}{d_1d_2}\geqslant 1
\quad\mbox{for}\quad d_1=d_2\quad\mbox{or}\quad d_1=d_2-1,
$$
$$
\frac43\frac{(d_1-1)(d_2-4)}{d_1d_2}\geqslant
1\quad\mbox{for}\quad d_1=d_2-2\quad\mbox{or}\quad d_1=d_2-3,
$$
$$
\frac43\frac{d_2-5}{d_2}\geqslant 1\quad\mbox{for}\quad
d_1\leqslant d_2-4.
$$
In each case we obtain a contradiction which completes the proof
of the proposition.

Now let us consider the quadratic points of the complete
intersection $F$ and its hyperplane sections.

{\bf Proposition 2.2.} {\it Assume that $o\in F$ is a quadratic
singularity of the variety $F$ and $\Delta\subset F$ is a
hyperplane section of the variety $F$, containing the point $o$,
and $o\in\Delta$ is a quadratic singularity of the variety
$\Delta$. Let $Y\subset\Delta$ be a prime divisor. Then the
following inequality holds:}
\begin{equation}\label{24.06.2020.1}
\frac{\mathop{\rm mult}_o}{\mathop{\rm deg}}\, Y\leqslant
\frac{4}{d_1d_2}.
\end{equation}

{\bf Proof.} Since the point $o$ is a quadratic singularity of the
variety $F$ and its hyperplane section $\Delta$, in the notations
of Subsection 0.2 for the linear span
$\langle\Delta\rangle\subset{\mathbb P}$, that is, for the
hyperplane that cuts out $\Delta$, we have
$\langle\Delta\rangle\neq\{\tau=0\}$. Note that there is a
uniquely determined hyperplane section of the variety $\Delta$ (by
the hyperplane $\tau|_{\langle\Delta\rangle}=0$), for which the
inequality (\ref{24.06.2020.1}) becomes an equality, that is, the
claim of our proposition is optimal. Now let us assume the
converse: for some prime divisor $Y\subset\Delta$ the inequality
$$
\frac{\mathop{\rm mult}_o}{\mathop{\rm deg}}\, Y>\frac{4}{d_1d_2}
$$
holds. Let us show that this assumption leads to a contradiction.
Since the rank of a quadratic form drops by at most 2 when the
form is restricted onto a hyperplane, the condition (R2) implies
that the rank of the quadratic singularity $o$ of the variety
$\Delta$ is at least 7. Let $D_1$ be the hyperplane section of the
variety $\Delta$, described above, for which the inequality
(\ref{24.06.2020.1}) becomes an equality. We have $D_1\neq Y$, so
that the algebraic cycle $(D_1\circ Y)$ of the scheme-theoretic
intersection is well defined, and in that cycle we can find an
irreducible component $Y_2$ of codimension 2 with respect to
$\Delta$, satisfying the inequality
$$
\frac{\mathop{\rm mult}_o}{\mathop{\rm deg}}\,
Y_2>\frac{8}{d_1d_2}.
$$
Now, we apply the technique of hypertangent divisors in the word
for word the same way as in the non-singular case, using the
regularity condition (R2). We obtain a subvariety of positive
dimension in ${\mathbb P}$, the multiplicity of which at the point
$o$ is higher than its degree. This contradiction completes the
proof of the proposition. Q.E.D.

Now let us consider bi-quadratic points of the complete
intersection $F$. Fix an arbitrary bi-quadratic point $o\in F$.

{\bf Proposition 2.3.} {\it Assume that $\Delta\subset F$ is the
section of the variety $F$ by an arbitrary linear subspace of
codimension 2 or 3, containing the point $o$, and $Y\subset\Delta$
is a prime divisor. If $\mathop{\rm codim}(\Delta\subset F)=2$,
then the inequality
$$
\frac{\mathop{\rm mult}_o}{\mathop{\rm deg}}\,
Y\leqslant\frac{6}{d_1d_2}
$$
holds. If $\mathop{\rm codim}(\Delta\subset F)=3$, then the
following inequality holds:}
$$
\frac{\mathop{\rm mult}_o}{\mathop{\rm deg}}\,
Y\leqslant\frac{8}{d_1d_2}.
$$

{\bf Proof} is similar to the non-singular and quadratic cases: we
assume that the prime divisor $Y$ does not satisfy the
corresponding inequality, after which use the technique of
hypertangent divisors, based on the regularity conditions
(R$2^2$.2) and (R$2^2$.3), to construct, by means of successive
intersections, a subvariety of positive dimension in ${\mathbb
P}$, that gives a contradiction. In this procedure only its
beginning, making use of the condition (R$2^2$.2), is
non-standard, and we will consider these first steps in detail.
Assume that $d_1\geqslant 4$, and that $\Delta$ is the section of
the variety $F$ by a linear subspace
$\langle\Delta\rangle\subset{\mathbb P}$ of codimension 2. The
hypertangent system $\Lambda_2$ is the pencil of quadrics,
generated by the quadratic forms $f_{1,2}$ and $f_{2,2}$. The
hypertangent divisor $D_{1,2}=\{f_{1,2}|_{\Delta}=0\}$ is
irreducible (it can not be a sum of hyperplane sections by the
condition (R$2^2$.1) or (R$2^2$.3)), and moreover, the equalities
$$
\mathop{\rm mult}\nolimits_oD_{1,2}=12\quad\mbox{and}\quad
\mathop{\rm deg}D_{1,2}=2d_1d_2
$$
hold, so that $D_{1,2}\neq Y$ and the scheme-theoretic
intersection $(D_{1,2}\circ Y)$ is well defined. Choosing in the
latter cycle an irreducible component $Y_2$ (of codimension 2 on
$\Delta$) with the maximal ratio of the multiplicity at the point
$o$ to the degree, we get
$$
\frac{\mathop{\rm mult}_o}{\mathop{\rm deg}}\,
Y_2>\frac{9}{d_1d_2}.
$$
Let $D_{2,2}=\{f_{2,2}|_{\Delta}=0\}$ be the second hypertangent
divisor. The condition (R$2^2$.2) implies that $(D_{1,2}\circ
D_{2,2})=D_{1,2}\cap D_{2,2}$ is an irreducible subvariety of
codimension 2 on $\Delta$. By the condition (R$2^2$.1) (or
(R$2^2$.3)) we have the equality
$$
\frac{\mathop{\rm mult}_o}{\mathop{\rm deg}}\, (D_{1,2}\circ
D_{2,2})=\frac32\cdot\frac32\cdot\frac{4}{d_1d_2}=\frac{9}{d_1d_2},
$$
which implies that $Y_2\neq(D_{1,2}\circ D_{2,2}$). However,
$Y_2\subset|D_{1,2}|$ by construction, so that $Y_2\not\subset
|D_{2,2}|$ and for that reason the effective algebraic cycle of
the scheme-theoretic intersection $(D_{2,2}\circ Y_2)$ is a well
defined cycle of codimension 3 on $\Delta$. Its component $Y_3$
with the maximal value of the ratio $\mathop{\rm
mult}_o\slash\mathop{\rm deg}$ satisfies the inequality
$$
\frac{\mathop{\rm mult}_o}{\mathop{\rm deg}}\,
Y_3>\frac{27}{2d_1d_2}.
$$
The remaining part of the procedure of applying the technique of
hypertangent divisors is absolutely standard and we skip it: we
use only the condition (R$2^2$.3) and in the hypertangent system
$\Lambda_3$ choose one general divisor, not two, after which we
continue as in the non-singular or quadratic case, completing the
proof.

If $d_1=3$, then the modified regularity condition (R$2^2$.2)
makes it possible to argue as above, intersecting $Y=Y_1$
successively with hypertangent divisors
$$
D_{1,2}=\{f_{1,2}|_{\Delta}=0\},\quad D_{2,2}=\{f_{2,2}|_{\Delta}
=0\},\quad D_{2,3}=\{f_{2,[2,3]}|_{\Delta}=0\},
$$
then skipping the hypertangent divisor $D_{2,4}$ and carrying on
in the standard way.

If $d_1=2$, then $f_1=f_{1,2}$, so that we intersect $Y$ with
hypertangent divisors $D_{2,2}$ and $D_{2,3}$, skip $D_{2,4}$ and
complete the construction in the standard way. This proves the
proposition in the case $\mathop{\rm codim}(\Delta\subset F)=2$.

In the case $\mathop{\rm codim}(\Delta\subset F)=3$ the assumption
$$
\frac{\mathop{\rm mult}_o}{\mathop{\rm deg}}\, Y>\frac{8}{d_1d_2}
$$
is so strong that the standard technique of hypertangent divisors
(in the linear system $\Lambda_2|_{\Delta}$ choose one general
divisor, after which argue as in the non-singular case) provides a
contradiction. Note that the condition (R$2^2$.3) consists of two
parts: when $\mathop{\rm codim}({\cal L}\subset{\mathbb
C}^{M+2})=2$ or 3, and in the proof one should use the one where
$$
\mathop{\rm codim\,}({\cal L}\subset{\mathbb C}^{M+2})=\mathop{\rm
codim\,}(\Delta\subset F).
$$
Q.E.D.

{\bf Remark 2.2.} For our purposes it is sufficient that in the
case $\mathop{\rm codim}(\Delta\subset F)=2$ a somewhat weaker
inequality
$$
\frac{\mathop{\rm mult}_o}{\mathop{\rm deg}}\,
Y\leqslant\frac{56}{9d_1d_2}
$$
were true (but for a hyperplane section $\Delta\subset F$ the
inequality
$$
\frac{\mathop{\rm mult}_o}{\mathop{\rm deg}}\,
Y\leqslant\frac{6}{d_1d_2}
$$
remains a requirement). However, such a partial relaxation of the
claim of Proposition 2.3 does not make its proof easier and does
not make the class of varieties covered by the main result
essentially larger (in the sense that it does not improve
essentially the estimate for the codimension of the set of
varieties that do not satisfy the regularity conditions), but
makes its statement more complicated. For that reason Proposition
2.3 is given in the form above.

%%%%%%%%%%%%%%%%%%%%%%%%%%%%%%%%%%%%%%%%%%%%%%%%%%%%%%%%%%%%%%%%%%
%%%%%%%%%%%%%%%%%%%%%%%%%%%%%%%%%%%%%%%%%%%%%%%%%%%%%%%%%%%%%%%%%%
%%%%%%%%%%%%%%%   SECTION 3

\section{Divisorially canonical complete \\ intersections}

In this section we prove Theorem 0.3, that is, the divisorial
canonicity of complete intersections of codimension 2, satisfying
the regularity conditions. Assuming that the claim of the theorem
is not true, we fix a non canonical singularity $E^*$ of a pair
$(F,\frac{1}{n}D)$ and one by one exclude the options when the
centre $B$ of the singularity $E^*$: is not contained in the set
of singular points $\mathop{\rm Sing} F$ (Subsection 3.1); is
contained in $\mathop{\rm Sing} F$, but is not contained in the
set of bi-quadratic singularities of the variety $F$ (Subsection
3.2); finally, is contained in the set of bi-quadratic points
(Subsections 3.3 -- 3.9). The last case is the hardest one and
consumes the largest part of the section. In the proof we use
certain local facts (shown in \S 5) and certain facts of
projective geometry (shown in \S 4). In Subsection 3.10 we briefly
discuss how one of these local facts simplifies a number of
previous papers.\vspace{0.3cm}

{\bf 3.1. Start of the proof. Non-singular points.} Let us start
the proof of Theorem 0.3. Let $F\in{\cal F}_{\rm reg}$ be a
regular complete intersection. Assume that there exists an
effective divisor $D\sim nH$, where $H\in\mathop{\rm Pic}F$ is the
class of a hyperplane section, such that the pair
$(F,\frac{1}{n}D)$ is not canonical, that is, for some divisor
$E^*$ over $F$ the Noether-Fano inequality holds:
$$
\mathop{\rm ord}\nolimits_{E^*}D>n\cdot a(E^*).
$$
The divisor $E^*$ (and other exceptional divisors, which will
emerge in the proof, satisfying the Noether-Fano inequality, or
its log version, $\mathop{\rm
ord}\nolimits_{E^*}D>n\cdot(a(E^*)+1$)), we will also call a non
canonical singularity of the pair $(F,\frac{1}{n}D)$ or a maximal
singularity of the divisor $D$ (respectively, a non log canonical
singularity of the pair $(F,\frac{1}{n}D)$ or a log maximal
singularity of the divisor $D$ in the case of the log version of
the Noether-Fano inequality). By linearity of the Noether-Fano
inequality in $D$ we may assume that $D$ is a prime divisor. Let
$B\subset F$ be the centre of the maximal singularity $E^*$ on
$F$.

In order to prove Theorem 0.3, let us show that the assumption on
the existence of such divisor $D$ (which has a maximal singularity
$E^*$) leads to a contradiction. By the symbol
$$
\mathop{\rm CS}(F,{\textstyle\frac{1}{n}} D)
$$
we denote the union of the centres of all maximal singularities of
the divisor $D$. Assuming $D$ to be fixed, let us assume in
addition that $B$ is an irreducible component of the closed set
$\mathop{\rm CS}(F,\frac{1}{n}D)$; in particular, in a
neighborhood of the generic point of the subvariety $B$ the pair
$(F,\frac{1}{n}D)$ is canonical outside $B$. Obtaining a
contradiction (``exclusion of the maximal singularity'') is
achieved in different ways, depending on whether a point $o\in B$
of general position is non-singular, quadratic or bi-quadratic
point of the complete intersection $F\subset{\mathbb P}$. Let us
show, in the first place, the following fact.

{\bf Proposition 3.1.} {\it The pair $(F,\frac{1}{n}D)$ is
canonical on the set of non-singular points
$F\backslash\mathop{\rm Sing}F$, that is,} $\mathop{\rm
CS}(F,\frac{1}{n}D)\subset\mathop{\rm Sing}F$.

{\bf Proof.} Assume the converse:
$$
B\not\subset\mathop{\rm Sing}F.
$$
In the first place, let us show that $B$ is of sufficiently high
codimension.

{\bf Lemma 3.1.} {\it The following inequality holds:}
$$
\mathop{\rm codim}(B\subset F)\geqslant 8.
$$

{\bf Proof of the lemma.} Assume the converse: $\mathop{\rm
codim}(B\subset F)\leqslant 7$. Then we have
$$
\mathop{\rm codim}(B\subset{\mathbb P})\leqslant 9.
$$
Let $P\subset{\mathbb P}$ be a general 11-dimensional linear
subspace and $F_P=F\cap P$ the corresponding section of the
variety $F$. Since, as we remember, the inequality
$$
\mathop{\rm codim}(\mathop{\rm Sing}F\subset F)\geqslant 10
$$
holds, we have $\mathop{\rm codim}(\mathop{\rm Sing}F\subset
{\mathbb P})\geqslant 12$, so that $F_P$ is a non-singular
complete intersection of codimension 2 in $P\cong{\mathbb
P}^{11}$, and moreover
$$
\mathop{\rm dim}B\cap P\geqslant 2.
$$
By \cite[Proposition 3.6]{Pukh06b} or \cite{Suzuki15}, for every
surface on $F_P$ the multiplicity of the effective divisor
$$
D_P=(D\circ F_P)=D|_{F\cap P}\sim nH_P
$$
(where the symbol $H_P$ means the class of a hyperplane section of
the variety $F_P\subset{\mathbb P}^{11}$) along this surface does
not exceed $n$. This contradicts the inequality
$$
\mathop{\rm mult}\nolimits_{B\cap P}D_P=
\mathop{\rm mult}\nolimits_BD>n,
$$
which is true, because $B\not\subset\mathop{\rm Sing}F$ is the
centre of a maximal singularity of the divisor $D$. Q.E.D. for the
lemma.

Let us come back to the proof of Proposition 3.1. Let $o\in B$ be
a point of general position, so that $o\not\in\mathop{\rm Sing}F$.
Let us consider a general 10-dimensional linear subspace
$P\subset{\mathbb P}$, containing the point $o$. Let $F_P=F\cap P$
be the corresponding section of the variety $F$. Setting
$D_P=(D\circ F_P)$, we get that the pair $(F_P,\frac{1}{n}D_P)$
has the point $o$ as an isolated centre of a non-canonical
singularity, that is, the point $o$ is an irreducible component of
the set $\mathop{\rm CS}(F_P,\frac{1}{n}D_P)$. It is well known
(see, for instance, \cite[Chapter 7, Proposition 2.3]{Pukh13a}),
that this implies that either the inequality
$$
\nu=\mathop{\rm mult}\nolimits_o D>2n
$$
holds, or on the exceptional divisor $E\cong{\mathbb P}^{M-1}$ of
the blow up $\varphi\colon F^+\to F$ of the point $o$ there is a
(uniquely determined) hyperplane $\Theta\subset E$, satisfying the
inequality
$$
\nu+\mathop{\rm mult}\nolimits_{\Theta}D^+>2n,
$$
where $D^+$ is the strict transform of the divisor $D$ on $F^+$.
We will consider only the second option because it formally
includes the first one. By the symbol $|H-\Theta|$ we denote the
(projectively two-dimensional) linear system of hyperplane
sections of $F$, a general element of which is a divisor $W\ni o$,
non-singular at this point and satisfying the equality
$$
W^+\cap E=\Theta
$$
(in other words, $W^+$ cuts out $\Theta$ on $E$). For such
hyperplane section $W$ the restriction $D_W=(D\circ W)$ is an
effective divisor on the factorial variety $W$, satisfying the
inequality
$$
\mathop{\rm mult}\nolimits_oD_W>2n,
$$
which can be re-written in the following way:
$$
\frac{\mathop{\rm mult}_o}{\mathop{\rm deg}}D_W>\frac{2}{d_1d_2}.
$$
The last inequality contradicts Proposition 2.1. This
contradiction completes the proof Proposition 3.1.

Note that we used the claim of Proposition 2.1 not in its full
power: for an effective divisor on a hyperplane section of the
variety $F$, whereas in Proposition 2.1 we had a diisor on the
section of the variety $F$ by a linear subspace of codimension 2.
In its full power the claim of Proposition 2.1 will be used below,
for the exclusion of the bi-quadratic case.\vspace{0.3cm}

%%%%%%%%%%%%%%%%%%%%%%%%%%%%%%%%%%%%%%%%%%%%%%%%%%%%%%%%%%%%%%%
%%%%%%%%%%%%%%%%%%%%%%%   subsection 3.2

{\bf 3.2. Quadratic singular points.} Denote by the symbol
$\mathop{\rm Sing}^{(2\cdot 2)}F$ the closed set of bi-quadratic
points of the variety $F$. Therefore, $\mathop{\rm
Sing}F\backslash\mathop{\rm Sing}^{(2\cdot 2)}F$ is the set of
quadratic singular points.

{\bf Proposition 3.2.} {\it The pair $(F,\frac{1}{n}D)$ is
canonical on the Zariski open set $F\backslash\mathop{\rm
Sing}^{(2\cdot 2)}F$, that is,}
$$
\mathop{\rm CS}(F,{\textstyle \frac{1}{n}}D)\subset
\mathop{\rm Sing}\nolimits^{(2\cdot 2)}F.
$$

{\bf Proof.} Assume the converse:
$$
B\subset\mathop{\rm Sing}F,\quad B\not\subset
\mathop{\rm Sing}\nolimits^{(2\cdot 2)}F.
$$
Let $o\in B$ be a point of general position. It is a quadratic
point of rank $\geqslant 11$ of the variety $F$. Now denote by the
symbol $P$ a general 5-dimensional linear subspace in ${\mathbb
P}$, containing the point $o$. Set $F_P=F\cap P$. The point $o$ is
a non-degenerate (in particular, isolated) quadratic singularity
of the three-dimensional variety $F_P$ (in fact, the only singular
point of that variety). Since $B\subset\mathop{\rm Sing}F$, we
have the inequality $\mathop{\rm codim}(B\subset F)\geqslant 10$,
so that the section $F_P$ can be constructed in two steps: first,
intersect $F$ with a general linear subspace $P'\ni o$ in
${\mathbb P}$ of dimension
$$
\mathop{\rm codim}(B\subset F)+2,
$$
and then with a general 5-dimensional subspace $P\subset P'$,
containing the point $o$. Therefore, we can apply inversion of
adjunction and conclude that the pair $(F_P,\frac{1}{n}D_P)$,
where $D_P=(D\circ F_P)\sim nH_P$ ($H_P$ is the class of a
hyperplane section of the variety $F_P$), is not log canonical,
but canonical out side the point $o$. Let
$$
\varphi_P\colon P^+\to P
$$
be the blow up of the point $o$ with the exceptional divisor
${\mathbb E}_P\cong{\mathbb P}^4$ and $F^+_P\subset P^+$ the
strict transform of the variety $F_P$ (obviously, $F^+_P$ is the
blow up of $F_P$ at the point $o$). Set
$$
E_P=F^+_P\cap{\mathbb E}_P.
$$
This is a non-singular two-dimensional quadric $\cong{\mathbb
P}^1\times{\mathbb P}^1$ in some hyperplane $\langle
E_P\rangle\cong{\mathbb P}^3$ in the projective space ${\mathbb
E}_P$. Obviously,
$$
a(E_P,F_P)=1,
$$
so that writing $D^+_P\sim nH_P-\nu E_P$ (where $D^+_P$ is the
strict transform of the divisor $D_P$ on $F^+_P$), we obtain two
options:

(Q1) $\nu>2n$, so that $E_P$ is a non log canonical singularity of
the pair $(F_P,\frac{1}{n}D_P)$,

(Q2) $n<\nu\leqslant 2n$, and then the closed set
$$
\mathop{\rm LCS}((F_P,{\textstyle \frac{1}{n}}D_P),F^+_P)
$$
--- the union of centres of all non log canonical singularities
of the {\it original} pair $(F_P,\frac{1}{n}D_P)$ on $F^+_P$
--- is a connected closed subset of the non-singular quadric
$E_P\subset\langle E_P\rangle$, which can be either

(Q2.1) a (possibly reducible) connected curve $C_P\subset E_P$, or

(Q2.2) a point $x_P\in E_P$.

Coming back to the original variety $F$, let us consider the blow
ups
$$
\varphi_{\mathbb P}\colon{\mathbb P}^+\to{\mathbb P}\quad
\mbox{and}\quad \varphi\colon F^+\to F
$$
of the point $o$ on ${\mathbb P}$ and $F$, respectively, where
$F^+\subset{\mathbb P}^+$ is the strict transform of $F$. The
exceptional divisors of those blow ups we denote by the symbols
${\mathbb E}$ and $E$, so that ${\mathbb E}\cong{\mathbb
P}^{M+1}$, and $E\subset {\mathbb E}$ is a quadratic hypersurface
of rank $\geqslant 11$ in the hyperplane $\langle
E\rangle\subset{\mathbb E}$. Write down for the strict transform
$$
D^+\sim nH-\nu E,
$$
where $\mathop{\rm mult}_oD=2\nu$. Now Proposition 2.2 implies
that the case (Q1) does not realize. It is not hard to exclude the
case (Q2.2), either: since $P\ni o$ is a general linear subspace,
in the case (Q2.2) we see that on the quadric $E$ there is a
linear subspace $\Lambda$ of codimension 2 (with respect to $E$),
such that $x_P=\Lambda\cap E_P$. However, on the quadric of rank
$\geqslant 11$ there are no linear subspaces of codimension 2 (it
is sufficient to have the inequality $\mathop{\rm rk}E\geqslant
7$). Therefore, the case (Q2.1) takes place. Since
$$
\mathop{\rm Pic}E={\mathbb Z}{\mathbb H},
$$
where ${\mathbb H}$ is the class of a hyperplane section of the
quadric $E$, every irreducible component of the connected curve
$C_P\subset E_P$ is of bi-degree of the form $(l,l)$ with respect
to the presentation $E_P\cong{\mathbb P}^1\times{\mathbb P}^1$,
where $l\geqslant 1$.

Now we will need the following local fact.

Let $o\in {\cal X}$ be a germ of an isolated singularity,
$\mathop{\rm dim}{\cal X}\geqslant 3$ and $\varphi_{\cal
X}\colon{\cal X}^+\to{\cal X}$ the blow up of the point $o$.
Assume that the following condition is satisfied:

$(G)$ the exceptional divisor ${\cal Q}=\varphi^{-1}_{\cal X}(o)$
is irreducible, reduced and non-singular in codimension 1, that
is,
$$
\mathop{\rm codim}(\mathop{\rm Sing}{\cal Q}
\subset{\cal Q})\geqslant 2,
$$
and the equality $a({\cal Q},{\cal X})=1$ holds.

Obviously, $\mathop{\rm Sing}{\cal X}^+\cap{\cal
Q}\subset\mathop{\rm Sing}{\cal Q}$. Assume, furthermore, that the
pair $({\cal X},\frac{1}{n}{\cal D})$ is not log canonical, but
canonical outside the point $o$, and ${\cal E}\neq{\cal Q}$, an
exceptional divisor over ${\cal X}$, is a non log canonical
singularity of that pair, the centre of which on ${\cal X}^+$ is a
prime divisor ${\cal W}\subset{\cal Q}$. Assume, finally, that
$$
{\cal D}^+\sim-\nu_{\cal D}{\cal Q},
$$
where $\nu_{\cal Q}\leqslant 2n$ and the symbol ${\cal D}^+$
stands for the strict transform of the effective divisor ${\cal
D}$ on ${\cal X}^+$. In these assumptions and notations the
following claim is true.

{\bf Theorem 3.1.} (i) {\it The prime divisor ${\cal W}$ comes
into the scheme-theoretic intersection $({\cal D}^+\circ{\cal Q})$
with the multiplicity} $\mu_{\cal W}> n$.

(ii) {\it The following inequality holds:}
$$
\nu_{\cal D}+\mathop{\rm mult}\nolimits_{\cal W}{\cal D}^+>2n.
$$

{\bf Proof} is given in \S 5.

The first claim of Theorem 3.1 implies that the curve $C_P$ is
irreducible and of bi-degree $(1,1)$, that is, it is a plane
section of the quadric $E_P$. This fact and the second claim of
Theorem 3.1 imply that on the quadric $E$ there is a hyperplane
section $\Lambda\subset E$, such that
$$
\nu+\mathop{\rm mult}\nolimits_{\Lambda}D^+>2n.
$$
Consider the pencil $|H-\Lambda|$ of hyperplane sections of the
variety $F$, a general divisor $W\in|H-\Lambda|$ of which contains
the point $o$, and its strict transform $W^+\subset F^+$ cuts out
$\Lambda$ on $E$:
$$
W^+\cap E=\Lambda.
$$
Obviously, the quadric $\Lambda$ is the exceptional divisor of the
blow up
$$
\varphi_ W\colon W^+\to W
$$
of the point $o$ on the variety $W$. Now set $D_W=(D\circ W)$, so
that for the strict transform of that divisor on $W^+$ we have
$$
D^+_W\sim nH_W-\nu_W\Lambda,
$$
where $\nu_W=\nu+\mathop{\rm mult}_{\Lambda}D^+>2n$ and $H_W$ is
the class of a hyperplane section of the variety $W$. Therefore,
for the effective divisor $D_W$ the inequality
$$
\frac{\mathop{\rm mult}_o}{\mathop{\rm deg}}D_W>\frac{4}{d_1d_2}
$$
holds, which contradicts Proposition 2.2. This completes the proof
of Proposition 3.2.\vspace{0.3cm}

%%%%%%%%%%%%%%%%%%%%%%%%%%%%%%%%%%%%%%%%%%%%%%%%%%%%%%%%%%%%%%%%%%%%%%%%%%%%%
%%%%%%%%%%%%%%%%%%%%%%%%%   subsection 3.3

{\bf 3.3. Bi-quadratic points: start of the proof.} It is most
difficult to exclude the bi-quadratic case. Let $B$ be an
irreducible component of the closed set $\mathop{\rm
CS}(F,\frac{1}{n}D)\subset\mathop{\rm Sing}^{(2\cdot 2)}F$, and
$o\in B$ a point of general position. Since we will have to
consider sections of the variety $F\subset{\mathbb P}$ by various
linear subspaces, it is convenient to introduce the following
notations. For $k\geqslant 1$ the symbol $P_k$ stands for a
$k$-dimensional linear subspace in ${\mathbb P}$, containing the
point $o$. This subspace is assumed to be general in its family,
which will be specified every time: either it is the family of all
linear subspaces of dimension $k$ in ${\mathbb P}$, containing the
point $o$, or it is a certain more special family, for instance,
the family of $k$-dimensional subspaces in some hyperplane,
containing the point $o$. In any case this generality will be
sufficient for the intersection
$$
F_k=(F\circ P_k)_{\mathbb P}=F\cap P_k
$$
(the symbol $(\,\circ\,)_{\mathbb P}$ means the scheme-theoretic
intersection in ${\mathbb P}$) to be an irreducible reduced
complete intersection of type $d_1\cdot d_2$ in ${\mathbb P}$.
Furthermore, the blow up of the point $o$ in the space ${\mathbb
P}$ is denoted by the symbol ${\mathbb P}^+$, and its exceptional
divisor by the symbol ${\mathbb E}$; we set $E=(F^+\circ{\mathbb
E})=F^+\cap{\mathbb E}$, where $F^+\subset{\mathbb P}^+$ is the
strict transform of $F$ on ${\mathbb P}^+$, obviously,
$E\subset{\mathbb E}$ is a complete intersection of two quadrics
in ${\mathbb E}\cong{\mathbb P}^{M+1}$. The strict transform of
the subspace $P_k$ on ${\mathbb P}^+$ is denoted by the symbol
$P^+_k$, so that $P^+_k\to P_k$ is the blow up of the point $o\in
P_k$. The exceptional divisor of that blow up is denoted by the
symbol ${\mathbb E}_k$, that is, ${\mathbb E}_k={\mathbb E}\cap
P^+_k\cong{\mathbb P}^{k-1}$. The strict transform of the
subvariety $F_k$ on $P^+_k$ we denote by the symbol $F^+_k$. In
all cases, considered below, the generality of the subspace $P_k$
is sufficient for the intersection
$$
E_k=(F^+_k\circ{\mathbb E}_k)=(F^+_k\circ{\mathbb E})=
F^+_k\cap{\mathbb E}_k
$$
to be an irreducible reduced complete intersection of two quadrics
in ${\mathbb E}_k$. The pull back on $F^+_k$ of a divisorial class
on $F_k$ will be denoted by the same symbol, and the strict
transform on $F^+_k$ of an effective divisor on $F_k$ will be
denoted by adding the upper index $+$. The restriction
$D|_{F_k}=(D\circ F_k)$ of the divisor $D$ onto $F_k$ we denote by
the symbol $D_k$. Therefore, the original pair $(F,\frac{1}{n}D)$
generates the pair $(F_k,\frac{1}{n}D_k)$ on the
$(k-2)$-dimensional variety $F_k$. Since
$$
\mathop{\rm codim}(\mathop{\rm Sing}F\subset F)\geqslant 10
$$
and the condition (R$2^2$.1) holds, for a general subspace
$P_{11}$ we get that the variety $F_{11}$ has a unique singular
point $o$, and moreover, $F^+_{11}$ is a non-singular variety and
$E_{11}\subset{\mathbb E}_{11}\cong{\mathbb P}^{10}$ is a
non-singular intersection of two quadrics. Furthermore, the point
$o$ is a connected component of the set
$$
\mathop{\rm LCS}(F_{11},{\textstyle \frac{1}{n}}D_{11}),
$$
that is, there is a non log canonical singularity of that pair,
the centre of which is the point $o$, whereas the pair
$(F_{11},\frac{1}{n}D_{11})$ is canonical outside the point $o$ in
a neighborhood of that point.

The same is true for the pair $(F_6,\frac{1}{n}D_6)$, where
$P_6\subset P_{11}$ is a general subspace (recall that all
subspaces $P_k$ contain the point $o$ by default). By Proposition
2.3 the inequality $\mathop{\rm mult}_o D\leqslant 6n$ holds, so
that
$$
D^+\sim D-\nu E,
$$
where $\nu\leqslant\frac32 n$, and for that reason $D^+_6\sim
D_6-\nu E_6$. Note that $a(E_6,F_6)=1$, so that $E_6$ is not a non
log canonical singularity of the pair $(F_6,\frac{1}{n}D_6)$. Let
$$
\mathop{\rm LCS}((F_6,{\textstyle \frac{1}{n}}D_6),E_6)
$$
mean the union of the centres of all non log canonical
singularities of the pair $(F_6,\frac{1}{n}D_6)$ over the point
$o$ on the exceptional divisor $E_6$. As we noted above,
$$
\mathop{\rm LCS}((F_6,{\textstyle \frac{1}{n}}D_6),E_6)\neq E_6,
$$
so that $\mathop{\rm LCS}((F_6,\frac{1}{n}D_6),E_6)$ is a proper
connected closed subset of the non-singular three-dimensional
variety $E_6$. This subset can

(B1) contain a surface $S(P_6)\subset E_6$,

(B2) be a connected curve $C(P_6)\subset E_6$,

(B3) be a point $x(P_6)\in E_6$.

Obviously, the case (B3) is impossible: if this case takes place,
then the non-singular 8-dimensional complete intersection of two
quadrics $E_{11}\subset{\mathbb P}^{10}$ would contain a linear
subspace of dimension 5, whereas for the numerical Chow group of
classes of cycles of codimension 3 we have
$$
A^3E_{11}\cong{\mathbb Z}{\mathbb H}^3_{11},
$$
where ${\mathbb H}_{11}$ is the class of a hyperplane section of
the variety $E_{11}\subset{\mathbb P}^{10}$.

Let us show that the case (B1) is impossible, too. Indeed, for the
singularity $o\in F_6$ all assumptions of Theorem 3.1 are
satisfied. Since the non-singular three-dimensional variety $E_6$
is a complete intersection of two quadrics in ${\mathbb
E}_6\cong{\mathbb P}^5$, we have
$$
A^1E_6=
\mathop{\rm Pic}E_6={\mathbb Z}{\mathbb H}_6,
$$
where ${\mathbb H}_6$ is the class of a hyperplane section with
respect to the embedding $E_6\subset{\mathbb E}_6$ (and in the
sequel everywhere the symbol ${\mathbb H}_k$ means the class of a
hyperplane section of the variety $E_k$, embedded in the
projective space ${\mathbb E}_k$). Now we have
$$
(D^+_6\circ E_6)\sim\nu{\mathbb H}_6,
$$
however, by Theorem 3.1 the effective divisor $(D^+_6\circ E_6)$
contains the prime divisor (a surface) $S(P_6)$ with a
multiplicity strictly higher than $n$. Since
$\nu\leqslant\frac32n$, we conclude that the surface $S(P_6)$ is a
hyperplane section of the variety $E_6$, and moreover, by the
claim (ii) os Theorem 3.1 the inequality
$$
\nu+\mathop{\rm mult}\nolimits_{S(P_6)}D^+_6>2n
$$
holds. Since $P_{11}$ and $P_6\subset P_{11}$ are linear subspaces
of general position, containing the point $o$, we can conclude
that there is a hyperplane section $\Lambda$ of the variety $E$,
satisfying the inequality
$$
\nu+\mathop{\rm mult}\nolimits_{\Lambda}D^+>2n.
$$
Now we argue in an almost word for word the same way as when we
excluded the case (Q2.1) in the proof of Proposition 3.2: consider
the hyperplane section $W$ of the variety $F$, containing the
point $o$ and cutting  out $\Lambda$ on $E$,
$$
W^+\cap E=\Lambda.
$$
(In the bi-quadratic case such hyperplane section is unique,
whereas in the quadratic case there is a pencil of such sections.)
Setting $D_W=(D\circ W)$, we obtain an effective divisor $D_W$ on
the hyperplane section $W$ of the variety $F$, satisfying the
inequality
$$
\frac{\mathop{\rm mult}_o}{\mathop{\rm deg}}D_W>\frac{8}{d_1d_2},
$$
which contradicts Proposition 2.3. This completes the exclusion of
the case (B1).

Therefore, we may assume that the case (B2) takes
place.\vspace{0.3cm}

%%%%%%%%%%%%%%%%%%%%%%%%%%%%%%%%%%%%%%%%%%%%%%%%%%%%%%%%%%%%%%%
%%%%%%%%%%%%%%%%%%%%%%%%%%%%   subsection 3.4

{\bf 3.4. Investigating the case (B2).} Recall that the pair
$(F_6,\frac{1}{n}D_6)$ is canonical outside the point $o$ (in a
neighborhood of that point), but has a non log canonical
singularity, the centre of which on $F^+_6$ is an irreducible
curve $C(P_6)$. Let us denote this curve, for compatibility of
notations, by the symbol $\Lambda_6$. Now we need one more local
fact.

Let $o\in{\cal X}$ be a germ of an isolated singularity,
satisfying the condition $(G)$ of Subsection 3.2. Assume that the
pair $({\cal X},\frac{1}{n}{\cal D})$ is not log canonical, but
canonical outside the point $o$ and ${\cal E}\neq{\cal Q}$ is a
non log canonical singularity of that pair, the centre of which on
${\cal X}^+$ is an irreducible subvariety ${\cal W}\subset{\cal
Q}$, which is not contained in $\mathop{\rm Sing}{\cal Q}$, and
$$
\mathop{\rm codim}({\cal W}\subset{\cal Q})\geqslant 2.
$$
Assume that ${\cal D}^+\sim-\nu_{{\cal D}}{\cal Q}$, where
$\nu_{{\cal D}}\leqslant 2n$. In these assumptions the following
claim is true. (The proof is given in \S 5.)

{\bf Theorem 3.2.} {\it At least one of the following two
inequalities holds:}

(1) $\mathop{\rm mult}_{{\cal W}}{\cal D}^+>n$,

(2) $\mathop{\rm mult}_{{\cal W}}({\cal D}^+\circ{\cal Q})>
3n-\nu_{{\cal D}}$.

Note that the scheme-theoretic intersection $({\cal D}^+\circ{\cal
Q})$ is the restriction of the effective divisor ${\cal D}^+$ onto
the exceptional divisor ${\cal Q}$, so that
$$
({\cal D}^+\circ{\cal Q})\sim\nu_{\cal D}H_{\cal Q},
$$
where $H_{\cal Q}=-({\cal Q}\circ{\cal Q})$ is the ``hyperplane
section'' of ${\cal Q}$.

Applying Theorem 3.2 to the pair $(F_6,\frac{1}{n}D_6)$, we obtain
the following fact.

{\bf Proposition 3.3.} {\it The following inequality holds:}
$$
\mathop{\rm mult}\nolimits_{\Lambda_6}D^+_6>n.
$$

{\bf Proof.} Assume the converse:
$$
\mathop{\rm mult}\nolimits_{\Lambda_6}D^+_6\leqslant n.
$$
By Theorem 3.2 we get the inequality
$$
\mathop{\rm mult}\nolimits_{\Lambda_6}(D^+_6\circ E_6)>3n-\nu.
$$
Recall that $F_6$ is the section of $F_{11}$ by a linear subspace
$P_6\subset P_{11}$ of general position, so that there exists an
irreducible subvariety $\Lambda_{11}\subset E_{11}$ of codimension
2, such that $\Lambda_6=\Lambda_{11}\cap P^+_6=\Lambda_{11}\cap
F^+_6$, that is, $\Lambda_6$ is the section of $\Lambda_{11}$ by a
linear 5-dimensional subspace of general position in ${\mathbb
E}_{11}\cong{\mathbb P}^{10}$, and the inequality
$$
\mathop{\rm mult}\nolimits_{\Lambda_{11}}(D^+_{11}\circ E_{11})>
3n-\nu
$$
is satisfied. However, $(D^+_{11}\circ E_{11})\sim\nu{\mathbb
H}_{11}$ is an effective divisor on the non-singular complete
intersection of two quadrics $E_{11}\subset{\mathbb
E}_{11}\cong{\mathbb P}^{10}$, which is cut out on $E_{11}$ by a
hypersurface of degree $\nu$, and $\Lambda_{11}\subset E_{11}$ is
an irreducible subvariety of dimension 6, so that by
\cite[Proposition 3.6]{Pukh06b} we have the inequality
$$
\nu>3n-\nu,
$$
and for that reason $\nu>\frac32n$, which is impossible. This
contradiction completes the proof of the proposition.

Getting back to the original variety $F\subset{\mathbb P}$, we
conclude that on the exceptional divisor $E\subset{\mathbb E}$
there is an irreducible subvariety $\Lambda\subset E$ of
codimension 2, satisfying the inequality
$$
\mathop{\rm mult}\nolimits_{\Lambda}D^+>n.
$$
\vspace{0.3cm}

%%%%%%%%%%%%%%%%%%%%%%%%%%%%%%%%%%%%%%%%%%%%%%%%%%%%%%%%%%%%%%%
%%%%%%%%%%%%%%%%%%%%%   subsection 3.5

{\bf 3.5. The secant variety.} Now we will need one simple but
non-trivial fact of the projective geometry, which we will prove
in \S 4.

Let $Q\subset{\mathbb P}^N$ be an irreducible complete
intersection of two quadrics, where $N\geqslant 16$. More
precisely, $Q=Q_1\cap Q_2$, where $Q_i=\{q_i=0\}\subset{\mathbb
P}^N$ is a quadratic hypersurface. Assume that the following two
conditions of generality are satisfied:

(C1) $\mathop{\rm max}\{\mathop{\rm rk}q_1,\mathop{\rm rk}q_2\}
\geqslant 16$,

(C2) $\mathop{\rm rk}(q_1,q_2)\geqslant 10$.

As we could see in Subsection 1.7, the condition (C2) implies the
inequality
$$
\mathop{\rm codim}(\mathop{\rm Sing}Q\subset Q)\geqslant 7.
$$
Let $X\subset Q$ be an irreducible subvariety of codimension 2.
Let us introduce the following notations. For a pair of points
$p\neq q$ in ${\mathbb P}^N$ the symbol $[p,q]$ means the line
joining these points. In order to simplify some formulas, it is
convenient to set $[p,p]=\emptyset$, so that when we use the
symbol $[p,q]$, we do not need to specify that $p\neq q$.
Furthermore, for a pair of distinct points $p\neq q$ on $Q$ by the
symbol $[p,q]_Q$ we denote the line $[p,q]$, {\it if it is
entirely contained in} $Q$. Otherwise, set $[p,q]=\emptyset$. For
convenience set $[p,p]_Q=\emptyset$. Now let us define the {\it
secant variety} of the subvariety $X$ on $Q$ by the formula
$$
\mathop{\rm Sec}(X\subset Q)=\overline{\bigcup_{p,q\in X}[p,q]_Q}
$$
(the line above means the closure).

The following claim is true.

{\bf Theorem 3.3.} {\it Precisely one of the following three
options takes place:

{\rm (S1)} $\mathop{\rm Sec}(X\subset Q)=Q$,

{\rm (S2)} $\mathop{\rm Sec}(X\subset Q)$ is a hyperplane section
of the variety $Q$, on which $X$ is cut out by a hypersurface of
degree $d_X\geqslant 2$ in ${\mathbb P}^N$,

{\rm (S3)} $X=\mathop{\rm Sec}(X\subset Q)$ is the section of the
variety $Q$ by a linear subspace of codimension 2 in ${\mathbb
P}^N$.

In any case $\mathop{\rm Sec}(X\subset Q)$ is the closure of the
union of all lines $[p,q]_Q$, such that} $[p,q]_Q\cap\mathop{\rm
Sing}Q=\emptyset$.

{\bf Proof} of Theorem 3.3 is given in \S 4.

Now let us consider the secant variety $\mathop{\rm
Sec}(\Lambda\subset E)$. The following claim is almost obvious.

{\bf Proposition 3.4.} {\it The support $|D^+|$ of the effective
divisor $D^+$ contains the closed set} $\mathop{\rm
Sec}(\Lambda\subset E)$.

{\bf Proof.} By the last claim of Theorem 3.3 it is sufficient to
show that for any pair of points
$$
p,q\in\Lambda\backslash\mathop{\rm Sing}E
$$
(in particular, $p,q\not\in\mathop{\rm Sing}F^+$) we have the
inclusion
$$
[p,q]_E\subset D^+.
$$
Here we can assume, of course, that $[p,q]_E\neq\emptyset$. Assume
the convers:
$$
[p,q]_E\not\subset D^+.
$$
Restricting the effective Cartier divisor $D^+$ onto the line
$[p,q]=[p,q]_E$, we get:
$$
([p,q]\cdot D^+)=\nu\geqslant\mathop{\rm mult}\nolimits_pD^++
\mathop{\rm mult}\nolimits_qD^+>2n.
$$
This contradiction completes the proof. Q.E.D.

Proposition 3.4 implies that for the secant variety $\mathop{\rm
Sec}(\Lambda\subset E)$ the option (S1) (in the notations of
Theorem 3.3) does not realize: $\mathop{\rm Sec}(\Lambda\subset
E)\neq E$. It is not hard to exclude the option (S3), too.

{\bf Proposition 3.5.} {\it The subvariety $\Lambda\subset E$ is
not a section of the variety $E$ by a linear subspace of
codimension 2 in} ${\mathbb E}$.

{\bf Proof.} Assume the converse:
$$
\Lambda=\mathop{\rm Sec}(\Lambda\subset E).
$$
Denote by the symbol $|H-\Lambda|$ the pencil of hyperplane
sections of the variety $F$, the general element $W$ in which
contains the point $o$ and its strict transform $W^+$ contains
$\Lambda$. Thus the strict transform $|H-\Lambda|^+$ of the pencil
$|H-\Lambda|$ on $F^+$ cuts out on $E$ a pencil of hyperplane
sections of the variety $E\subset{\mathbb E}$, containing
$\Lambda$. Let $W\in|H-\Lambda|$ be a general element of the
pencil.

For the restriction $D_W=(D\circ W)=D|_W$ of the divisor $D$ onto
the hyperplane section $W$ we have:
$$
D^+_W\sim nH_W-\nu E_W,
$$
where $H_W$ is the class of a hyperplane section of
$W\subset{\mathbb P}^{M+1}$ and $E_W=(E\circ W^+)=E\cap W^+$ is
the exceptional divisor of the blow up of the point $o$ on $W$,
and moreover,
$$
\mathop{\rm mult}\nolimits_{\Lambda}D^+_W=
\mathop{\rm mult}\nolimits_{\Lambda}D^+> n.
$$
The divisor $D_W$ is effective, but can be reducible and contain
as a component the subvariety
$$
F_{\Lambda}=\mathop{\rm Bs}|H-\Lambda|,
$$
which is the uniquely determined section of $F$ by a linear
subspace of codimension 2 in ${\mathbb P}$, containing the point
$o$, such that
$$
E\cap F^+_{\Lambda}=\Lambda.
$$
The subvariety $F_{\Lambda}$ is obviously a hyperplane section of
the variety $W$. Now write down:
$$
D_W=\Xi_W+aF_{\Lambda},
$$
where $a\in{\mathbb Z}_+$ and the effective divisor $\Xi_W$ on $W$
does not contain $F_{\Lambda}$ as a component. Since
$F^+_{\Lambda}\sim H_W-E_W$ and $\mathop{\rm
mult}_{\Lambda}F^+_{\Lambda}=1$, we get
$$
\Xi^+_W\sim (n-a)H_W-(\nu-a)E_W,
$$
where $\mathop{\rm mult}_{\Lambda}\Xi^+_W>n-a$ and $\nu-a>n-a$. In
particular, the divisor $\Xi^+_W$ is non-zero. Now we can
construct the well defined effective divisor
$$
D_{\Lambda}=(\Xi_W\circ F_{\Lambda})
$$
on the variety $F_{\Lambda}\subset{\mathbb P}^M$, where
$D_{\Lambda}\sim(n-a)H_{\Lambda}$ (where $H_{\Lambda}$ is the
class of a hyperplane section of $F_{\Lambda}$) and, moreover,
$$
D^+_{\Lambda}\sim(n-a)H_{\Lambda}-\nu_{\Lambda}\Lambda,
$$
where the coefficient $\nu_{\Lambda}$ satisfies the inequality
$$
\nu_{\Lambda}\geqslant(\nu-a)+
\mathop{\rm mult}\nolimits_{\Lambda}\Xi^+_W>2(n-a).
$$

So we constructed an effective divisor $D_{\Lambda}$ on the
section $F_{\Lambda}$ of the variety $F$ by a linear subspace of
codimension 2 in ${\mathbb P}$, satisfying the inequality
$$
\frac{\mathop{\rm mult}_o}{\mathop{\rm deg}}D_{\Lambda}>
\frac{8}{d_1d_2},
$$
which contradicts Proposition 2.3 and completes the proof of
Proposition 3.5.

We conclude that of the three options listed in Theorem 3.3, the
second one (S2) takes place: $\mathop{\rm Sec}(\Lambda\subset E)$
is a hyperplane section of the variety $E\subset{\mathbb E}$, on
which $\Lambda$ is cut out by a hypersurface of degree
$d_{\Lambda}\geqslant 2$.\vspace{0.3cm}

%%%%%%%%%%%%%%%%%%%%%%%%%%%%%%%%%%%%%%%%%%%%%%%%%%%%%%%%%%%%%%%%%
%%%%%%%%%%%%%%%%%%%%   subsection 3.6

{\bf 3.6. Restriction onto a hyperplane section.} Let $R$ be the
uniquely determined hyperplane section of the variety $F$, such
that $o\in R$ and
$$
E_R=E\cap R^+=(E\circ R^+)=\mathop{\rm Sec}(\Lambda\subset E).
$$
By Proposition 3.4, the support $|D^+|$ of the effective divisor
$D^+$ contains $E_R$. Now let us estimate the multiplicity of the
divisor $D^+$ at a point of general position of the subvariety
$E_R$. In order to do it, we will need one local fact more.

Let $o\in{\cal X}$ be a germ of a three-dimensional isolated
non-degenerate bi-quadratic singularity: if $\varphi_{\cal
X}\colon{\cal X}^+\to{\cal X}$ is the blow up of the point $o$,
then ${\cal X}^+$ and the exceptional surface ${\cal
E}=\varphi^{-1}(o)$ are non-singular, and ${\cal E}\cong{\cal
Q}_1\cap{\cal Q}_2\subset{\mathbb P}^4$ is a del Pezzo surface of
degree 4. Let $L\subset{\cal E}$ be a line and $p\neq q$
--- two distinct points on $L$. Furthermore, let ${\cal D}$ be an
effective divisor on ${\cal X}$ and ${\cal D}^+\sim-\nu_{\cal
D}{\cal E}$ its strict transform on ${\cal X}^+$, where $\nu_{\cal
D}\in{\mathbb Z}_+$. Assume that
$$
\mathop{\rm mult}\nolimits_p{\cal D}^+=
\mathop{\rm mult}\nolimits_q{\cal D}^+=\mu\geqslant 1.
$$

{\bf Theorem 3.4.} {\it The following inequality holds:}
$$
\mathop{\rm mult}\nolimits_L{\cal D}^+\geqslant
\frac13(2\mu-\nu_{\cal D}).
$$

{\bf Proof} is given in \S 5. It is absolutely similar to the
proof of Lemma 4.2 in \cite{Pukh15a} and Lemma 4.5 in
\cite[Subsection 3.7]{Pukh00d}, however we can not simply refer to
those statements, because the normal sheaf of the line $L$ with
respect to ${\cal X}^+$ is not the same as in those lemmas, and
for that reason the arguments of \cite{Pukh00d,Pukh15a} require a
small modification.

Now let $p,q\in\Lambda$ be a general pair of distinct points, such
that the line $[p,q]$ is contained in $E$ and does not intersect
the set $\mathop{\rm Sing}E$. Considering the section $F_5=F\cap
P_5$ of the variety $F$ by a general 5-dimensional linear subspace
$P_5\ni o$, such that
$$
F^+_5\supset[p,q],
$$
and applying Theorem 3.4, we get that the multiplicity of the
divisor $D^+$ along $E_R$ is at least $\frac13(2\mu-\nu)$, where
$\mu=\mathop{\rm mult}_{\Lambda}D^+>n$. For that reason for the
effective divisor
$$
D_R=D|_R=(D\circ R)
$$
(the pair $(F,R)$ is canonical and the divisor $D$ by assumption
is irreducible, so that $D\neq R$) we have
$$
D^+_R\sim nH_R-\nu_RE_R,
$$
where $H_R$ is the class of a hyperplane section of the variety
$R\subset{\mathbb P}^{M+1}$, and for $\nu_R$ the inequality
$$
\nu_R \geqslant\nu+\frac13(2\mu-\nu)=\frac23(\nu+\mu)>\frac43n
$$
holds. Unfortunately, this inequality is insufficient to obtain a
contradiction by means of the estimates of \S 2, as it was done in
the previous cases. However, we made a step forward: for the pair
$(R,\frac{1}{n}D_R)$ we have the inequality $\nu_R>\frac43n$,
which is essentially stronger than the inequality $\nu>n$, which
is satisfied for the original pair $(F,\frac{1}{n}D)$. By
Proposition 2.3 we still have $\nu_R\leqslant\frac32n$. Now in
order to exclude the maximal singularity, the centre of which is
contained in the closed set $\mathop{\rm Sing}^{(2\cdot 2)}F$ of
bi-quadratic points, we will study the singularities of the pair
$(R,\frac{1}{n}D_R)$ at the point $o$ and show that the whole
procedure of exclusion that was carried out above for the original
pair $(F,\frac{1}{n}D)$, can be (and even with some
simplifications, if we use the already known facts) repeated for
the new pair $(R,\frac{1}{n}D_R)$, which will complete the proof
of Theorem 0.3.

The hyperplane in ${\mathbb P}$ that cuts out the subvariety $R$
on $F$ is the linear span of $R$ and for that reason is denoted by
the symbol $\langle R\rangle$. The hyperplane in ${\mathbb E}$
that cuts out the subvariety $E_R$ on $E$ is denoted by the symbol
${\mathbb E}_R$. Obviously,
$$
{\mathbb E}_R=\langle E_R\rangle=
{\langle R\rangle}^+\cap{\mathbb E},
$$
where ${\langle R\rangle}^+\subset{\mathbb P}^+$ is the strict
transform on the blow up of the point $o$.

Recall the well known fact: the dimension of the singular set of a
complete intersection in the projective space jumps by at most 1
when we take a hyperplane section, see \cite{IP,Pukh00a}.
Therefore, the inequality
$$
\mathop{\rm codim}(\mathop{\rm Sing}R\subset R) \geqslant 8
$$
holds so that the section of the variety $R\subset{\langle
R\rangle}\cong{\mathbb P}^{M+1}$ by a general 9-dimensional linear
subspace has no singularities.

Now it is convenient to extend the system of notations, introduced
at the beginning of Subsection 3.3, to the sections of the variety
$R$. For a sufficiently general $k$-dimensional linear subspace
$P_k\ni o$ the intersection $R\cap P_k$ is denoted by the symbol
$R_k$. If $P_k\subset\langle R\rangle$ by construction, then
$R_k=F_k$. In this case, obviously,
$$
\mathbb E_k={\mathbb E}\cap P^+_k={\mathbb E}_R\cap P^+_k.
$$
The strict transform of the subvariety $R_k$ on $P^+_k$ is denoted
by the symbol $R^+_k$, and then we have
$$
E_k=F^+_k\cap{\mathbb E}_k=R^+_k\cap{\mathbb E}_k.
$$
For the restriction of the divisor $D$ onto $F_k$ we have
$$
D_k=D|_{F_k}=D_R|_{R_k}.
$$
We will will always specify the inclusion $P_k\subset\langle
R\rangle$ because the last equalities are very useful. A general
$k$-dimensional linear subspace (in ${\mathbb P}$ or $\langle
R\rangle$), which is not required to contain the point $o$, is
denoted by the symbol $\Pi_k$.\vspace{0.3cm}

%%%%%%%%%%%%%%%%%%%%%%%%%%%%%%%%%%%%%%%%%%%%%%%%%%%%%%%%%%%%%%%%
%%%%%%%%%%%%%%%%%%%%%   subsection 3.7

{\bf 3.7. Non-singular points of the variety $R$.} First of all,
let us show an analog of Proposition 3.1 for the pair
$(R,\frac{1}{n}D_R)$.

{\bf Proposition 3.6.} {\it The pair $(R,\frac{1}{n}D_R)$ is
canonical outside the closed set} $\mathop{\rm Sing}R$.

{\bf Proof} is similar to the proof of Proposition 3.1 and we will
only trace its main steps. Assume the converse:
$$
\mathop{\rm CS}(R,{\textstyle \frac{1}{n}}D_R)
\not\subset\mathop{\rm Sing}R.
$$
Let $B^*\subset R$ be the centre of a maximal singularity $E^*_R$
of the pair $(R,\frac{1}{n}D_R)$, such that
$B^*\not\subset\mathop{\rm Sing}R$ and $B^*$ is an irreducible
component of the closed set $\mathop{\rm CS}(R,\frac{1}{n}D_R)$.
First of all, since the section $R\cap\Pi_9$ of the variety $R$ by
a general 9-dimensional subspace $\Pi_9\subset\langle R\rangle$ is
a non-singular complete intersection of codimension 2 in
$\Pi_9\cong{\mathbb P}^9$, we can argue as in the proof of Lemma
3.1 and conclude that
$$
\mathop{\rm codim}(B^*\subset R)\geqslant 6.
$$
Let $p\in B^*$ be a point of general position; in particular,
$p\not\in\mathop{\rm Sing}R$ (the symbol $o$ continues to denote
our fixed bi-quadratic point). Carrying on as in the proof of
Proposition 3.1, consider a general 8-dimensional linear subspace
$\Pi_8\subset\langle R\rangle$, containing the point $p$. The pair
$$
(R\cap\Pi_8,{\textstyle \frac{1}{n}}D|_{R\cap\Pi_8})
$$
has the point $p$ as an isolated centre of a non-canonical
singularity, so that by \cite[Chapter 7, Proposition 2.3]{Pukh13a}
either $\mathop{\rm mult}_p D_R>2n$, or on the exceptional divisor
$E(p,R)\cong{\mathbb P}^{M-2}$ of the blow up $R^{(p)}\to R$ of
the point $p$ there is a uniquely determined hyperplane
$\Theta(p)\subset E(p,R)$, satisfying the inequality
$$
\mathop{\rm mult}\nolimits_pD_R+
\mathop{\rm mult}\nolimits_{\Theta(p)}D^{(p)}_R>2n,
$$
where $D^{(p)}_R$ is the strict transform of $D_R$ on $R^{(p)}$.
The second option include the first one. By the symbol
$|H_R-\Theta(p)|$ denote the projectively two-dimensional linear
system of hyperplane sections of the variety $R$, a general
element $W$ in which contains the point $p$, is non-singular at
that point and satisfies the equality
$$
W^+\cap E(p,R)=\Theta(p).
$$
For a general divisor $W\in|H_R-\Theta(p)|$ the restriction
$D_R|_W=(D_R\circ W)$ satisfies the inequality $\mathop{\rm
mult}_p (D_R\circ W)>2n$, so that
$$
\frac{\mathop{\rm mult}_o}{\mathop{\rm deg}}(D_R\circ W
>\frac{2}{d_1d_2}.
$$
However, $(D_R\circ W)$ is an effective divisor on the section $W$
of the variety $F$ by a linear subspace of codimension 2,
containing the point $p$, and $W$ is non-singular at that point.
We get a contradiction with the claim of Proposition 2.1. This
completes the proof of Proposition 3.6.

{\bf Corollary 3.1.} {\it For a general 9-dimensional subspace
$\Pi_9\subset\langle R\rangle$ the pair
$$
\left(R\cap\Pi_9,\frac{1}{n}D|_{R\cap\Pi_9}\right)
$$
is canonical.}\vspace{0.3cm}

%%%%%%%%%%%%%%%%%%%%%%%%%%%%%%%%%%%%%%%%%%%%%%%%%%%%%%%%%%%%%%%
%%%%%%%%%%%%%%%%%%%%%%%%%%%%%   subsection 3.8

{\bf 3.8. Singularities of the pair $(R,\frac{1}{n}D_R)$ at the
point $o$.} Let us go back to the study of the bi-quadratic point
$o$. Recall that for a section $F_{11}$ of the variety $F$ by a
general 11-dimensional linear subspace $P_{11}\ni o$ the point $o$
is a connected component of the closed set $\mathop{\rm
LCS}(F_{11},\frac{1}{n}D_{11})$, and
$$
D^+_{11}\sim D_{11}-\nu E_{11}\sim nH_{11}-\nu E_{11}.
$$
Furthermore, because of the generality of the subspace $P_{11}$
$$
\mathop{\rm mult}\nolimits_{[E_{11}\cap E_R]}D^+_{11}=
\mathop{\rm mult}\nolimits_{E_R}D^+,
$$
where $E_{11}\cap E_R$ is a hyperplane section of the variety
$E_{11}\subset{\mathbb E}_{11}\cong{\mathbb P}^{10}$. Consider the
linear subspace
$$
P_{10}=P_{11}\cap\langle R\rangle.
$$
On one hand, $P_{10}$ is a general 10-dimensional linear subspace
in $\langle R\rangle$, containing the point $o$. On the other
hand, $P_{10}$ is a hyperplane in $P_{11}$, so that $F_{10}=F\cap
P_{10}$ is a hyperplane section of the variety $F_{11}$. However,
$P_{10}$ is a specially selected --- generally speaking, uniquely
determined --- hyperplane in $P_{11}$, that is, it is {\it not} a
hyperplane of general position in $P_{11}$, containing the point
$o$. For that reason, by inversion of adjunction we can only claim
that
$$
o\in \mathop{\rm LCS}(F_{10},{\textstyle\frac{1}{n}}D_{10}),
$$
where $D_{10}=(D\circ F_{10})=(D_R\circ F_{10})\sim
nH_{10}-\nu_RE_{10}.$ It is clear that $E_{10}=E_{11}\cap E_R$.
The variety $F_{10}$ is a section of the variety $F_{11}$ by a
certain specially selected hyperplane, therefore $\mathop{\rm
Sing} F_{10}$ can be larger than $P_{10}\cap\mathop{\rm Sing}
F_{11}$, in particular, larger than $\{o\}$. However, we have the
inequality
$$
\mathop{\rm codim}(\mathop{\rm Sing}R\subset
\langle R\rangle)\geqslant 10,
$$
which implies that the singularities of the variety
$F_{10}=R_{10}$ are zero-dimensional; in particular, the point $o$
is an {\it isolated} singularity of that variety.

{\bf Proposition 3.7.} {\it The closed set $\mathop{\rm
CS}(F_{10},\frac{1}{n}D_{10})$ is zero-dimensional. In particular,
the pair $(F_{10},\frac{1}{n}D_{10})$ is canonical outside the
point $o$ in a neighborhood of that point, and the point $o$
itself is the centre of a non log canonical singularity of the
pair} $(F_{10},\frac{1}{n}D_{10})$.

{\bf Proof.} Consider a hyperplane $\Pi_9\subset P_{10}$ of
general position (in particular, not containing the point $o$).
Since $P_{10}\subset\langle R\rangle$ is a general 10-dimensional
subspace, containing the point $o$, the subspace $\Pi_9\subset
P_{10}$ is a general 9-dimensional subspace, without any
restrictions. Applying Corollary 3.1, we obtain the first claim of
the proposition. The remaining two follow from it directly. Q.E.D.
for the proposition.

Now we have what we wanted: the pair $(R,\frac{1}{n}D_R)$ at the
point $o$ satisfies the properties which are completely similar to
the properties of the original pair $(F,\frac{1}{n}D)$ at that
point.

Since for a general subspace $P_6\subset\langle R\rangle$ the pair
$(R_6=F_6,\frac{1}{n}D_6)$ is not log canonical at the point $o$,
but canonical in a punctured neighborhood of that point, and
moreover, $a(E_6,R_6)=1$ and $\nu_R\leqslant\frac32 n$, for the
closed set
$$
\mathop{\rm LCS}((R_6,{\textstyle\frac{1}{n}}D_6),E_6)
$$
we have the three options $(B1)_R$, $(B2)_R$ and $(B3)_R$,
identical to the options $(B1)$, $(B2)$ and $(B3)$ of Subsection
3.3. The case $(B3)_R$ can be excluded in an almost word for word
the same way as the case $(B3)$: assuming that it takes place, we
get that the non-singular 7-dimensional complete intersection of
two quadrics $E_{10}\subset{\mathbb P}^9$ (for a general subspace
$P_{10}\subset\langle R\rangle$) contains a linear subspace of
codimension 3 with respect to $E_{10}$, whereas for the numerical
Chow group we have
$$
A^3E_{10}\cong{\mathbb Z}{\mathbb H}^3_{10},
$$
where ${\mathbb H}_{10}$ is the class of a hyperplane section of
the variety $E_{10}\subset{\mathbb P}^9$, a contradiction. Almost
in the word for word the same way as in Subsection 3.3 the case
$(B1)$ was excluded, we now exclude the case $(B1)_R$: from the
inequality $\nu_R\leqslant\frac32n$ we see that the surface
$S(P_6)$ is a hyperplane section of the variety $E_6$, whence,
applying the claim (ii) of Theorem 3.1 and taking into account
that $P_6\subset P_{10}\subset\langle R\rangle$ are subspaces of
general position, we conclude that there is a hyperplane section
$\Lambda(R)\subset E_R$ of the variety $E_R$, satisfying the
inequality
$$
\nu_R+\mathop{\rm mult}\nolimits_{\Lambda(R)}D^+_R>2n.
$$
Now we consider the hyperplane section $W$ of the variety $R$,
containing the point $o$ and cutting out $\Lambda(R)$ on
$E_R\colon W^+\cap E_R=\Lambda(R)$. As a result, we obtain an
effective divisor $(D\circ W)_R$ on the hyperplane section $W$ of
the variety $R$, satisfying the inequality
$$
\frac{\mathop{\rm mult}_o}{\mathop{\rm deg}}(D\circ W)_R>
\frac{8}{d_1d_2},
$$
which contradicts Proposition 2.3.

Therefore, the cases $(B1)_R$ and $(B3)_R$ are impossible and we
may assume that the case $(B2)_R$ takes place.\vspace{0.3cm}

%%%%%%%%%%%%%%%%%%%%%%%%%%%%%%%%%%%%%%%%%%%%%%%%%%%%%%%%%%%%%%%%%%
%%%%%%%%%%%%%%%%%%%%%   subsection 3.9

{\bf 3.9.  Investigating the case $(B2)_R$.} Our arguments are
similar to the arguments of Subsection 3.4. The pair
$(R_6,\frac{1}{n}D_6)$ is canonical outside the point $o$ in a
neighborhood of that point, but has a non log canonical
singularity, the centre of which on $R^+_6$ is an irreducible
curve $\Lambda_6(R)$.

{\bf Proposition 3.8.} {\it The following inequality holds:}
$$
\mathop{\rm mult}\nolimits_{\Lambda_6(R)}D^+_6>n.
$$

{\bf Proof} is completely similar to the proof of Proposition 3.3
(with the obvious replacements of $\nu$ by $\nu_R$, $F_{11}$ by
$F_{10}=R_{10}$, $P_{11}$ by $P_{10}\subset\langle R\rangle$,
$\Lambda_{11} \subset E_{11}$ by $\Lambda_{10}\subset E_{10}$,
$\Lambda_6$ by $\Lambda_6(R)$ etc.) and we do not repeat this
argument. Q.E.D. for the proposition.

Carrying on as in Subsection 3.4, we ascend back to the variety
$R$ and conclude that on the exceptional divisor
$E_R\subset{\mathbb E}_R$ there is an irreducible subvariety
$\Lambda(R)\subset E_R$ of codimension 2, satisfying the
inequality
$$
\mathop{\rm mult}\nolimits_{\Lambda(R)}D^+_R>n.
$$

The next step of our proof is parallel to the arguments of
Subsection 3.5. Since the rank of a quadratic form drops by at
most 2 when the form is restricted onto a hyperplane, Theorem 3.3
applies to the subvariety $\Lambda(R)$ on the complete
intersection of two quadrics $E_R\subset{\mathbb E}_R$.

{\bf Proposition 3.9.} {\it The support $|D^+_R|$ of the effective
divisor $D^+_R$ contains the closed set}
$$
\mathop{\rm Sec}(\Lambda(R)\subset E_R).
$$

{\bf Proof} is identical to the proof of Proposition 3.4 and we do
not give it. Q.E.D.

Obviously, $\mathop{\rm Sec}(\Lambda(R)\subset E_R)\neq E_R$.

{\bf Proposition 3.10.} {\it The subvariety $\Lambda(R)\subset
E_R$ is not a section of the variety $E_R$ by a linear subspace of
codimension 2 in} ${\mathbb E}_R$.

{\bf Proof} is completely similar to the proof of Proposition 3.5.
However, the notations need to be changed, so we will trace
briefly the arguments. Instead of the pencil $|H-\Lambda|$ (in the
proof of Proposition 3.5) we consider the pencil
$|H_R-\Lambda(R)|$, a general element $W$ in which is a hyperplane
section of the variety $R$, containing the point $o$ and such that
$W^+\supset\Lambda(R)$. We set again
$$
D_W=(D_R\circ W)=(D\circ W)
$$
and get $D^+_W\sim nH_W-\nu_RE_W$, where $H_W$ is the class of a
hyperplane section of $W$ and
$$
E_W=(E_R\circ W^+)=(E\circ W^+)=E\cap W^+,
$$
where $\mu_R=\mathop{\rm mult}_{\Lambda(R)}D^+_W=\mathop{\rm
mult}_{\Lambda(R)}D^+_R>n$. Instead of $F_{\Lambda}$ consider the
subvariety
$$
R_{\Lambda(R)}=\mathop{\rm Bs}|H_R-\Lambda(R)|.
$$
It is the section of the variety $R$ by a linear subspace of
codimension 2 in $\langle R\rangle$, containing the point $o$, and
moreover, $R^+_{\Lambda(R)}\cap E_R=\Lambda(R)$. Obviously, the
subvariety $R_{\Lambda(R)}$ is a hyperplane section of the variety
$W$. Now we write down
$$
D_W=\Xi_W+aR_{\Lambda(R)},
$$
where $a\in{\mathbb Z}_+$ and the effective divisor $\Xi_W$ does
not contain $R_{\Lambda(R)}$ as a component. We get:
$$
\Xi^+_W\sim(n-a)H_W-(\nu_R-a)E_W,
$$
where $\mathop{\rm mult}_{\Lambda(R)}\Xi^+_W>n-a$ and
$\nu_R-a>n-a$, so that
$$
D_{\Lambda(R)}=(\Xi_W\circ R_{\Lambda(R)})
$$
is a well defined effective divisor on the variety
$R_{\Lambda(R)}$, satisfying the inequality
$$
\frac{\mathop{\rm mult}_o}{\mathop{\rm deg}}D_{\Lambda(R)}>
\frac{8}{d_1d_2},
$$
which contradicts Proposition 2.3, which is now used in its full
power, and completes the proof of Proposition 3.10.

Therefore, $\mathop{\rm Sec}(\Lambda(R)\subset E_R)$ is a
hyperplane section of the variety $E_R$, on which $\Lambda(R)$ is
cut out by a hypersurface of degree $d_{\Lambda(R)}\geqslant 2$.
Let $Z$ be the uniquely determined hyperplane section of the
variety $R$, such that $o\in Z$ and
$$
E_Z=Z^+\cap E_R=(Z^+\circ E_R)=
\mathop{\rm Sec}(\Lambda(R)\subset E_R).
$$
Applying Theorem 3.4 and arguing as in Subsection 3.6, for the
effective divisor
$$
D_Z=D|_Z=D_R|_Z=(D_R\circ Z)
$$
we have
$$
D^+_Z\sim nH_Z-\nu_ZE_Z,
$$
where $H_Z$ is the class of a hyperplane section of the variety
$Z\subset{\mathbb P}^M$, and for $\nu_Z$ the inequality
$$
\nu_Z\geqslant\nu_R+\frac13(2\mu_R-\nu_R)=\frac23(\nu_R+\mu_R)>
\frac{14}{9}n>\frac32n
$$
holds (recall that $\mu_R=\mathop{\rm mult}_{\Lambda(R)}D^+_R>n$).
This contradicts Proposition 2.3 (which is again appplied in its
full power) and completes the exclusion of the bi-quadratic case.

Q.E.D. for Theorem 0.3.\vspace{0.3cm}

%%%%%%%%%%%%%%%%%%%%%%%%%%%%%%%%%%%%%%%%%%%%%%%%%%%%%%%%%%%%%%%%%
%%%%%%%%%%%%%%%%%%%%%%%%%%   subsection 3.10

{\bf 3.10. A remark on hypersurfaces of index 1.} We conclude with
a discussion of how Theorem 3.1 simplifies the proofs of a number
of known results. In \cite{Pukh05} the divisorial canonicity of
Zariski general non-singular hypersurfaces $F\subset{\mathbb
P}^{M+1}$ of degree $M+1$ was shown for $M\geqslant 5$ (the symbol
$F$ is free after the proof of Theorem 0.3 is completed).
Moreover, as it was shown, for example in \cite{Pukh15a}, the
codimension of the set of hypersurfaces, violating the conditions
of general position (the regularity conditions) at at least one
non-singular point, is estimated from below by a function,
quadratic in $M$ (more precisely, growing as $\frac12 M^2$).
However, every non-trivial family of hypersurfaces contains
singular hypersurfaces. For that reason in \cite{Pukh09b} a study
of hypersurfaces $F\subset{\mathbb P}^{M+1}$ of degree $M+1$ with
quadratic singularities was started: for $M\geqslant 8$
hypersurfaces with non-degenerate quadratic points were
considered. The main result of \cite{Pukh09b} is the divisorial
canonicity (that is to say, the canonicity of every pair
$(F,\frac{1}{n}D)$, where $D\sim nH_F$ is an effective divisor,
$H_F$ is the class of a hyperplane section of $F$) of the
hypersurface $F$, in the assumption that the non-singular points
satisfy the regularity conditions of \cite{Pukh05}, and every
singular point $o\in F$ is a non-degenerate quadratic singularity
and satisfies the following conditions. Let $(z_1,\dots,z_{M+1}$)
be a system of affine coordinates with the origin at the point $o$
and
$$
f=q_2+\dots+q_{M+1}=0
$$
the equation of the hypersurface $F$ with respect to $z_*$,
decomposed into components, homogeneous in $z_*$. It is assumed
that
\begin{itemize}

\item the sequence $q_2,\dots,q_{M+1}$ is regular and

\item for each $k\in\{2,3,4,5\}$ and any linear subspace
$P$ of codimension 2, containing the point $o$, the set
$$
F\cap P\cap\{q_2=0\}\cap\dots\cap\{q_k=0\}
$$
of codimension $k+1$ with respect to $F$ is irreducible and has
multiplicity $(k+1)!$ at the point $o$.
\end{itemize}

It is very difficult to estimate the codimension of the set of
hypersurfaces, for which such strong and complicated regularity
conditions are violated. Essentially, in \cite{Pukh09b} divisorial
canonicity is shown for hypersurfaces $F$ with a {\it unique}
singular point of general position, which allows to apply this
result only to Fano fibre spaces $V/{\mathbb P}^1$ over the
one-dimensional base. Note that the proof given in \cite{Pukh09b}
is technically very difficult, this applies both to \S 2 (``The
local analysis of a divisor at the quadratic point''), and to \S 3
(``Exclusion of a maximal singularity'').

In particular, the main theorem of \cite{Pukh09b} could not be
applied directly to the study of birational geometry of Fano
hypersurfaces of index 2 in ${\mathbb P}^{M+2}$: for that work,
one needs divisorial canonicity of {\it every} hyperplane section,
but they form a $(M+2)$-dimensional family and so the quadratic
singularities of sections can degenerate considerably. So in
\cite[Subsections 2.3-2.5]{Pukh16a}) the arguments of
\cite{Pukh09b} were modified and the conditions of general
position were relaxed to an extent, such that for a general
hypersurface of index 2 these conditions were satisfied for {\it
every} hyperplane section. However, the proof of divisorial
canonicity in \cite{Pukh16a} was based on the approach used in
\cite{Pukh09b} and so was still very complicated.

The next step in simplifying the proof of divisorial canonicity
was made in \cite{Pukh15a}. Now for the equation
$f=q_2+\dots+q_{M+1}$ it was required that the rank of the
quadratic form $q_2$ is at least 8, the divisor
$$
\{q_3|_{\{q_2=0\}}=0\}
$$
on the quadric $\{q_3=0\}$ is not a sum of 3 (not necessarily
distinct) hyperplane sections from the same pencil and, finally,
for any subspace $\Pi\subset{\mathbb C}^{M+1}_{z_1,\dots,z_{M+1}}$
of codimension $c\in\{0,1,2\}$ the sequence
$$
q_2|_{\Pi},\dots,q_{M-c}|_{\Pi}
$$
is regular, see \cite[\S 3, Subsection 3.3]{Pukh15a}. However, the
exclusion of a maximal singularity, the centre of which is
contained in the set of quadratic points, still was a complicated
technical argument, see \cite[Subsection 4.3]{Pukh15a}. The same
argument was used in \cite{Pukh2018b}, too, where in order to
describe the birational geometry of Fano hypersurfaces of index 2
with quadratic singularities one needs every hyperplane section to
be divisorially canonical.

Theorem 3.1, stated in Subsection 3.2 of the present section and
shown in \S 5, turns the exclusion of a maximal singularity for
the quadratic case into an easy exercise. It is sufficient to
assume that the equation $q_2+\dots+q_{M+1}=0$ of the hypersurface
satisfies two conditions of general position: $\mathop{\rm
rk}q_2\geqslant 7$ and the sequence
$$
q_2|_{\Pi},\dots,q_k|_{\Pi}
$$
is regular for $k=\ulcorner\frac14(3M-1)\urcorner$ for any
hyperplane $\Pi\subset{\mathbb C}^{M+1}_{z_1,\dots,z_{M+1}}$.
Repeating the proof of Proposition 3.2 (with simplifications), we
exclude the maximal singularity, the centre of which is contained
in the set of quadratic singularities of the hypersurface $F$.

%%%%%%%%%%%%%%%%%%%%%%%%%%%%%%%%%%%%%%%%%%%%%%%%%%%%%%%%%%%%%%%%%%%%%
%%%%%%%%%%%%%%%%%%%%%%%%%%%%%%%%%%%%%%%%%%%%%%%%%%%%%%%%%%%%%%%%%%%%%
%%%%%%%%%%%%%%%   SECTION 4

\section{Projective geometry}

In this section we prove Theorem 3.3. In Subsection 4.1 we remind
the classification of secant varieties in the projective space for
closed subsets of codimension 2. In Subsection 4.2 we study the
intersections of irreducible subvarieties on a quadric with a
general linear subspace of maximal dimension on that quadric. On
that basis, in Subsection 4.3 we classify the secant varieties on
a quadric for irreducible subvarieties of codimension 2. After
that in Subsection 4.4 we complete the proof of Theorem
3.3.\vspace{0.3cm}

{\bf 4.1. The secant variety in the projective space.} We use the
notations that were introduced in Subsection 3.5 for the statement
of Theorem 3.3:
$$
X\subset Q=Q_1\cap Q_2\subset{\mathbb P}^N
$$
is an irreducible subvariety of codimension 2 on a complete
intersection of two quadrics, satisfying the conditions (C1) and
(C2). We will prove Theorem 3.3 in three steps: first, we consider
a subvariety of codimension 2 in the projective space (when the
classification of secant varieties is trivial), then we show
certain facts about the intersection of an irreducible subvariety
on a quadric with a linear subspace of maximal dimension (making
the claims of \S 4 in \cite{Pukh2018b} stronger and more precise)
and, finally, we study the secant varieties on a complete
intersection of two quadrics. We start with the subvarieties of
the projective space.

Let $Z\subset{\mathbb P}^m$ be an irreducible closed set of
codimension 2. Set
$$
\mathop{\rm Sec } (Z)=\overline{\bigcup_{(p,q)\in Z\times Z}[p,q]}
$$
(recall that $[p,p]=\emptyset$ by definition, see Subsection 3.5).
Considering the obvious case of an irreducible curve in ${\mathbb
P}^3$, it is easy to see that there are three options:

(1.1) $Z=\mathop{\rm Sec } (Z)$ is a linear subspace of
codimension 2 in ${\mathbb P}^m$,

(1.2) $\mathop{\rm Sec } (Z)$ is a hyperplane in ${\mathbb P}^m$
and $Z$ is a hypersurface of degree $d_Z\geqslant 2$ in that
hyperplane,

(1.3) $\mathop{\rm Sec } (Z)={\mathbb P}^m$.

Moreover, if $\Xi\subset{\mathbb P}^m$ is a closed set of
codimension $\geqslant 4$ in the case (1.1), of codimension
$\geqslant 3$ in the case (1.2) and of codimension $\geqslant 2$
in the case (1.3), then for a general pair of points $(p,q)\in
Z\times Z$ the line $[p,q]$ does not meet $\Xi$, so that
$\mathop{\rm Sec } (Z)$ is the closure of the union of all lines
$[p,q]$ for $(p,q)\in Z\times Z$, such that
$[p,q]\cap\Xi=\emptyset$.

Assume now that subset $Z\subset{\mathbb P}^m$ is reducible:
$$
Z=\bigcup_{i\in I}Z_i,
$$
where $Z_i\subset{\mathbb P}^m$ are irreducible subvarieties of
codimension 2, and $Z_i\neq Z_j$ for $i\neq j$ and $|I|\geqslant
2$. It is easy to see that in this case for the secant variety
$\mathop{\rm Sec } (Z)$ (which is swept out by the lines $[p,q]$,
where the points $p,q$ can lie on the same or different components
of the set $Z$) there are the following options:

(2.1) $\mathop{\rm Sec } (Z)$ is a hyperplane in ${\mathbb P}^m$,
containing all subvarieties $Z_i$,

(2.2) $\mathop{\rm Sec } (Z)$ is a union of several ($\geqslant
2$) hyperplanes in ${\mathbb P}^m$,

(2.3) $\mathop{\rm Sec } (Z)={\mathbb P}^m$.

It is easy to make the option (2.2) more precise, but we will not
need a more precise description. Note that the simplest example of
the case (2.2) is given by 3 lines in ${\mathbb P}^3$, passing
through one point, but not contained in one plane: in that case
the secant variety is a union of three planes. Note also that if
$\Xi\subset{\mathbb P}^m$ is a closed set of codimension
$\geqslant 3$ in the cases (2.1) and (2.2) and codimension
$\geqslant 2$ in the case (2.3), then $\mathop{\rm Sec } (Z)$ is
the closure of the union of all lines $[p,q]$ for $(p,q)\in
Z\times Z$, such that $[p,q]\cap\Xi=\emptyset$.\vspace{0.3cm}

%%%%%%%%%%%%%%%%%%%%%%%%%%%%%%%%%%%%%%%%%%%%%%%%%%%%%%%%%%%%%%%%%%
%%%%%%%%%%%%%%%%%%%%   subsection 4.2

{\bf 4.2. Intersections of subvarieties on a quadric.} Now let
$G\subset{\mathbb P}^m$ be an irreducible quadric. Its vertex
subspace $\mathop{\rm Sing} G\subset{\mathbb P}^m$ is of dimension
$m-\mathop{\rm rk}G$, where $\mathop{\rm rk}G=\mathop{\rm rk}g$ is
the rank of the quadratic form $g$, defining this quadric. If
$\mathop{\rm rk} G=m+1$, then the dimension of $\mathop{\rm
Sing}G=\emptyset$ is $(-1)$. By the symbol $G^{\rm sm}$ we denote
the set $G\backslash\mathop{\rm Sing}G$ of non-singular points.
Every point of the quadric lies on at least one linear subspace
$P\subset G$ of maximal dimension $m-\ulcorner\frac12\mathop{\rm
rk} G\urcorner$. Let ${\cal L}$ be the family of all linear
subspaces $P\subset G$ of maximal dimension.

If $\mathop{\rm rk}G\in 2{\mathbb Z}$, then ${\cal L}={\cal
L}_1\sqcup{\cal L}_2$ consists of two connected components, each
of which is a non-singular irreducible projective variety; if
$\mathop{\rm rk}G\not\in 2{\mathbb Z}$, then ${\cal L}$ is a
non-singular irreducible projective variety.

Let $Y\subset G$ be an irreducible subvariety, which is not
contained entirely in $\mathop{\rm Sing}G$, and $P\in{\cal L}$ a
linear subspace of general position (obviously,
$P\supset\mathop{\rm Sing}G$). The following claim makes
Proposition 4.1 in \cite{Pukh2018b} stronger and more precise.

{\bf Proposition 4.1.} (i) {\it The closed set
$$
\overline{Y\cap(P\backslash\mathop{\rm Sing}G)}
$$
is either empty or each of its irreducible components is of
codimension $\mathop{\rm codim}(Y\subset G)$ in $P$.

{\rm (ii)} Assume that the inequality
$$
\mathop{\rm codim}(Y\subset G)\leqslant \left[\frac12\mathop{\rm
rk}G\right]-2
$$
holds. Then $Y\cap P$ is a non-empty closed set, each irreducible
component of which is of codimension $\mathop{\rm codim}(Y\subset
G)$ in $P$, and the algebraic cycle $(Y\circ P)^{\rm sm}$ of the
scheme-theoretic intersection of $Y$ and $P$ on the non-singular
part $G^{\rm sm}$ contains each irreducible component with
multiplicity 1.

{\rm (iii)} Assume that the inequality
$$
\mathop{\rm codim}(Y\subset G)\leqslant \left[\frac12\mathop{\rm
rk}G\right]-3
$$
holds. Then the closed set $Y\cap P$ is irreducible and the
scheme-theoretic intersection $(Y\circ P)^{\rm sm}$ is reduced.
Moreover, the inclusion
$$
\mathop{\rm Sing}(Y\cap P\cap G^{\rm sm})\subset \mathop{\rm
Sing}Y\cap P\cap G^{\rm sm}
$$
holds, so that outside the set $\mathop{\rm Sing}G$ the subvariety
$Y\cap P$ either is non-singular, or the codimension of its
singular set with respect to $P$ is}
$$
\mathop{\rm codim}(\overline{(\mathop{\rm Sing}
Y\cap G^{\rm sm})}\subset G).
$$

{\bf Proof.} (i) Assume that
$$
Y\cap(P\backslash\mathop{\rm Sing} G)\neq\emptyset.
$$
Consider the projection
$$
\pi_P\colon{\mathbb P}^m\backslash P
\to{\mathbb P}^{\ulcorner\frac12\rm rk G\urcorner-1}
$$
from the subspace $P$. The closures of the fibres of the
projection $\pi_P$ are subspaces $\Lambda\supset P$ of dimension
$\mathop{\rm dim}P+1$. For such subspace we have
$$
G\cap\Lambda=P\cup G(P,\Lambda),
$$
where $G(P,\Lambda)\in{\cal L}$, and if $P\in{\cal L}$ is a
subspace of general position and
$\Lambda=\overline{\pi^{-1}_P(s)}$ for a point $s\in{\mathbb
P}^{\ulcorner\frac12\rm rk G\urcorner-1}$ of general position,
then $G(P,\Lambda)$ is a subspace of general position, too. For
every fibre $\Lambda$ of the projection $\pi_P$ the intersection
$$
P_{\Lambda}=P\cap G(P,\Lambda)
$$
is a hyperplane in $P$ (and in $G(P,\Lambda)$), containing the
vertex space $\mathop{\rm Sing}G$, and it is easy to check that,
varying the point $s\in{\mathbb P}^{\ulcorner\frac12\rm rk
G\urcorner-1}$, we obtain the whole linear system of hyperplanes
in $P$, containing $\mathop{\rm Sing}G$ (and the same is true, of
course, about $G(P,\Lambda)$). Restricting the morphism $\pi_P$
onto the irreducible quasi-projective variety $Y\backslash(Y\cap
P)$, we see that for a point $s$ of general position the closed
affine set
$$
Y\cap[G(P,\Lambda)\backslash P_{\Lambda}]
$$
either is empty, or each of its irreducible components is of
codimension $\mathop{\rm codim}(Y\subset G)$ in $G(P,\Lambda)$.
The subspaces $P$ and $G(P,\Lambda)$ are symmetric objects in our
construction, so that we can conclude that for a general point $s$
each irreducible component of the closed affine set
$$
Y\cap[P\backslash P_{\Lambda}]
$$
is of codimension $\mathop{\rm codim}(Y\subset G)$ in $P$ (because
this set is non-empty by our assumption and the mobility of the
hyperplane $P_{\Lambda}\subset P$, pointed out above). But this
means precisely that each irreducible component of the set
$$
\overline{Y\cap(P\backslash\mathop{\rm Sing}G)}
$$
is of codimension $\mathop{\rm codim}(Y\subset G)$. The claim (i)
is shown.

Let us prove the claim (ii). Since
$$
\mathop{\rm codim}(Y\subset {\mathbb P}^m)=
\mathop{\rm codim}(Y\subset G)+1,
$$
the codimension of each component of the non-empty closed set
$Y\cap P$ with respect to the space $P$ does not exceed the number
$\mathop{\rm codim}(Y\subset G)+1$. The assumption of the part
(ii) implies that
$$
\mathop{\rm dim}P-\mathop{\rm codim}(Y\subset G)-1\geqslant
\mathop{\rm dim}\mathop{\rm Sing}G+1,
$$
so that each irreducible component of the set $Y\cap P$ does not
coincide with $\mathop{\rm Sing}G$ and is not contained in that
set. Applying the claim (i) shown above and taking into account
that (in the notations of the proof of the claim (i) above) each
component of a fibre of general position of the projection
$\pi_P|_{Y\backslash(Y\cap P)}$ is reduced at the generic point,
we obtain the claim (ii).

Let us consider the part (iii). It consists of three claims: on
irreducibility, reducedness and on the singularities. Let us
consider them one by one. In order to prove that the set $Y\cap P$
is irreducible, we argue in the same way as in the proof of
\cite[Proposition 4.1]{Pukh2018b}. Assume the converse: for a
general subspace $P\in{\cal L}$
$$
Y\cap P=\bigcup_{i\in I}Y_i(P),
$$
where $|I|\geqslant 2$ and $Y_i(P)$ are all distinct irreducible
components of the set $Y\cap P$. By the part (i) each of them is
of codimension $\mathop{\rm codim}(Y\subset G)$ with respect to
$P$. The assumption of the part (iii) implies the inequality
$$
\mathop{\rm dim} Y_i(P)\geqslant\mathop{\rm dim}\mathop{\rm Sing}G+3.
$$
Therefore (in the notations of the proof of the part (i)) for a
general fibre $\Lambda=\overline{\pi^{-1}_P(s)}$ the closed set
$Y_i(P)\cap P_{\Lambda}$ is irreducible (since the image of the
set $Y_i(P)$ under the projection from the subspace $\mathop{\rm
Sing}G$ has the dimension $\geqslant 2$) and the equality
$$
\mathop{\rm codim}((Y_i(P)\cap P_{\Lambda})\subset P)= \mathop{\rm
codim}(Y\subset G)+1
$$
holds. Moreover, for a general fibre $\Lambda$ of the projection
$\pi_P$ the irreducible closed sets $Y_i(P)\cap P_{\Lambda}$,
$i\in I$, are all distinct (this is true because the hyperplane
$P_{\Lambda}\subset P$ varies in the linear system of hyperplanes,
containing $\mathop{\rm Sing}G$). Therefore, the irreducible
components $Y_i(P)$ are identified by their hyperplane sections
$Y_i(P)\cap P_{\Lambda}$.

Let us use again the symmetry of the subspaces $P$ and
$G(P,\Lambda)$ in our construction. We get that there is a
one-to-one correspondence between the irreducible components of
the intersection $Y\cap P$ and the irreducible components of the
intersection $Y\cap G(P,\Lambda)$: two components correspond to
each other, if their intersections with $P_{\Lambda}$ are equal.
Therefore, we can write:
$$
Y\cap G(P,\Lambda)=\bigcup_{i\in I}Y_i(G(P,\Lambda)).
$$
But then for each $i\in I$ the set
$$
Y_i=\overline{\bigcup_{s\in U}Y_i(G(P,\overline{\pi^{-1}_P(s)}))},
$$
where the union is taken over points of a non-empty Zariski open
subset $U\subset{\mathbb P}^{\ulcorner\frac12\rm rk
G\urcorner-1}$, is an irreducible component of the original set
$Y$, and for $i\neq j$, $i,j\in I$, we have $Y_i\neq Y_j$. This
contradiction with the irreducibility of the set $Y$ completes the
proof of the first, set-theoretic, claim of the part (iii).

That the scheme-theoretic intersection $(Y\circ P)^{\rm sm}$ is
reduced, follows from the fact that a fibre of general position of
the morphism $\pi_P|_{Y\backslash P}$ is reduced: as in the proof
of the claim (i), we obtain the reducedness of the
scheme-theoretic intersection
$$
(Y\circ[G(P,\Lambda)\backslash P_{\Lambda}]),
$$
whence by the symmetry of the subspaces $P$ and $G(P,\Lambda)$ in
our construction, we get the reducedness of the scheme-theoretic
intersection of $Y$ and $P\backslash P_{\Lambda}$ for a general
hyperplane $P_{\Lambda}\subset P$, containing the subspace
$\mathop{\rm Sing}G$, which implies that the intersection $(Y\circ
P)^{\rm sm}$ is reduced.

Finally, let us show the last claim of the part (iii), about the
singularities of the quasi-projective variety $Y\cap P\cap G^{\rm
sm}$. By Bertini's theorem for a point of general position
$s\in{\mathbb P}^{\ulcorner\frac12{\rm rk} G\urcorner-1}$ we have
the equality
$$
\mathop{\rm Sing}[Y\cap(G(P,\Lambda)\backslash P_{\Lambda})]=
(\mathop{\rm Sing}Y)\cap(G(P,\Lambda)\backslash P_{\Lambda}),
$$
where $\Lambda=\overline{\pi_P^{-1}(s)}$. Using the symmtery of
our construction again, we get
$$
\mathop{\rm Sing}[Y\cap(P\setminus P_{\Lambda})]=
(\mathop{\rm Sing}Y)\cap(P\backslash P_{\Lambda})
$$
for a general hyperplane $P_{\Lambda}\supset\mathop{\rm Sing}G$ in
$P$. From here, taking into account of the part (i), follows the
last claim of the part (iii).

Proof of Proposition 4.1 is complete.\vspace{0.3cm}

%%%%%%%%%%%%%%%%%%%%%%%%%%%%%%%%%%%%%%%%%%%%%%%%%%%%%%%%%%%%%%%%%
%%%%%%%%%%%%%%%%%%%%   subsection 4.3

{\bf 4.3. Secant varieties on a quadric.} Now let
$G\subset{\mathbb P}^m$ be a quadric of rank $\mathop{\rm
rk}G\geqslant 8$. For distinct points $p,q\in G$ the symbol
$[p,q]_G$ means the line $[p,q]$, if it is entirely contained in
$G$; otherwise, we set $[p,q]_G=\emptyset$. As usual,
$[p,p]_G=\emptyset$.

For a closed subset $Y\subset G$ set
$$
\mathop{\rm Sec } (Y\subset G)=
\overline{\bigcup_{(p,q)\in Y\times Y}[p,q]_G}.
$$
For this variety $Y$, we take an irreducible subvariety of
codimension 2. By the symbol
$$
\mathop{\rm Sec}\nolimits^* (Y\subset G)
$$
we denote an irreducible component of the closed set $\mathop{\rm
Sec}(Y\subset G)$, satisfying the following property: for a
general linear subspace $P\in{\cal L}$ of maximal dimension,
$\mathop{\rm Sec}^* (Y\subset G)$ contains an irreducible
component of the secant variety $\mathop{\rm Sec} (Y\cap P)$ (in
the sense of Subsection 4.1).

{\bf Proposition 4.2.} {\it One of the following three options
takes place:

(1) $Y=\mathop{\rm Sec } (Y\cap G)$ is the section of the quadric
$G$ by a linear subspace of codimension 2 in ${\mathbb P}^m$,

(2) $\mathop{\rm Sec } (Y\cap G)$ is a hyperplane section of the
quadric $G$ (a quadric of rank $\geqslant 6$), and $Y$ is cut out
on it by a hypersurface of degree} $d_Y\geqslant 2$,

(3) $\mathop{\rm Sec } (Y\cap G)=G$.

{\bf Proof.} Let $P\subset G$ be a general linear subspace of
maximal dimension $m-\ulcorner\frac12\mathop{\rm rk}G\urcorner$.
By the assumption on the rank of the quadric $G$ we can apply the
claim (ii) of Proposition 4.1: each component of the set $Y\cap P$
is of codimension 2 in $P$, and the scheme-theoretic intersection
of $Y$ and $P$ is reduced at the generic point of each of these
components.

Assume that $Y\cap P$ is a linear subspace of codimension 2. Let
$\Pi\subset P$ be a general two-dimensional plane; in particular,
$\Pi\cap\mathop{\rm Sing}G=\emptyset$. Then $Y$ and $\Pi$
intersect transversally at one point, which is non-singular on
$G$. Consider a general 7-dimensional linear subspace
$R\subset{\mathbb P}^m$, containing the plane $\Pi$. Obviously,
$G_R=G\cap R$ is a non-singular 6-dimensional quadric, so that for
the numerical Chow group of classes of cycles of codimension 2 we
have
$$
A^2G_R={\mathbb Z}H^2_R
$$
(where $H_R$ is the class of a hyperplane section of the quadric
$G_R$), and $Y_R=Y\cap R$ is an irreducible subvariety of
codimension 2 in $G_R$ and
$$
Y_R\sim m_YH^2_R
$$
for some $m_Y\geqslant 1$. In fact
$$
1=(Y_R\cdot\Pi)_{G_R}=m_Y,
$$
so that $\mathop{\rm deg}Y_R=\mathop{\rm deg}Y=2$ and $Y$ is a
quadratic hypersurface in its linear span $\langle
Y\rangle\cong{\mathbb P}^{m-2}$. Therefore,
$$
Y=G\cap\langle Y\rangle
$$
is the section of the quadric $G$ by a linear subspace of
codimension 2 in ${\mathbb P}^m$.

Assume now that $\mathop{\rm Sec} (Y\cap P)$ is a hyperplane in
$P$. Then, if $\mathop{\rm Sec} (Y\subset G)\neq G$, then
$\mathop{\rm Sec}^* (Y\subset G)$ is a prime divisor on $G$, which
is cut out on the factorial quadric $G$ by a hypersurface of
degree $s_Y\geqslant 1$ in ${\mathbb P}^m$. The irreducible
subvariety $\mathop{\rm Sec}^* (Y\subset G)$ of codimension 1
satisfies the condition (iii) of Proposition 4.1 and for that
reason
$$
\mathop{\rm Sec}\nolimits^* (Y\subset G)\cap P
$$
is an irreducible hypersurface in $P$, containing the hyperplane
$\mathop{\rm Sec} (Y\cap P)$. Therefore,
$$
\mathop{\rm Sec}\nolimits^*(Y\subset G)\cap P=
\mathop{\rm Sec} (Y\cap P)
$$
and $s_Y=1$, so that $\mathop{\rm Sec} (Y\subset G)=\mathop{\rm
Sec}^*(Y\subset G)$ is a hyperplane section of the quadric $G$.
This hyperplane section by the inequality $\mathop{\rm rk}
G\geqslant 8$ is factorial, so that $Y$ is cut out on it by a
hypersurface of degree $d_Y\geqslant 1$. However, the case $d_Y=1$
is impossible: in that case $\mathop{\rm Sec}(Y\subset G)=Y$.
Therefore, $d_Y\geqslant 2$.

Assume now, that $\mathop{\rm Sec} (Y\cap P)$ is a union of
$\geqslant 2$ hyperplanes. Arguing as above, we get that if
$\mathop{\rm Sec} (Y\subset G)\neq G$, then the closed set
$$
\mathop{\rm Sec}\nolimits^* (Y\subset G)\cap P
$$
is one of those hyperplanes. Since $Y$ can not be contained in two
{\it distinct} hyperplanes in ${\mathbb P}^m$ (in that case we
would get the equality $\mathop{\rm Sec}^*(Y\subset G)=Y$ again,
contrary to the assumption), we conclude that
$$
\mathop{\rm Sec}(Y\subset G)=G.
$$
Finally, if $\mathop{\rm Sec}(Y\cap P)=P$, then obviously, we get
$\mathop{\rm Sec}(Y\subset G)=G$ again.

This completes the proof of Proposition 4.2. Q.E.D.

{\bf Remark 4.1.} Since
$$
\mathop{\rm codim}(\mathop{\rm Sing}G\subset P)=\mathop{\rm rk}G-
\ulcorner{\textstyle \frac12}\mathop{\rm rk}G\urcorner,
$$
the assumption on the rank of the quadric $G$ implies the
inequality
$$
\mathop{\rm codim}(\mathop{\rm Sing}G\subset P)\geqslant 4,
$$
so that the subvariety $\mathop{\rm Sec}(Y\subset G)$ is swept out
by the secant lines $[p,q]_G$ for pairs $p,q\in Y$, such that
$[p,q]_G\cap\mathop{\rm Sing}G=\emptyset$.\vspace{0.3cm}

%%%%%%%%%%%%%%%%%%%%%%%%%%%%%%%%%%%%%%%%%%%%%%%%%%%%%%%%%%%%%%%%%
%%%%%%%%%%%%%%%%%%   subsection 4.4

{\bf 4.4. The secant varieties on the intersection of two
quadrics.} Let us prove Theorem 3.3. We use the notations of
Subsection 3.5. For certainty, we assume that $\mathop{\rm
rk}q_1\geqslant 16$, where $Q_1=\{q_1=0\}$. By the claim (iii) of
Proposition 4.1 for a general subspace $P\subset Q_1$ of maximal
dimension the intersection $Q\cap P=Q_2\cap P$ is an irreducible
quadratic hypersurface, and the codimension $\mathop{\rm
codim}(\mathop{\rm Sing}(Q\cap P)\subset P)$ is at least
$$
\mathop{\rm min}\{\mathop{\rm codim}(\mathop{\rm Sing}Q_1\subset P),
\mathop{\rm codim}(\mathop{\rm Sing}Q\subset Q_1)\}\geqslant 8,
$$
so that $\mathop{\rm rk}q_2|_P\geqslant 8$. Let $X\subset Q$ be an
irreducible subvariety of codimension 2, then $\mathop{\rm
codim}(X\subset Q_1)=3$ and we can apply the claim (iii) of
Proposition 4.1 again. We get that $X\cap P$ is an irreducible
subvariety of codimension 3 in $P$. It is contained in the quadric
$Q\cap P$ of rank $\geqslant 8$ and is of codimension 2 with
respect to that quadric. Applying Proposition 4.2, we get that one
of the following three options takes place:

(1) $X\cap P=\mathop{\rm Sec}((X\cap P)\subset(Q\cap P))$ is the
section of the quadric $Q\cap P$ by a linear subspace of
codimension 2 in $P$,

(2) $\mathop{\rm Sec}((X\cap P)\subset(Q\cap P))$ is a hyperplane
section of the quadric $Q\cap P$, on which the irreducible
subvariety $X\cap P$ is cut out by a hypersurface of degree
$d_X\geqslant 2$ in $P$,

(3) $\mathop{\rm Sec}((X\cap P)\subset(Q\cap P))=Q\cap P$.

Moreover, by Remark 4.1 and the arguments of Subsection 4.1
$\mathop{\rm Sec}((X\cap P)\subset(Q\cap P))$ is swept out by the
lines $[p,q]_{Q\cap P}$, not intersecting the set $\mathop{\rm
Sing}Q_1\cup\mathop{\rm Sing}Q$.

In the case (3) we have $\mathop{\rm Sing}(X\subset Q)=Q$. Assume
that the case (2) takes place, where $\mathop{\rm Sec}(X\subset
Q)\neq Q$. Again let us introduce the notation $\mathop{\rm
Sec}^*(X\subset Q)$: it is an irreducible component of
$\mathop{\rm Sec}(X\subset Q)$, the intersection of which with a
general subspace $P\subset Q_1$ of maximal dimension contains the
set $\mathop{\rm Sec}((X\cap P)\subset(Q\cap P))$. Obviously,
$\mathop{\rm Sec}^*(X\subset Q)$ is a prime divisor on $Q$, and by
the part (iii) of Proposition 4.1 the scheme-theoretic
intersection of that divisor with $P$ is irreducible, reduced and
contains a hyperplane section of the quadric $Q\cap P$, and for
that reason it is equal to that hyperplane section. The variety
$Q\subset{\mathbb P}^N$ is factorial, so that
$$
\mathop{\rm Sec}\nolimits^*(X\subset Q)\sim s_XH_Q,
$$
where $s_X\geqslant 1$ and $H_Q$ is the class of a hyperplane
section of the variety $Q$. Restricting onto $P$, we get $s_X=1$,
so that
$$
\mathop{\rm Sec}\nolimits^*(X\subset Q)=\mathop{\rm Sec}(X\subset Q)
$$
is a hyperplane section of $Q$, and $X$ is a prime divisor on that
section. However, the singular set of any hyperplane section of
the complete intersection $Q$ is of codimension $\geqslant 5$ with
respect to that section, and that section itself is an irreducible
complete intersection of two quadrics in ${\mathbb P}^{N-1}$, so
that it is a factorial variety. Therefore, $X$ is cut out on
$\mathop{\rm Sec}(X\subset Q)$ by some hypersurface of degree
$d_X\geqslant 2$ in ${\mathbb P}^{N-1}$.

It remains to consider the case (1).

Let $\Pi\subset[Q\cap P]$ be a two-dimensional plane of general
position and $R$ the section of the variety $Q\subset{\mathbb
P}^N$ by a general 7-dimensional linear subspace, containing
$\Pi$. Then $R\subset{\mathbb P}^7$ is a non-singular complete
intersection of two quadrics, so that for the numerical Chow group
of classes of cycles of codimension 2 we have
$$
A^2R={\mathbb Z}H^2_R
$$
and $X_R=X\cap R=(X\circ R)_Q\sim d_{X,R}H^2_R$ for
$d_{X,R}\geqslant 1$. In particular, we have
$$
\mathop{\rm deg}X=\mathop{\rm deg}X_R=4d_{X,R}.
$$
However, $(X_R\cdot\Pi)=d_{X,R}=1$, so that $\mathop{\rm deg}X=4$.

The variety $X$ is of dimension $N-4$, so that it is contained in
a hyperplane in ${\mathbb P}^N$. Therefore, $X$ is a prime divisor
on a hyperplane section of the variety $Q$, and that section, as
we noted above, is a factorial complete intersection of
codimension 2 in ${\mathbb P}^{N-1}$. This already implies that
$X$ is a section of the variety $Q$ by a linear subspace of
codimension 2.

The fact that the secant variety $\mathop{\rm Sec}(X\subset Q)$ is
swept out by the lines $[p,q]_Q$, not intersecting the singular
set $\mathop{\rm Sing}Q$, follows from the remarks above.

Proof of Theorem 3.3 is complete.

%%%%%%%%%%%%%%%%%%%%%%%%%%%%%%%%%%%%%%%%%%%%%%%%%%%%%%%%%%%%%%%%%%%
%%%%%%%%%%%%%%%%%%%%%%%%%%%%%%%%%%%%%%%%%%%%%%%%%%%%%%%%%%%%%%%%%%%
%%%%%%%%%%%%%%%   SECTION 5

\section{The local facts}

In this section we prove the local claims, that were used in \S 3
in the proof of the divisorial canonicity of complete
intersections: Theorem 3.1 (in Subsections 5.1-5.3), Theorem 3.2
(in Subsections 5.4) and Theorem 3.4 (in Subsections
5.5).\vspace{0.3cm}

{\bf 5.1. The oriented graph of the singularity ${\cal E}$.} Let
us start the proof of Theorem 3.1. We use the notations of
Subsection 3.2: $\varphi_{\cal X}\colon{\cal X}^+\to{\cal X}$ is
the blow up of the isolated singularity $o\in{\cal X}$, the
exceptional divisor ${\cal Q}=\varphi^{-1}_{\cal X}(o)$ satisfies
the condition $(G)$ and the pair $({\cal X},\frac{1}{n}{\cal D})$
is not log canonical, but canonical outside the point $o$. This
pair has a log maximal singularity ${\cal E}$, which by the
inequality $\nu_{\cal D}\leqslant 2n$ can not coincide with the
exceptional divisor ${\cal Q}$. By assumption, the centre ${\cal
W}$ of the singularity ${\cal E}$ on ${\cal X}^+$ is a prime
divisor on ${\cal Q}$. The inequality
$$
\mathop{\rm ord}\nolimits_{\cal E}{\cal D}>n\cdot(a_{\cal E}+1)
$$
holds, where $a_{\cal E}=a({\cal E},{\cal X})$ is the discrepancy
of ${\cal E}$ with respect to ${\cal X}$. The first (and main)
claim of Theorem 3.1 states: the subvariety ${\cal W}\subset{\cal
Q}$ comes into the scheme-theoretic intersection $({\cal
D}^+\circ{\cal Q})$ with the multiplicity $\mu_{\cal W}>n$.
Consider the resolution of the singularity ${\cal E}$ in the sense
of \cite[Chapter 2, Section 1]{Pukh13a}: the sequence of blow ups
$$
\varphi_{i,i-1}\colon{\cal X}_i\to{\cal X}_{i-1},
$$
$i=1,\dots,N$, where ${\cal X}_0={\cal X}$, ${\cal X}_1={\cal
X}^+$ and $\varphi_{1,0}=\varphi_{\cal X}$, the birational
morphism $\varphi_{i,i-1}$ blows up the centre $B_{i-1}$ of the
singularity ${\cal E}$ on ${\cal X}_{i-1}$ (so that $B_0=o$ and
$B_1={\cal W}$); the exceptional divisor of the blow up
$\varphi_{i,i-1}$ is denoted by the symbol $E_i$ (there is no
danger to mix up these notations with the notations of \S 1), so
that $E_1={\cal Q}$, and, finally, the exceptional divisor
$E_N\subset{\cal X}_N$ of the last blow up realizes the
singularity ${\cal E}$ (that is, the discrete valuations
$\mathop{\rm ord}_E$ and $\mathop{\rm ord}_{E_N}$ coincide). By
our assumptions, the varieties $B_1,\dots,B_{N-1}$ are
subvarieties of codimension 2 in ${\cal X}_1,\dots,{\cal
X}_{N-1}$, respectively, and
$$
B_i\not\subset\mathop{\rm Sing}{\cal X}_i
$$
for $i=1,\dots,N-1$, so that over a non-empty Zariski subset of
the subvariety $B_{i-1}$ the exceptional divisor $E_i$ is alocally
trivial ${\mathbb P}^1$-bundle for $i=2,\dots,N$. Let
$\Gamma=\Gamma_{\cal E}$ be the oriented graph of the resolution
of the singularity ${\cal E}$ (see \cite[Chapter 2, Section
1]{Pukh13a}): the indices
$$
1,\dots,N
$$
are its vertices and a pair of vertices $i,j$ is joined by an
oriented edge (an {\it arrow}) $i\to j$, if $i>j$ and
$$
B_{i-1}\subset E^{i-1}_j,
$$
where the upper index $a$ means the strict transform on ${\cal
X}_a$, in particular, $E^{i-1}_j$ is the strict transform of the
exceptional divisor $E_j\subset{\cal X}_j$ on ${\cal X}_{i-1}$.
Let, as usual, the symbol $p_{ij}$ mean for $i\neq j$ the number
of paths in the oriented graph $\Gamma$ from the vertex $i$ to the
vertex $j$ (so that $p_{ij}=0$ for $i<j$), and for all
$i=1,\dots,N$ set $p_{ii}=1$. In order to simplify the notations,
we write $p_i$ instead of $p_{Ni}$. Set
$$
\mu_i=\mathop{\rm mult}\nolimits_{B_{i-1}}{\cal D}^{i-1}
$$
for $i=2,\dots,N$ and $\mu_1=\nu_{\cal D}$. Now our assumption that
${\cal E}$ (or $E_N$) is a log maximal singularity of the pair
$({\cal X},\frac{1}{n}{\cal D})$, takes the explicit form of the
{\it log Noether-Fano inequality}
$$
\sum^N_{i=1}p_i\mu_i>n\cdot\left(\sum^N_{i=1}p_i+1\right).
$$
Let $L_i\subset E_i$ be the fibre of the projection $E_i\to
B_{i-1}$ over a general point of the variety $B_{i-1}$.
Intersecting the strict transform ${\cal D}^N$ with the strict
transform $L^N_i$, we see that the multiplicities $\mu_i$ for
$i\geqslant 2$ satisfy the inequalities
\begin{equation}\label{12.12.2019.1}
\mu_i\geqslant\sum_{j\to i}\mu_j.
\end{equation}
Recall that by assumption $\mu_1\leqslant 2n$.

Since the variety ${\cal X}_1$ is non-singular at the general
point of the subvariety ${\cal W}=B_1$ of codimension 2, we get
the inequality
$$
\mu_{\cal W}\geqslant\sum_{i\to 1}\mu_i.
$$
Therefore, the claim (i) of Theorem 3.1 is implied by the
following fact.

{\bf Proposition 5.1.} {\it The following inequality holds:}
$$
\sum_{i\to 1}\mu_i>n.
$$

We will prove Proposition 5.1 in two steps: first, we reduce it to
a certain claim of convex geometry, which, in its turn, will be
reduced to a combinatorial claim about the graph
$\Gamma$.\vspace{0.3cm}

%%%%%%%%%%%%%%%%%%%%%%%%%%%%%%%%%%%%%%%%%%%%%%%%%%%%%%
%%%%%%%%%%%%%%%%%%%%   subsection 5.2

{\bf 5.2. Proof of Proposition 5.1: some convex geometry.} Recall
the two well known properties of the graph $\Gamma$. If $i\to j$
and $l$ is a vertex between $i$ and $j$, that is,
$$
j<l<i,
$$
then $i\to l$. (It follows from the fact that the image of
$B_{i-1}$ on ${\cal X}_{l-1}$ is $B_{l-1}$.) Since the centres of
the blow ups $B_1,\dots,B_{N-1}$ are subvarieties of codimension
2, each of them can be contained in the strict transforms of at
most two previous exceptional divisors, which implies the second
property: from every vertex $i$ emerge at most two arrows, one of
which is $i\to(i-1)$ (for $i\geqslant 2$).

Now let us consider the real space ${\mathbb R}^N$ with
coordinates $t_1,\dots,t_N$. Define the linear functions
$$
\lambda^*_0(t_1,\dots,t_N)=\sum^N_{i=1}p_i\,t_i
$$
and for $i=1,\dots,N$
$$
\lambda_i(t_1,\dots,t_N)=t_i-\sum_{j\to i}t_j
$$
(so that, in particular, $\lambda_N=t_N$). Set
$$
\lambda(t_1,\dots,t_N)=\sum_{i\to 1}t_i=t_2+\dots+t_k,
$$
where the number $k\in\{2,\dots,N\}$ is defined by the condition
that in the graph $\Gamma$ there are arrows
$$
2\to 1,\quad \dots,\quad k\to 1,
$$
and either $k=N$, or $(k+1)\nrightarrow 1$.

Consider the hyperplane
$$
\Pi^*=\{\lambda^*_0(t_1,\dots,t_N)=\sum^N_{i=1}p_i+1\}
$$
and the compact convex set $\Delta^*\subset\Pi^*$, given by the
set of inequalities
$$
\lambda_i\geqslant 0 \quad\mbox{for}\quad
i=1,\dots,N,\quad\mbox{and}\quad t_1\leqslant 2.
$$
Considering the point $\frac{1}{n}(\mu_1,\dots,\mu_N)\in{\mathbb
R}^N$, we see that Proposition 5.1 follows from the inequality
$$
\mathop{\rm min}_{\Delta^*}\lambda\geqslant 1.
$$
Furthermore, set
$$
\lambda_0(t_2,\dots,t_N)=\sum^N_{i=2}p_i\,t_i
$$
and re-write the equation of the hyperplane $\Pi^*$ in the
following form:
$$
\lambda_0(t_2,\dots,t_N)=\sum^N_{i=2}p_i+1+p_1(1-t_1).
$$
From here we conclude that it is sufficient to consider the worst
case $t_1=2$, that is, to show the inequality
$$
\mathop{\rm min}_{\Delta^*\cap\{t_1=2\}}\lambda\geqslant 1.
$$
Therefore we consider the real space ${\mathbb R}^{N-1}$ with the
coordinates $t_2,\dots,t_N$. Taking into account the equality
$$
p_1=\sum_{i\to 1}p_i=p_2+\dots+p_k,
$$
we get that the set $\Delta^*\cap\{t_1=2\}\subset{\mathbb
R}^{N-1}$ is the simplex $\Delta$ in the hyperplane
$$
\Pi=\{\lambda_0(t_2,\dots,t_N)=\sum^N_{i=k+1}p_i+1\}
$$
(if $k=N$, then in the right hand part of the equation we have 1),
given by the system of linear inequalities
\begin{equation}\label{14.12.2019.2}
\lambda_i(t_2,\dots,t_N)\geqslant 0\quad\mbox{for}\quad
i=2,\dots,N.
\end{equation}
Therefore, in order to prove Proposition 5.1, it is sufficient to
show the inequality
\begin{equation}\label{13.12.2019.1}
\mathop{\rm min}_{\Delta}\lambda\geqslant 1.
\end{equation}
It is easy to see that if $k=N$, then
$$
p_2=\dots=p_k=1
$$
and $\lambda_0=\lambda$, so that the needed inequality
(\ref{13.12.2019.1}) holds in a trivial way. So we assume that
$N\geqslant k+1$.

{\bf Lemma 5.1.} {\it The following inequality holds:}
$$
\mathop{\rm min}_{\Delta\cap\{t_N=0\}}\lambda\geqslant 1.
$$

{\bf Proof.} Since the inequality (\ref{13.12.2019.1}) is true in
the trivial case $k=N$, we can use induction on $N$. For the
uppermost vertex $N$ there are two options: (1) there is only one
arrow, $N\to (N-1)$, coming out of $N$,  (2) there are are two
arrows, $N\to(N-1)$ and $N\to l$ for some $l\leqslant N-2$, coming
out of $N$. Let us consider the second one, the first one is
simpler.

Let $\Gamma_1$ and $\Gamma_2$ be the subgraphs of the graph
$\Gamma$ with the vertices $1,\dots,N-1$ and $1,\dots,l$,
respectively. Since the first step of every path from $N$ to $i<N$
comes through either $(N-1$) or $l$ (for $i\leqslant l$), we get
$p_i=p_{N-1,i}+p_{l,i}$. For any point
$$
(b_2,\dots,b_{N-1},0)\in\Delta\cap\{t_N=0\},
$$
taking into account the equality $p_N=p_{N,N}=1$, we get
$$
\lambda_0(b_2,\dots,b_{N-1},0)=\sum^{N-1}_{i=2}p_{N-1,i}b_i+
\sum^l_{i=2}p_{l,i}b_i
=\sum^{N-1}_{i=k+1}p_{N-1,i}+\sum^l_{i=k+1}p_{l,i}+p_N+1=
$$
$$
=\left(\sum^{N-1}_{i=k+1}p_{N-1,i}+1\right)+
\left(\sum^l_{i=k+1}p_{l,i}+1\right)
$$
(if $l\leqslant k$, then the sum over $k+1\leqslant i\leqslant l$
is equal to zero). Therefore, at least one of the following two
inequalities holds:
$$
\sum^{N-1}_{i=2}p_{N-1,i}\,b_i\geqslant\sum^{N-1}_{i=k+1}p_{N-1,i}+1
$$
or
$$
\sum^l_{i=2}p_{l,i}\,b_i\geqslant\sum^l_{i=k+1}p_{l,i}+1.
$$
Multiplying the vector $(b_2,\dots,b_{N-1})$ by
$\frac{1}{1+\varepsilon}$ for some $\varepsilon>0$, we can ensure
that this inequality becomes an equality and apply the induction
hypothesis. If the option (1) takes place, the proof is obvious.
Q.E.D. for the lemma.

The minimum of the function $\lambda$ on the simplex $\Delta$ is
attained at one of its vertices, which are obtained by replacing
the inequality signs by the equality signs in all inequalities
(\ref{14.12.2019.2}) but one. By the lemma shown above it is
sufficient to consider the vertex given by the equations
$\lambda_i=0$ for $i=2,\dots,N-1$. The coordinates
$(a_2,\dots,a_N)$ of that vertex are easy to compute, going ``from
the top to the bottom''. Say,
$\lambda_{N-1}(t_{N-1},t_N)=t_{N-1}-t_N$ and for that reason
$a_{N-1}=a_N$ and in general, the equalities
$$
a_i=p_i\,a_N \quad\mbox{for}\quad i=2,\dots,N-1,
$$
hold, where $a_N$ is computed from the relation
$$
\left(\sum^N_{i=2}p^2_i\right)a_N=\sum^N_{i=k+1}p_i+1.
$$
Therefore, Proposition 5.1 follows now from the combinatorial fact
below.

{\bf Proposition 5.2.} {\it The following inequality holds:}
\begin{equation}\label{14.12.2019.4}
(p_2+\dots+p_k)\left(\sum^N_{i=k+1}p_i+1\right)
\geqslant\sum^N_{i=2}p^2_i.
\end{equation}
\vspace{0.3cm}

%%%%%%%%%%%%%%%%%%%%%%%%%%%%%%%%%%%%%%%%%%%%%%%%%%%%%%%%%%%%%%%%
%%%%%%%%%%%%%%%%%%%%   subsection 5.3

{\bf 5.3. Proof of Proposition 5.2.} We argue by induction on the
number of vertices of $\Gamma$, reducing the claim of the
proposition to a similar claim for the subgraph with the vertices
$a,\dots,N$. As a first example we consider the case $k=2$. Here
the inequality (\ref{14.12.2019.4}) takes the form of the
inequality
$$
p_2\left(\sum^N_{i=3}p_i+1\right)\geqslant\sum^N_{i=2}p^2_i.
$$
Assume that the vertices $3,\dots,l$ are connected with the vertex
2 by arrows, but $(l+1)\nrightarrow 2$. Then
$$
p_2=p_3+\dots+p_l
$$
and by the induction hypothesis
$$
(p_3+\dots+p_l)\left(\sum^N_{i=l+1}p_i+1\right)\geqslant
\sum^N_{i=3}p^2_i.
$$
Adding $p^2_2$ to both parts of the last inequality, we complete
the proof in the case $k=2$.

Assume that $k\geqslant 3$. Furthermore, assume that
$$
(k+1)\nrightarrow(k-1).
$$
In that case there are no arrows from the vertices $i\geqslant
k+1$ to the vertices $j\leqslant k-1$, therefore $p_2=\dots=p_k$.
Assume that there are arrows
$$
(k+1)\to k,\quad\dots,\quad l\to k,
$$
but $(l+1)\nrightarrow k$, so that $p_2=\dots=p_k=p_{k+1}+\dots
+p_l$. The inequality (\ref{14.12.2019.4}) takes the form
$$
(k-1)p_2\left(p_2+\sum^N_{i=l+1}p_i+1\right)\geqslant(k-1)p^2_2
+\sum^N_{i=k+1}p^2_i.
$$
This inequality obviously follows from the estimate
$$
(p_{k+1}+\dots+p_l)\left(\sum^N_{i=l+1}p_i+1\right)\geqslant
\sum^N_{i=k+1}p^2_i,
$$
which is true by the induction hypothesis.

Finally, let us consider the case
$$
(k+1)\to(k-1).
$$
Assume that in $\Gamma$ there are arrows
$$
k\to(k-1),\quad (k+1)\to(k-1),\quad \dots,\quad l\to(k-1),
$$
but $(l+1)\nrightarrow(k-1)$. In that case
$$
p_2=\dots=p_{k-1}=p_k+\dots+p_l
$$
and the left hand side of the inequality (\ref{14.12.2019.4}) can
be transformed in the following way:
$$
((k-3)p_2+p_{k-1}+p_k)\left(\sum^l_{i=k+1}p_i+
\sum^N_{i=l+1}p_i+1\right)=
$$
$$
=((k-3)p_2+p_k)\left(\sum^N_{i=k+1}p_i+1\right)+
p_{k-1}\left(\sum^l_{i=k+1}p_i\right)+
p_{k-1}\left(\sum^N_{i=l+1}p_i+1\right).
$$
Taking into account that by the induction hypothesis
$$
p_{k-1}\left(\sum^N_{i=l+1}p_i+1\right)\geqslant\sum^N_{i=k}p^2_i,
$$
and by Lemma 2.7 in \cite[Chapter 2]{Pukh13a}
$$
\sum^N_{i=k+1}p_i+1\geqslant p_{k-1},
$$
one can estimate the left hand side of the inequality
(\ref{14.12.2019.4}) from below by the expression
$$
((k-3)p_2+p_k)p_{k-1}+(p_{k-1}-p_k)p_{k-1}+\sum^N_{i=k}p^2_i=
$$
$$
=(p_2+\dots+p_{k-1})p_{k-1}+\sum^N_{i=k}p^2_i=\sum^N_{i=2}p^2_i,
$$
which is what we need.

Proof of Proposition 5.2 and, therefore, of Proposition 5.1 and of
the claim (i) of Theorem 3.1 is complete.

{\bf Remark 5.1.} Assume in addition that ${\cal X}^+$ and ${\cal
Q}$ are factorial, and
$$
\mathop{\rm Pic}{\cal Q}={\mathbb Z}H_{\cal Q},
$$
where $H_{\cal Q}$ is the intersection class $-({\cal Q}\circ{\cal
Q})$, considered as a divisor on ${\cal Q}$. In that case the
claim (i) of Theorem 3.1 implies that ${\cal W}\sim H_{\cal Q}$ is
the ``hyperplane section'' of the exceptional divisor ${\cal Q}$
(because $\nu_{\cal D}\leqslant 2n$).

{\bf Proof} of the part (ii) of Theorem 3.1 is carried out in the
word for word the same way the proof of Proposition 9 in \cite[\S
3]{Pukh05} or Proposition 2.4 in \cite[Chapter 7]{Pukh13a}. In the
notations of Subsection 5.1 the claim (ii) of Theorem 3.1 takes
the form of the inequality
$$
\mu_1+\mu_2>2n.
$$
We will not repeat these arguments here. Q.E.D. for Theorem
3.1.\vspace{0.3cm}

%%%%%%%%%%%%%%%%%%%%%%%%%%%%%%%%%%%%%%%%%%%%%%%%%%%%%%%%%%%%%%%%%%
%%%%%%%%%%%%%%%%%%%%%   subsection 5.4

{\bf 5.4. Proof of Theorem 3.2.} Recall that in the assumptions of
Theorem 3.2 the centre ${\cal W}\subset{\cal Q}$ of the log
maximal singularity ${\cal E}$ on ${\cal X}^+$ is an irreducible
subvariety of codimension $\geqslant 2$ (with respect to ${\cal
Q}$), which is not contained in $\mathop{\rm Sing}{\cal Q}$, that
is, ${\cal Q}$ and ${\cal X}^+$ are non-singular at the general
point of ${\cal W}$. The rest of the assumptions are the same as
in Theorem 3.1.

Consider again the resolution of the singularity ${\cal E}$. In
the present case it is convenient to somewhat change the notations
introduced in Subsection 5.1: we set ${\cal X}_0={\cal X}^+$ and
${\cal X}_{-1}={\cal X}$, that is, shift the numbering of the blow
ups $\varphi_{i,i-1}$ and varieties ${\cal X}_i$ by one downwards.
In particular, $E_0={\cal Q}$ and $B_0={\cal W}$. Therefore, the
resolution takes the form of the sequence of blow ups
$$
{\cal X}_{-1}={\cal X}\leftarrow{\cal X}_0={\cal X}^+
\leftarrow{\cal X}_1\leftarrow\dots\leftarrow{\cal X}_N.
$$
Let the blow ups $\varphi_{i,i-1}\colon{\cal X}_i\to{\cal
X}_{i-1}$ correspond to the centres $B_{i-1}\subset{\cal X}_{i-1}$
of codimension $\geqslant 3$ (with respect to ${\cal X}_{i-1}$)
for $i=0,1,\dots,L$, and for $i\geqslant L+1$ we have $\mathop{\rm
codim}(B_{i-1}\subset{\cal X}_{i-1})=2$. Consider again the graph
$\Gamma$ of the resolution of the singularity ${\cal E}$: now its
set of vertices is
$$
0,1,\dots,N.
$$
We use the notations $p_{ij}$ and $p_i=p_{Ni}$ for the number of
paths in the graph $\Gamma$, introduced in Subsection 5.1. The
only difference is that the indices $i$ and $j$ can take the value
0. Define the number $k\geqslant 1$ by the condition
$$
B_{i-1}\subset E^{i-1}_0
$$
for $i=1,\dots,k$ (that is, $i\to 0$), so that
$$
p_0=p_1+\dots+p_k.
$$
Setting $\mu_0=\nu_{\cal D}$ and $\mu_i=\mathop{\rm
mult}_{B_{i-1}}{\cal D}^{i-1}$, we get the following well known
explicit for of the log Noether-Fano inequality:
$$
p_0\mu_0+\sum^N_{i=1}p_i\mu_i>\left(p_0+\sum^N_{i=1}p_i\delta_i+1\right)n,
$$
where $\delta_i=\mathop{\rm codim}(B_{i-1}\subset{\cal
X}_{i-1})-1$.

{\bf Lemma 5.2.} {\it If $\mu_1\leqslant n$, then} $k\geqslant
L+1$.

{\bf Proof.} Assume the converse: $k\leqslant L$. From the log
Noether-Fano inequality we get
$$
\sum^k_{i=1}p_i\mu_i+\sum^L_{i=k+1}p_i\mu_i+\sum^N_{i=L+1}p_i\mu_i>
$$
$$
>\left(\sum^k_{i=1}p_i(\delta_i+1-{\frac{\mu_0}{n}})+
\sum^L_{i=k+1}p_i\delta_i+\sum^N_{i=L+1}p_i+1\right)n.
$$
Since $\mu_1\leqslant n$, the more so $\mu_i\leqslant n$ for all
$i\geqslant 1$, so that inequality, written down above, is
impossible: for $i\in\{1,\dots,k\}$ we have
$$
\delta_i+1-\frac{\mu_0}{n}\geqslant 3-\frac{\mu_0}{n}\geqslant 1,
$$
for $i\in\{k+1,\dots,L\}$ we have $\delta_i\geqslant 2$. This
contradiction completes the proof of the lemma.

Since the inequality $\mu_1>n$ is precisely the inequality (1) of
Theorem 3.2, we will assume, starting from this moment, that
$\mu_1\leqslant n$, so that by the lemma shown above we have
$k\geqslant L+1$. Now the log Noether-Fano inequality can be
re-written in the following form:
\begin{equation}\label{25.12.2019.2}
\begin{array}{c}
\displaystyle \sum^L_{i=1}p_i\mu_i+\sum^k_{i=L+1}p_i\mu_i+
\sum^N_{i=k+1}p_i\mu_i>\\
\\ \displaystyle
>\left(\sum^L_{i=1}p_i(\delta_i+1-\frac{\mu_0}{n})+
\left(2-\frac{\mu_0}{n}\right)
\sum^k_{i=L+1}p_i+\sum^N_{i=k+1}p_i+1\right)n.
\end{array}
\end{equation}
Note that $\delta_{E,i}=\delta_i-1$ for $i=1,\dots,L$ is the
discrepancy of the exceptional divisor $E_i\cap E^i_0=E_i\cap{\cal
Q}^i$ of the blow up
$$
\varphi^E_{i,i-1}\colon E^i_0\to E^{i-1}_0
$$
of the subvariety $B_{i-1}\subset E^{i-1}_0$ with respect to
$E^{i-1}_0$:
$$
\delta_{E,i}=\mathop{\rm codim}(B_{i-1}\subset E^{i-1}_0)-1.
$$

{\bf Lemma 5.3.} {\it From every vertex
$$
a\in\{L+1,\dots,k\}
$$
only one arrow goes out, $a\to(a-1)$. In particular, the subgraph
of the graph $\Gamma$ with the vertices $L+1,\dots,k$ is a chain,
if} $k\geqslant L+2$.

{\bf Proof.} By construction, for those values of $a$ we have
$\mathop{\rm codim}(B_{a-1}\subset{\cal X}_{a-1})=2$. Since
$B_{a-1}\subset E_{a-1}$ (for any $a$) and $B_{a-1}\subset
E^{a-1}_0$ (since $a\leqslant k$), there is the only option:
$$
B_{a-1}= E_{a-1}\cap E^{a-1}_0.
$$
Therefore, $\varphi_{a,a-1}$ blows up a divisor on $E^{a-1}_0$ and
for that reason $\varphi^E_{a,a-1}$ are the identity maps, so that
$$
E^L_0\cong E^{L+1}_0\cong\dots\cong E^k_0.
$$
All the previous blow ups $\varphi^E_{i,i-1}$ with $i\leqslant L$
blow up subvarieties of codimension $\geqslant 2$ on $E^{i-1}_0$,
so that they are non-trivial and
$$
E_L\cap E^L_0\neq(E_i\cap E^i_0)^L
$$
for $i\leqslant L-1$, so that for $a\in\{L+1,\dots,k\}$
$$
B_{a-1}=E_{a-1}\cap E^{a-1}_0\neq(E_i\cap E^i_0)^{a-1}
$$
for $i\leqslant L-1$ and for that reason, since $B_{a-1}\subset
E^{a-1}_0$, we get
$$
B_{a-1}\not\subset E^{a-1}_i,
$$
that is, $a\nrightarrow i$. Q.E.D. for the lemma.

Since for any vertices $i<j<l$ from $l\to i$ it follows that $j\to
i$, the lemma implies that, apart from the arrow $(k+1)\to k$,
from the vertex $(k+1)$ at most one more arrow $(k+1)\to(k-1)$ can
come out. This implies the following claim.

{\bf Proposition 5.3.} (i) {\it If in the graph $\Gamma$ there is
an arrow $b\to a$, where $a\leqslant k$ and $b\geqslant k+1$, then
$a=k-1$.

{\rm (ii)} The following equalities hold:
$$
p_L=\dots=p_{k-1}\quad\mbox{and}\quad p_{k-1}=p_k+ \sum_{k\neq
i\to(k-1)}p_i.
$$

{\rm (iii)} For $i=1,\dots,L$ the following equalities hold:}
$$
p_i=p_L\cdot p_{Li}.
$$

{\bf Proof.} The claim (i) is obvious by Lemma 5.3, (ii) follows
from (i), because every path to the vertex $a\in\{L,\dots,k-1\}$
must go through $(k-1)$. Finally, the claim (iii) follows from the
fact that every path from the vertex $N$ to the vertex
$i\in\{1,\dots,L\}$ goes through the vertex $L$, and for that
reason is a composition of a path from $N$ to $L$ and some path
from $L$ to $i$. Q.E.D. for the proposition.

Let us complete the proof of Theorem 3.2. Taking into account that
by assumption $\mu_i\leqslant n$ for $i\geqslant 1$, from the
inequality (\ref{25.12.2019.2}) we get:
$$
p_L\left(\sum^L_{i=1}p_{Li}\mu_i+\sum^{k-1}_{i=L+1}\mu_i\right)+
p_k\mu_k>
p_L\sum^L_{i=1}p_{Li}(\delta_{E,i}+2-\frac{\mu_0}{n})n,
$$
whence, taking into account the obvious inequality $p_k\leqslant
p_{k-1}=p_L$ we immediately get the estimate
\begin{equation}\label{16.12.2020.1}
\sum^L_{i=1}p_{Li}\mu_i+\sum^k_{i=L+1}\mu_i>
\sum^L_{i=1}p_{Li}\left(\delta_{E,i}+2-\frac{\mu_0}{n}\right)n.
\end{equation}
Setting for $i\leqslant k-1$
$$
\varphi_{k,i}=\varphi_{i+1,i}\circ\dots\circ\varphi_{k,k-1}
\colon{\cal X}_k\to{\cal X}_i
$$
and $\varphi_{k,k}=\mathop{\rm id}_{{\cal X}_k}$, consider the
effective divisor
$$
{\cal D}^k=\varphi^*_{k,0}{\cal D}^+-\sum^k_{i=1}\mu_i\varphi^*_{k,i}E_i.
$$
Since ${\cal D}^+$ does not contain ${\cal Q}=E_0$ as a component,
the restriction of ${\cal D}^k$ onto $E^k_0$ is an effective
divisor
$$
({\cal D}^k\circ E^k_0)=(\varphi^E_{k,0})^*{\cal D}^+_E-
\sum^k_{i=1}\mu_i(\varphi^E_{k,i})^*(E_i\circ E^i_0),
$$
where the meaning of the symbols $\varphi^E_{k,i}$ is obvious and
${\cal D}^+_E=({\cal D}^+\circ E_0)$ is the restriction of ${\cal
D}^+$ onto ${\cal Q}$. In particular, the last exceptional divisor
$(E_L\circ E^L_0)$ is subtracted from $(\varphi^E_{k,0})^*{\cal
D}^+_E$ with the multiplicity
$$
\sum^L_{i=1}p_{Li}\mu_i+\sum^k_{i=L+1}\mu_i.
$$
Now we set for $i\in\{1,\dots,L\}$
$$
\nu_{E,i}=\mathop{\rm mult}\nolimits_{B_{i-1}}{\cal D}^{i-1}_E,
$$
and get that $(\varphi_{L,0}^E)^*{\cal D}^+_E$ (recall that for
$a\in\{L+1,\dots,k\}$ the maps $\varphi^E_{a,a-1}$ are the
identity maps, so that $\varphi^E_{L,0}=\varphi^E_{k,0}$) contains
$(E_L\circ E^L_0)$ with the multiplicity
$$
\sum^L_{i=1}p_{Li}\nu_{E,i}\geqslant
\sum^L_{i=1}p_{Li}\mu_i+\sum^k_{i=L+1}p_{Li}\mu_i.
$$
From the inequality (\ref{16.12.2020.1}) we get:
$$
\sum^L_{i=1}p_{Li}\nu_{E,i}>\sum^L_{i=1}p_{Li}
(\delta_{E,i}n+2n-\mu_0).
$$
Taking into account that $\delta_{E,i}\geqslant 1$ and the
multiplicities $\nu_{E,i}$ do not increase, we conclude that
$$
\nu_{E,1}=\mathop{\rm mult}\nolimits_{\cal W}{\cal D}^+_E>3n-\mu_0.
$$
Since $\mu_0=\nu_{\cal D}$, this completes the proof of Theorem
3.2.\vspace{0.3cm}

%%%%%%%%%%%%%%%%%%%%%%%%%%%%%%%%%%%%%%%%%%%%%%%%%%%%%%%%%%%%%%%%%%%
%%%%%%%%%%%%%%%%%%%%%%%%%%   subsection 5.5

{\bf 5.5. Proof of Theorem 3.4.} Here $o\in{\cal X}$ is a germ of
a three-dimensional non-degenerate bi-quadratic singularity. We
use the notations of Subsection 3.6: ${\cal E}\subset{\mathbb
P}^4$ is a non-singular del Pezzo surface of degree 4,
$L\subset{\cal E}$ is a line, $p\neq q$ are distinct points on
$L$. Since
$$
(L^2)_{\cal E}=-1\quad\mbox{and} \quad ({\cal E}\cdot L)_{{\cal
X}^+}=-1,
$$
we have ${\cal N}_{L\slash{\cal X}^+}\cong{\cal
O}_L(-1)\oplus{\cal O}_L(-1)$, so that for the blow up
$\varphi_L\colon{\cal X}_L\to{\cal X}^+$ of the line $L$ we get:
the exceptional divisor
$$
{\cal E}_L=\varphi^{-1}_L(L)\cong{\mathbb P}^1\times{\mathbb P}^1
$$
is isomorphic to the direct product of $L$ (the first factor) and
the fibre ${\mathbb P}^1$. The classes of the fibres  of the
projections of the surface ${\cal E}_L$ onto the direct factors we
denote by the symbols $s_L$ and $f_L$ (where $f_L$ is the class of
the fibre of the projection $\pi_L\colon{\cal E}_L\to L$), so that
$$
A^1{\cal E}_L=\mathop{\rm Pic}{\cal E}_L=
{\mathbb Z}s_L\oplus{\mathbb Z}f_L.
$$

{\bf Lemma 5.4.} {\it The class of the scheme-theoretic
intersection $({\cal E}_L\circ{\cal E}_L)$ as a divisor on the
surface ${\cal E}_L$ is} $-s_L-f_L$.

{\bf Proof.} Let $\widetilde{{\cal E}}=\varphi^*_L{\cal E}-{\cal
E}_L$ be the strict transform of the surface ${\cal E}$ on ${\cal
X}_L$, $\widetilde{{\cal E}}\cong{\cal E}$. The surfaces
$\widetilde{{\cal E}}$ and ${\cal E}_L$ intersect each other
transversally along a non-singular curve $C$, which is a section
of the projection $\pi_L$. Obviously,
$$
({\cal E}_L\cdot C)_{{\cal X}_L}=(C^2)_{\widetilde{\cal E}}=
(L^2)_{\cal E}=-1.
$$
From the direct decomposition of the normal sheaf ${\cal
N}_{L\slash{\cal X}^+}$ (or from the numerical equalities
$$
(C^2)_{{\cal E}_L}=(\widetilde{{\cal E}}\cdot C)_{{\cal X}_L}=
(\varphi^*_L{\cal E}\cdot C)_{{\cal X}_L}-
({\cal E}_L\cdot C)_{{\cal X}_L}=({\cal E}\cdot L)_{{\cal X}^+}-
({\cal E}_L\cdot C)_{{\cal X}_L}=0)
$$
it follows that $C\sim s_L$ on ${\cal E}_L$. Write down
$$
({\cal E}_L\circ{\cal E}_L)\sim -s_L+xf_L
$$
for some $x\in{\mathbb Z}$. From the equalities given above it
follows that $x=({\cal E}^2_L\cdot\widetilde{{\cal E}})_{{\cal
X}_L}$. However,
$$
({\cal E}^2_L\cdot\widetilde{{\cal E}})_{{\cal X}_L}=
(C^2)_{\widetilde{\cal E}}=-1,
$$
which proves the lemma.

Now let us consider the effective divisor ${\cal D}$ on ${\cal X}$
and its strict transforms ${\cal D}^+\sim-\nu_{{\cal D}}{\cal E}$
on ${\cal X}^+$ and ${\cal D}_L\sim-\nu_{{\cal D}}{\cal
E}-\nu_L{\cal E}_L$ on ${\cal X}_L$, where
$$
\nu_L=\mathop{\rm mult}\nolimits_L{\cal D}^+.
$$
Obviously,
$$
{\cal D}_L|_{{\cal E}_L}\sim\nu_{{\cal D}}f_L+\nu_L(s_L+f_L)=
\nu_Ls_L+(\nu_{\cal D}+\nu_L)f_L.
$$
On the other hand, the strict transform ${\cal D}_L$ has the
multiplicity $(\mu-\nu_L)$ along the fibres $\varphi^{-1}_L(p)$
and $\varphi^{-1}_L(q)$ of the projection $\pi_L$, so that the
1-cycle ${\cal D}_L|_{{\cal E}_L}$ contains these fibres with the
multiplicity $\geqslant(\mu-\nu_L)$. For that reason the
inequality
$$
\nu_{\cal D}+\nu_L\geqslant 2(\mu-\nu_L)
$$
holds (since the pseudo-effective cone $A^1_+{\cal E}_L$ of the
surface ${\cal E}_L\cong{\mathbb P}^1\times{\mathbb P}^1$ is,
obviously, ${\mathbb Z}_+s_L\oplus{\mathbb Z}_+f_L$). From there
we conclude that the inequality
$$
\nu_L\geqslant\frac13(2\mu-\nu_{{\cal D}})
$$
holds.

Proof of Theorem 3.4 is completed.

\begin{flushleft}
Department of Mathematical Sciences,\\
The University of Liverpool
\end{flushleft}

\noindent{\it pukh@liverpool.ac.uk}

\end{document}